\documentclass[12pt]{article}

\usepackage{amsmath,amsthm,amssymb}

\theoremstyle{plain}
\newtheorem{Thm}{Theorem}%[section]
\newtheorem{Prop}{Proposition}[section]
\newtheorem{Cor}{Corollary}
\newtheorem{lem}{Lemma}[section]
\theoremstyle{remark}
\newtheorem{rem}{\indent \sc Remark}[section]

\makeatletter

\@addtoreset{equation}{section}

\def\R{\mathbb{R}}
\def\Z{\mathbb{Z}}

\def\v2{\vskip2mm}
\def\n{\noindent}

\def\sgn{{\rm sgn}}

\def\F{{\cal F}}
\def\a{\alpha}
\def\b{\beta}
\def\e{\varepsilon}
\def\de{\delta}
\def\ga{\gamma}
\def\k{\kappa}
\def\la{\lambda}
\def\fa{\varphi}
\def\th{\theta}

\def\om{\omega}
\def\De{\Delta}
\def\Ga{\Gamma}
\def\La{\Lambda}
\def\Om{\Omega}

\def\pf{{\it Proof.}}
\def\v2{\vskip2mm}
\def\n{\noindent}

\def\0{{\bf 0}}

\def\tst12{{\textstyle \frac12}}
\def\1{{\bf 1}}
\def\sg{{\rm sgn\,}}

\def\n{\noindent}
\def\beq{\begin{eqnarray*}}
\def\eeq{\end{eqnarray*}}

\def\beqn{\begin{equation}}
\def\eeqn{\end{equation}}

\begin{document}
Stochastic processes and their Applications (2019), https://doi.org/10.1016/j.spa.2019.02.006
\vskip40mm
\begin{center}
{\bf \Large Asymptotically stable random walks
of index $1<\alpha<2$ killed on a finite set}  \\
\vskip8mm
{K\^ohei UCHIYAMA} \\
\vskip2mm
Department of Mathematics, Tokyo Institute of Technology \\
Oh-okayama, Meguro Tokyo 152-8551\\
%e-mail: \,uchiyama@math.titech.ac.jp
\end{center}

\vskip10mm

{\it running head}:   random walks  killed on a finite set 

{\it key words}:  one dimensional random walk; 
first passage time; killed at the origin; in a domain of attraction, transition probability, stable process


\vskip20mm

\begin{abstract}
For a random walk on the integer lattice $\mathbb{Z}$ that is  
  attracted to a strictly stable process with index $\alpha\in (1, 2)$  we  obtain the asymptotic form  of the transition probability for the  walk killed when it hits a finite set. The asymptotic forms obtained are valid uniformly in a   natural  range of the space and time variables. The situation  is relatively simple when the  limit stable process has jumps  in both positive and negative  directions;   in the other case when the jumps are one sided  rather interesting matters are involved   and  detailed analyses are necessitated.    
 \end{abstract}
\vskip50mm

\newpage
\section{Introduction}

Let   $S_n=X_1+\cdots+ X_n$ be a  random walk on the integer lattice  $\Z$ started at   $S_0 \equiv 0$, where  the increments $ X_1,  X_2, \ldots$ are  independent and identically distributed  random variables  defined on some probability  space $(\Om, \F, P)$ and taking values in $\Z$. Let $E$  indicate the expectation under $P$  as usual and $X$ be a random variable having the same law as $ X_1$.  
We suppose throughout the paper that 
\v2
1) \,  the walk $S_n$ is   in the  domain of attraction of a strictly stable law of index $1<\a<2$ or, what amounts to the same thing (cf \cite{F}),  if $\phi(\th):= Ee^{i\th X}$, then 
\beqn\label{f_hyp}\lim_{\th \to \pm 0} \frac{1-\phi(\th)}{|\th|^\a L(1/\th)} =  e^{\pm i\pi\ga/2}
\eeqn
where  $L(x)$  is a positive even function on the real line  $\R$ slowly varying at infinity  and  $\ga$ is a real number such that  $|\ga|\leq 2-\a$. 
\v2\n
For simplicity we also suppose  that 
\v2
2) the walk is  {\it strongly aperiodic} in the sense of  Spitzer \cite{S}, namely for any $x\in \Z$, $P[S_n =x] >0$ for all sufficiently large $n$. 

\v2
The essential assumption is of course the condition  1), which  entails $EX =0$ so that  the walk is recurrent (see Appendix (A)  for an equivalent  condition in terms of the tails of distribution function of $X$ and some related facts), the condition 2)  giving rise to little loss of generality  (see Remark \ref{rem2.2} and the comment given in  the end of Section 7). 

Under the assumption  (\ref{f_hyp}) with $1\leq \a\leq 2$ ($ L\equiv 1$ is assumed  if $\a=1$) 
Belkin \cite{B1} shows that the law  of a normalized sum $S_n$  started at zero and conditioned to avoid a given finite set until time $n$ converges to a probability law.
Our result in the present paper  is seen as a local version (including the corresponding \lq conditional local limit theorem' as a part)  of his result   with   the starting position $x$ allowed to be arbitrary and the convergence uniform in  $x$  subject to a certain natural constraint.

Assuming the  condition 2)  in addition and restricting the exponent to $1<\a<2$  we  shall obtain  precise  asymptotic forms of the mass function of the hitting time of the origin and of  the transition probability for the  walk killed when it hits the origin.  
The estimates obtained is  uniform for the space variables within the natural space-time region $x=O(c_n)$, where $c_n$ is a norming sequence associated with the walk. For $\ga<|2-\a|$ when the  limit stable process has jumps  in both positive and negative  directions, the transition probability of the killed walk will be  shown to behave roughly  like that of the limit stable process killed at the origin throughout the region.
  For $\ga=\pm(2-\a)$  the situation differs and  involves interesting matters: we shall  identify the space-time region of the same asymptotic behaviour  as for the limit stable process and observe that   in the most of the remaining part that is unbounded the  transition probability for the stable process is negligibly small as compared with the one for the walk, provided that the walk is neither left- nor right-continuous. Here the asymptotics is described by means  of the potential function, $a(x)$ say, whose property varies in a significant way depending on  the distribution function $F$ of $X$ on its lighter side.     
 We shall   extend the results  to the case when  the walk is killed on hitting a finite set instead of the origin.   The corresponding  results for the walks with finite variance are obtained   by  \cite{U1dm}, \cite{U1dm_f}.

%The left tail of $F$ gets lighter,  the drift becomes stronger
%then the probability that  the walk started at $x$ avoids to visit the origin   until a given time $n$ tends to zero as $n\to\infty$ in  a certain  order of magnitude  that is  common to all  $x>0$. We shall  present  the phenomenas which reflect this fact   explicitly.

 Among others the following  two  steps  are crucial in  the proof  of our main results. First we 
 extend  Belkin's result on the conditional law  to the case when  the starting position is allowed to tend to infinity along with  $n$ in a certain reasonable  way. Second   we deduce a local limit theorem from the  integral theorem by Belkin with the help of  Gnedenko's local limit theorem,  in which the deduction rests on    the  idea  devised by Denisov and Wachtel \cite{DW}.
   
 In the  classical papers \cite{KS} and \cite{K}  Kesten and Spitzer, studying a  multidimensional lattice walk  killed on a finite set,    obtained various ratio limit theorems under a mild  assumption on the walk; Kesten \cite{K} especially obtained   an exact asymptotic result for the ratio of the transition probability to a sum of the mass functions of the hitting time which is directly related to our results  (cf.  Remark  \ref{rem4C}  at the end of the next section).    
 For a non-lattice walk on $\R^d$ Port and Stone \cite{PS} established the existence of the potential operator  and applied the result to  study  the ratio limit theorems similar to those treated  in \cite{KS}. 
  \cite{Kai0}  and \cite{Kai} studied the one dimensional  random walk having zero mean and finite variance and killed at zero: in   \cite{Kai0}  the local limit theorem corresponding to the integral version of \cite{B1} is established, while  \cite{Kai} shows that a normalized process of  the walk conditioned to be killed at time  $n$ converges in law  to a Brownian excursion. \cite{B2} obtained a functional limit theorem for   the conditioned walk (on $\Z$) in the same setting as in  \cite{B1}. In  \cite{Vys}   the largest gap within the range of  the walk
conditioned to avoid a bounded set was studied.    Stable  or L\'evy processes  conditioned to avoid zero were studied by  \cite{DKW}, \cite{KKPW} and \cite{Pt}.   Our problem is studied by the present author in an article  \cite{Uattrc} under the restrictive condition that $L\equiv 1$ in 1) by  which  Fourier analysis  is effectively applied;  the results obtained there have provided a guide line for the present work especially in case $\ga=|2-\a|$.

\section{Statements of results}

 We  first introduce fundamental objects  that appear in the description of our results and state some well known facts concerning them. 
Put $p^n(x)=P[S_n=x]$, $p(x)=p^1(x)$ ($x\in \Z$) and define the potential function 
 $$a(x)=\sum_{n=0}^\infty[p^n(0)-p^n(-x)];$$
 the series on the RHS is convergent  and   $a(x)/|x|\to 1/\sigma^2$ and $a(x+y)-a(x) \to \pm y/\sigma^2$ as $x\to \pm \infty$, where $\sigma^2 =E X^2$ (that is infinite under the present setting): cf. Sections 28 and 29 of Spitzer \cite{S}.  
To make  expressions concise we  use the notation
$$a^\dagger(x) = \1(x=0) + a(x),$$
where $\1({\cal S})$ equals  1 or 0 according as  a statement ${\cal S}$ is true or false.
If $S$ is {\it left-continuous} (i.e., $p(x)=0$ for $x\leq2$), then $a(x)= 0$ for all $x\geq 0$ (under $\sigma^2=\infty$), and similarly for  right-continuous walks, whereas if it is neither left- nor right-continuous,
 $a(x)>0$ whenever  $x\neq 0$.   (See Appendices (B) and (C)  for additional facts related to $a$.)

 We write  $S^x_n$ for  $x+S_n$, the walk started at $S^x_0=x\in \Z$. For a subset $B\subset \R$,   put  
 $$\sigma^x_{B}=\inf \{n\geq 1: S^x_n\in B\},$$
  the time of the first entrance of  the walk $S^x$ into $B$. 
To avoid the overburdening of notation we write $S^x_{\sigma_B}$    for  $S^x_{\sigma^x_B}$ and 
 $S_{\sigma_B}$ for $S_{\sigma^0_B}$.
 % sometimes $\sigma B$ is written for $\sigma_B$, e.g., $S^x_{\sigma [0,\infty)}$ for $S^x_{\sigma_{[0,\infty)}}$. 
 
 When the  spatial variables  become indefinitely large the asymptotic results are naturally  expressed by means of the   stable process appearing in the scaling limit and we need to introduce relevant quantities.   
Let $Y_t$ be a stable process started at zero with characteristic exponent
$$\psi(\th) =  e^{ i (\sg \th) \pi\ga/2}|\th|^\a  \qquad (|\ga| \leq 2-\a, \ga \;\mbox{is real})$$
so that $Ee^{i\th Y_t} = e^{- t\psi(\th)}$, where  $\sg \th =1$ if $\th >0$, $0$ if $\th=0$ and  $-1$ if $\th <0$.  ($\ga$ has the same sign as the skewness parameter so that  the extremal case $\ga=2-\a$ corresponds to the spectrally  positive case.)    In the assumption (\ref{f_hyp}) we can suppose the function  $L$---of which  only asymptotics at infinity is significant---so aptly
 chosen as to be differentiable and satisfy  $L'(x)/L(x) = o(1/|x|)$ as $|x|\to\infty$, and then take   positive numbers $c_n$, $n\geq 0$  which are increasing and satisfy  
%Given  a positive number $c_\circ$, we can choose a sequence of positive numbers  $B_n$  so that as $n\to\infty$
 \beqn\label{NC}
 n c_n^{-\a}L(c_n) \to 1. 
 \eeqn
It then follows that 
 $$\lim_{n\to\infty} \log E e^{i\th S_n/c_n} = - \psi(\th),$$
in other words, the law of  $S_n/c_n$ converges to the stable law whose
 characteristic function equals $e^{- \psi(\th)}$.
Denote by 
 $\mathfrak{p}_t(x)$ and  $\mathfrak{f}^{\, x}(t)$  the  density of the distribution of $Y_t$ and of the  first hitting time to the origin of $Y^x_t := x+Y_t$, respectively:
 $$ \mathfrak{p}_t(x)= P[Y_t\in dx]/dx,  \qquad  \mathfrak{f}^{\,x}(t) = (d/dt)P[ \exists s\leq t,   Y^x_s =0];$$
 there exist  jointly continuous versions of these densities  (for $t>0$) and we shall always choose such ones. 
  It follows that $S_{\lfloor nt\rfloor}/c_n \Rightarrow Y_{ t}$ (weak convergence of distribution) and  by Gnedenko's local limit theorem \cite{GK}  as $n\to\infty$
 \beqn\label{llt} 
 p^n(x) = \frac{ \mathfrak{p}_{1} (x/c_n)  +o(1)}{c_n},
 \eeqn 
where $o(1)$ is uniform for $x\in \Z$ and  $\lfloor b\rfloor$ denotes the integer part of a real number $b$.

 For real numbers $s, t$, 
 $s\vee t= \max\{s,t\}$ and  $s\wedge t=\min\{s,t\}$, $t_+=t\vee 0$, $t_- = (-t)_+$  and  $\lceil t\rceil$ denotes the smallest integer that is not less than  $t$;  for positive sequences $(s_n)$ and $(t_n)$, $s_n\sim t_n$ and $s_n\asymp t_n$  mean, respectively,  that the ratio  $s_n/t_n$ approaches unity and    that  $s_n/t_n$ is bounded away from zero and infinity. There will  arise an expression like $s_n\sim Ct_n$ with $C= 0$ which we regard as meaning  $s_n/t_n \to 0$. 
 For two real functions $f(x)$ and $g(y)$  with  $g(y) > 0$ each of the expressions 
 $$f(x) <\!<g(y) \quad \mbox{and}\quad -g(y) <\!< -f(x)$$
means that either $f(x) \leq 0$ or  $[f(x)]_+/g(y)$ tends to zero as  $x$ or  $y$ (or both) tends to $+\infty$ or $-\infty$ depending on the situation where it occurs. 
 % (see a remark given right after Theorem \ref{thm5} for an extended usage).
 We use the letters $x,y,z$ and  $w$ to represent integers  which indicate points assumed by the walk when discussing  matters on the random walk, while the same letters may stand for real numbers when  the stable process is dealt with;   we  shall sometimes use the Greek  letters $\xi$, $\eta$ etc. to denote the real variables  the stable process may assume. For convenience the normalizing sequence $c_n$ is extended to a continuous function on $[0,\infty)$ by linear interpolation. $C, C_1, C'$ etc. always designate positive constants that may change from line to line, whether they are significant or not.  

 The rest of this section is divided into four subsections. In the first three  we deal with the  special case  $B=\{0\}$, i.e.,  when the  killing takes place  at and only at the origin. In the first subsection {\bf 2.1}   we state certain fundamental results. For  $\ga< |2-\a|$ they  provide the precise asymptotic form of the transition probability of the killed walk within  the whole range $|x|\vee |y| =O(c_n)$, whereas for $\ga=|2-\a|$  the asymptotic form they provide
 is confined to a very small part of the range. The results for the remaining region
 are addressed in the subsections {\bf 2.2} and {\bf 2.3}. 
 The general case of finite sets---closely  parallel  to the special case $B=\{0\}$---are dealt with in the last  subsection {\bf 2.4}.

%When both the walk and the stable process get into the discussion at the same time and the differentiation is needed or suitable 

\v2\v2
{\bf 2.1.  Fundamental results for the walk killed at the origin.}
\v2

Let $f^x(n)$ denote the probability that the walk started at $x$  visits the origin at  $n$ for the first time: 
$$f^x(n) = P[\sigma^x_{\{0\}} =n].$$
Put
$$\k_{\a,\ga} =\k_{\a, -\ga} = \frac{(\a-1)\sin \frac{\pi}{\a}}{\Ga(\frac1{\a})\sin \frac{ \pi(\a -\ga)}{2\a} } = \frac{(1-\frac1{\a})\sin \frac{\pi}{\a}}{\mathfrak{p}_1(0) \pi}$$
(see (\ref{p0}) for the second equality); in particular if $\ga=|2-\a|$, $\k_{\a,\ga} =(\a-1)/\Ga(1/\a)$.  

We know  (cf. \cite[Lemma 2.1]{B1}) that as $n\to\infty$
\beqn\label{S_f}
P[\sigma^0_{\{0\}} >n] \sim \k c_n/n \qquad (\k= \k_{\a,\ga}/(1-{\textstyle \frac1{\a}})),
\eeqn
which is shown by the standard method based on Karamata's Tauberian theorem (a similar result is given in  \cite[Theorem 5.1]{PS} for non-lattice case). Our first theorem is the local version of (\ref{S_f}).

\v2

 %Thm1
\begin{Thm}\label{thm1} \,  For any admissible $\ga$,  as  $n\to\infty$
$$f^0(n) \sim   \frac{ \k_{\a,\ga}c_n}{n^{2}}.$$
\end{Thm}
\v2

For a non-empty subset $B\subset \Z$ put
\beqn\label{trans}
Q^n_B(x,y) := P[S^x_n=y, \sigma^x_B \geq n],
\eeqn
which entails  $Q^n_B(x,y)=  P[S^x_n =y, \sigma_B^x=n]$  ($n\geq1, y\in B, x\in \Z$) as well as $Q^0_B(x,y)=\1(x=y)$;  and similarly for a closed set  $\De \subset \R$
\[
 \mathfrak{p}^\De_t(\xi,\eta) := P[Y^\xi_t\in d\eta, \sigma_\xi^\De >t]/d\eta.
 \]
where $Y^\xi_t =\xi +Y_t$ and $\sigma^\Delta_\xi$ is the first entrance time of $Y^\xi$ into $\De$.  By the scaling law for stable processes we have
\[
  \mathfrak{p}^\De_n(x,y) = n^{-1/\alpha} \mathfrak{p}^{\De/n^{1/\alpha}}_1(x/n^{1/\a},y/n^{1/\a}).
\]  

\v2

 Let $c_n$ satisfy (\ref{NC}) and   write  $x_n$ for $x/c_n$  and similarly for  $y_n$.
 We are primarily interested in the asymptotic form of $Q^n_B(x,y)$ in the space-time region  $|x|\vee |y| <Mc_n$ with an arbitrarily given constant $M$.  It differs in a significant way according as 
 $|\ga| <2-\a$ or $|\ga| =2-\a$. First we state the result for the  former case that is formulated in a rather neat form.

%
%Thm2
\begin{Thm}\label{thm2}\,  If $|\ga|<2-\a$, then  for any $M>1$, uniformly for $|x| \vee |y|<Mc_n$, as $n\to\infty$
\beqn\label{eq_thm2}
Q^n_{\{0\}}(x,y) \sim \left\{\begin{array}{lr}
{\displaystyle  a^\dagger (x)f^{0}(n)a^\dagger(-y) } \quad &( |x_n| \vee |y_n| \to 0),\\[2mm]
{\displaystyle \mathfrak{f}^{\,x_n}(1)a^\dagger(-y)/n} &  (y_n \to 0, |x_n|>1/M),
\\[2mm]
{\displaystyle  a^\dagger (x)\, \mathfrak{f}^{\,-y_n}(1)/n } \quad &( x_n \to 0, |y_n|>1/M),\\[2mm]
\mathfrak{p}^{\{0\}}_{1} (x_n,y_n)/c_n & (|x_n| \wedge |y_n| > 1/M).
\end{array} \right.
\eeqn
\end{Thm}
\v2

 If $L\equiv 1$, i.e., the walk is in the domain of  normal attraction, then, in view of asymptotics of $a$ given later in Remark \ref{rem2.2},  (\ref{eq_thm2})  restricted to the case  $|x|\wedge |y| \to\infty$ reduces to $Q^n_{\{0\}}(x,y) \sim \mathfrak{p}^{\{0\}}_1(x_n,y_n)/c_n = \mathfrak{p}^{\{0\}}_n(x,y)$ (see (\ref{xi/eta})), so that the transition probability of the killed walk behaves quite similarly to that of its limit process if $|\ga|<2-\a$.

 For $\ga=\pm(2-\a)$ the above formula (\ref{eq_thm2})  holds  in and only in the partial range  of $|x|\vee |y|< Mc_n$ that is identified in the next proposition.    The range of validity  is described by means of the function
 $$\La_n(x) := a^\dagger(x)c_n/n.$$
Slightly abusing notation, we write  $|x_n|<\!< \La_n(x)$ instead of $|x|/a^\dagger(x) <\!< c_n^2/n$ with the convention: if the walk is left-continuous (resp. right-continuous) so that $a(x)=0$ for $x>0$ (resp. $x<0$),  this condition  is understood to be violated    for $x > 0$ (resp.  $x<0$). 
We have only to give results for  $\ga=2-\a$ when $a(x)/a(-x)\to 0$ ($x\to\infty$), the results for $\ga=-2+\a$ being obtained  from them  by  exchanging $x$ and $y$.

%%%%%%%

\begin{figure}[t]
 \begin{center}
\begin{picture}(405,145)(30,45)

%(Left)
\put(85,40){{\footnotesize $\gamma=2-\alpha$}}

\put(25,120){\vector (1,0){150}}
\put(167,125){{\footnotesize  {\sc x}}}
\put(100,60){\vector (0,1){126}}
\put(104,178){{\footnotesize {\sc y}}}

\thicklines
%\put(40,145){\line (1,0){118}}
\put(40,95){\line (1,0){118}}
%\put(75,70){\line(0,1){100}}
\put(125,70){\line(0,1){100}}

\put(130,99){{\footnotesize $ y_n = -1/M$}}
%\put(160,150){{\footnotesize $ y_n = $}}
%\put(165,140){{\footnotesize $ 1/M$}}
\put(128,72){{\footnotesize $x_n = 1/M$}}

%vertical dashed line
%\thinlines
\put(106, 70){\line(0,1){4}}
\put(106, 76){\line(0,1){4}}
\put(106, 82){\line(0,1){4}}
\put(106, 88){\line(0,1){4}}
\put(106, 94){\line(0,1){4}}
\put(106, 100){\line(0,1){4}}
\put(106, 106){\line(0,1){4}}
\put(106, 112){\line(0,1){4}}
\put(106, 118){\line(0,1){4}}
\put(106, 124){\line(0,1){3}}
\put(106, 129){\line(0,1){3}}
\put(106, 134){\line(0,1){4}}
\put(106, 140){\line(0,1){4}}
\put(106, 146){\line(0,1){4}}
\put(106, 152){\line(0,1){4}}
\put(106, 158){\line(0,1){4}}
\put(106, 164){\line(0,1){4}}

%horizontal dashes line

\put(40,112){\line(1,0){2}}
\put(44,112){\line(1,0){4}}
\put(50,112){\line(1,0){4}}
\put(56,112){\line(1,0){4}}
\put(62,112){\line(1,0){4}}
\put(68,112){\line(1,0){4}}
\put(74,112){\line(1,0){4}}
\put(80,112){\line(1,0){4}}
\put(86,112){\line(1,0){4}}
\put(92,112){\line(1,0){4}}
\put(98,112){\line(1,0){4}}
\put(104,112){\line(1,0){4}}
\put(110,112){\line(1,0){4}}
\put(116,112){\line(1,0){4}}
\put(122,112){\line(1,0){4}}
\put(128,112){\line(1,0){4}}
\put(134,112){\line(1,0){4}}
\put(140,112){\line(1,0){4}}
\put(146,112){\line(1,0){4}}
\put(152,112){\line(1,0){4}}

\thinlines
%(1)
\put(107,112){\line(0,1){57}}
\put(109,112){\line(0,1){57}}
\put(111,112){\line(0,1){57}}
\put(113,112){\line(0,1){57}}
\put(115,112){\line(0,1){57}}
\put(117,112){\line(0,1){57}}
\put(119,112){\line(0,1){57}}
\put(121,112){\line(0,1){57}}
\put(123,112){\line(0,1){57}}

%(2)
\put(106,97){\line(-1,0){66}}
\put(106,99){\line(-1,0){66}}
\put(106,101){\line(-1,0){66}}
\put(106,103){\line(-1,0){66}}
\put(106,105){\line(-1,0){66}}
\put(106,107){\line(-1,0){66}}
\put(106,109){\line(-1,0){66}}
\put(106,111){\line(-1,0){66}}
%\put(106,101){\line(-1,0){66}}

\put(40,114){$....................$}
\put(42,117){$...................$}
\put(42,122){$...................$}
\put(40,125){$....................$}
\put(42,128){$...................$}
\put(40,131){$....................$}
\put(40,134){$....................$}
\put(78,137){$........$}
\put(42,137){$........$}
\put(40,140){$........${\footnotesize  (1)}$........$}
\put(40,143){$........$}
\put(80,143){$........$}
\put(42,146){$........$}
\put(78,146){$........$}
\put(42,149){$...................$}
\put(40,152){$....................$}
\put(42,155){$...................$}
\put(40,158){$....................$}
\put(42,161){$....................$}
\put(40,164){$....................$}
\put(42,167){$...................$}

\put(126,114){$.........$}
\put(128,117){$.........$}
\put(127,122){$.........$}
\put(128,125){$.........$}
\put(126,128){$.........$}
\put(128,131){$.........$}
\put(126,134){$.........$}
\put(128,137){$.........$}
\put(126,140){$...${\footnotesize  (2)}$...$}
\put(128,143){$..$}
\put(150,143){$..$}
\put(126,146){$...$}
\put(148,146){$...$}
\put(128,149){$.........$}
\put(126,152){$.........$}
\put(128,155){$.........$}
\put(126,158){$.........$}
\put(128,161){$.........$}
\put(126,164){$.........$}
\put(127,167){$.........$}

\put(40,92){$....................$}
\put(42,89){$...................$}
\put(82,86){$.......$}
\put(42,86){$........$}
\put(84,83){$......$}
\put(40,83){$........$}
\put(40,80){$........${\footnotesize  (2$^*$)}$.......$}
\put(42,77){$...................$}
\put(40,74){$....................$}
\put(40,71){$....................$}

%%%%
%%%

\put(230,160){{\footnotesize {\sc Figure 1}:  The dotted regions  indicate the }}

\put(230,150){{\footnotesize   range of validity for (\ref{eq_thm2}) described by (1), (2) }}

\put(230,140){{\footnotesize   and  (2$^*$)  in  Proposition \ref{prop1} and the striped    }}
\put(230,130){{\footnotesize
regions  the range that is added by Theorem \ref{thm3}   }}
\put(230,120){{\footnotesize  both  in the extreme case $\ga = 2-\a$. The vertical }}
\put(230,110){{\footnotesize     broken  line roughly indicates   the boundary for  }}
\put(230,100){{\footnotesize    the regime $x_n <\!< \La_n(x)$ and the horizontal  one }}
\put(230,90){{\footnotesize  that for $-y_n <\!< \La_n(-y)$.  }}

%The vertical broken  line roughly  indicates where the boundary of the regime $x_n <\!< \La_n(x)$ (resp. $$) is located.

\end{picture}
\end{center}
% \caption{{\bf u} is  the  point of intersection where the half line  $\{(\eta,a): \eta>0\}$ meets
%  the  line passing through $\x$ and   tangential to  the circle   $\partial U(a)$.}

\vspace*{0cm}
\end{figure}
\v2\v2

%%%%%%%%%%%

%prop2.1
\begin{Prop}\label{prop1}
  If $\ga =2-\a$ (when the limiting  stable process has no negative jumps), then  for any $M>1$, (\ref{eq_thm2}) holds uniformly in $x, y$ under  each of the following constraints:
\v2
{\rm (1)} $\quad -M< x_n <\!< \La_n(x) \quad \mbox{and}\quad -\La_n(-y)  <\!< y_n < M$;
\v2
{\rm (2)} $\quad  1/M < x_n <M \quad \mbox{and}\quad -\La_n(-y)  <\!< y_n < M$;
\v2
{\rm (2$^*$)} $\quad -M< x_n <\!< \La_n(x) \quad \mbox{and}\quad - M< y_n < -1/M$.
\v2
[Recall the expression $-\La_n(-y)  <\!< y_n$ is the same as $- y_n <\!< \La_n(-y) $.]
\end{Prop}
\v2

By  (\ref{NC})  together with   Remark \ref{rem2.1}(a)  given in the next subsection it follows that  
 $$\La_n(\pm c_n) \to  \kappa^a_{\a,\ga,\pm} \quad \mbox{ with $\kappa^a_{\a,\ga,\pm} >0$ or $=0$ according as $\pm \ga \neq 2-\alpha$ or $=2-\a$}.$$
   Hence  for $x\geq 0$, the condition  $x_n <\!<  \La_n(x)$  always entails $x_n \to 0$ and if $\ga <2-\a$ the converse is true, so that for
 $\ga< 2-\a$ and $x\geq0$, the condition $x <\!< \La_n(x)$ is the same as $x_n \to 0$; 
  in case $\ga =2-\a$ it  largely depends on the behaviour of $a(x)$ as $x\to\infty$ (cf. \cite{Upot} for how it behaves) but still permits $x$ to grow indefinitely  large  unless the walk is left-continuous; and similarly for $x\leq 0$. From 
 Theorem \ref{thm2} and Proposition \ref{prop1} (see Corollary \ref{cor1} below) we infer that $\La_n(x) \sim \kappa^{-1}P[\sigma^x_{\{0\}}>n]$ ($\k =(1-\frac1\a)^{-1}\kappa_{\a,\ga}$) whenever $|x_n| <\!< \La_n(x)$.

%Rem2.1
\begin{rem}\label{rem0}
(a) \,  For    fixed  $x, y$,   (\ref{eq_thm2}) (reducing to $Q^n_{\{0\}}(x,y) \sim a^\dagger(x)a^\dagger(-y) f^0(n)$) follows from a special case of the result of Kesten \cite{K}   (see Remark \ref{rem4C} given at the end of this section).  When $\ga=0$ (i.e., the limit stable process is symmetric) and $L(x) \to 1$, the  asymptotic form of  $f^0(n)$ is derived  in  \cite{K} in which an asymptotic form  for $\a=1$ is also obtained,  which reads $f^0(n) \sim \pi /n(\log n)^2$.  

\v2
(b) \, Belkin \cite[Theorem 2.1]{B1} shows the  conditional limit theorem that may read 
\beqn\label{Belkin}
\lim_{n\to\infty} P[ S_n/c_n > \xi \,|\, \sigma^0_{\{0\}} >n] = \int_\xi^\infty h(t)dt \quad \mbox{for} \quad \xi\in \R,
\eeqn
where $h$ is a probability density on $\R$ which is bounded and  continuous and whose characteristic function $\phi_h(\th):=\int_{-\infty}^\infty  e^{i\th\xi}h(\xi)d\xi $ is given by
\beqn\label{eqL03}\phi_h(\th) = 1- \psi(\th)\int_0^1 \xi^{\frac1{\a}-1}e^{  \psi(\th)(\xi -1)}d\xi.
\eeqn
 Our results, Theorem \ref{thm2} and Proposition \ref{prop1},  imply the local version of this limit theorem. On taking $x=0$ in the third case of (\ref{eq_thm2}),  with the help of (\ref{S_f}) and (\ref{Belkin}) we especially find that for $\xi\in \R$,
   $ P[ S_n = \lfloor c_n\xi \rfloor\,|\, \sigma^0_{\{0\}} >n]  \sim h(\xi)/c_n$ ($n\to\infty$)
   %  locally uniformly in $\xi$ on $\R\setminus \{0\}$% and uniformly for $x$ satisfying $x_n <\!<\La_n(x)$; 
   and
   \beqn\label{h/f/B}
   h(\xi) = [(1-1/\a)/\kappa_{\a,\ga}] \mathfrak{f}^{\,-\xi}(1).
   \eeqn
 (This identity is verified directly from (\ref{eqL03}) in our proofs of Theorem \ref{thm2}  and Proposition \ref{prop1}.)
  \end{rem}
 %Here the limit is independent of $x$ for which the restriction  $x_n <\!<\La_n(x)$ is crucial in case $\ga=|2-\a|$
% where  this restriction is not implied by $x_n \to 0$  (see Corollary \ref{cor-1}). 

%Our results also  entails the local limit theorem that says that   for all admissible $\ga$ and for   $|x_n|, |y_n| \in [1/M, M]$,  $Q^n_{\{0\}}(x,y) \sim \mathfrak{p}_{1}^{\{0\}}(x_n,y_n)/c_n$ except for the case  
% when $\ga=|2-\a|$,  $\ga x>0$ and $\ga y<0$.  In this exceptional case we have    $\mathfrak{p}_{1}^{\{0\}}(x_n,y_n)=0$ and some estimates of $Q^n_{\{0\}}(x,y)$ are given in the next subsection (see Theorem \ref{thm4} and Proposition \ref{prop2.4}(ii)).

    We state a consequence on 
  $f^x(n)$ that directly follows from the results given above by virtue of the identity
  $ f^x(n) =  Q_{\{0\}}^n(x,0).$  

 %Cor1
\begin{Cor}\label{cor1} \,   For each admissible $\ga$ and  $M>1$,  as  $n\to\infty$
\beqn \label{eq_thm1}
f^x(n) \sim \left\{\begin{array}{lc} {\displaystyle  a^\dagger(x) f^0(n) } \qquad  \mbox{if} \;\; |x_n| <\!< \La_n(x),\\[2mm]
 \mathfrak{f}^{\, x_n}(1)/n \qquad  \mbox{uniformly for}\; \;  |x_n| \in[1/M, M].
\end{array}\right.
\eeqn
\end{Cor}
\v2
By (\ref{eq_thm1}) the first three formulae of (\ref{eq_thm2}) are written  as
$$ Q^n_{\{0\}}(x,y) \sim \left\{ \begin{array} {ll} f^x(n)a^\dagger(-y)  \quad (y_n\to 0), \\
a^\dagger(x) f^{-y}(n) \quad (x_n\to 0),
\end{array}\right. \qquad \mbox{uniformly for} \;\; |x_n|\vee|y_n|<M
$$
valid if $\ga\neq |2-\a|$ or if $\ga=2-\a$ and  $(x,y)$ is in  the range specified by Proposition \ref{prop1}.

Taking  (\ref{h/f/B}) into account, from these formulae one readily deduces the following local conditional limit theorem. 
%Cor2
\begin{Cor}\label{loc_clt} \, If $\ga<|2-\a|$,  as $n\to \infty$ uniformly for $|x| <\!< c_n$ and $|y| <Mc_n$,
$$
P[S^x_n =y\, |\, \sigma^x_{\{0\}}>n] \sim \left\{\begin{array}{lr}
(1-\frac1{\a})\La_n(-y)/c_n &\quad (y_n \to 0),   \\ [2mm]
h(y_n)/c_n & (|y| \geq c_n/M).
\end{array}\right.
$$
[For $|y_n|\geq 1$ we have  $O(nL(y)/c_n |y|^\a)$ as an upper bound of the  LHS (see   Lemma \ref{L6.1}).]
\end{Cor}

%For $\ga< 2-\a$ and $x\geq0$ the first case of (\ref{eq_thm1}) may read
%$$ f^x(n)\sim a^\dagger(x) f^0(n) \qquad \mbox{as} \;\; n\to\infty\;\;  \mbox{and}\;\; x_n\to  0.$$
% The asymptotic formulae given in 
 Corollary \ref{loc_clt}  entails that  the conditional distribution  
 $P[S^x_n/c_n \in d\xi\,|\, \sigma^x_{\{0\}}>n]$ converges to $h(\xi)d\xi$ ($n\to\infty$) \lq uniformly' for $|x|<\!< c_n$. 
 For $\ga =2-\a$,  this remains true  for $-1 <\!< x_n <\!<\La_n(x) $ but does not for $\e \La_n(x)\leq x_n <\!<1$ any more for each $\e>0$    as is revealed in the next subsection  (see (\ref{crs_ov}), (\ref{mix})) in which we address the issue for those  $x,y$ that  satisfy   $\e \La_n(x) < x_n <\!< 1$, $-\La_n(-y) <\!< y_n<M$, hence   are  excluded   from  Proposition \ref{prop1} (see {\sc Figure 1})---we have included Proposition \ref{prop1} in this section, it being proved by essentially the same arguments that derive Theorems \ref{thm1} and \ref{thm2}. %Before proceeding to  further results here 
.   

\v2\v2

%2.2  Case  $\ga=2-a$ with $\e \La_n(x) < x_n$.}
{\bf 2.2. Case  $\ga=2- \a$ with $\e \La_n(x) \leq x_n <\!< 1$, $-\La_n(-y)<\!<y_n <M$.}
\v2

  Let $U_{{\rm d}}(x)$, $x=0, 1, 2,\ldots$
 denote the renewal function of the strictly descending ladder height process: $U_{{\rm d}}(0)=1$ and for $x\geq 1$,
 $$U_{{\rm d}}(x) = 1+ \sum_{z=1}^x u(z) \quad\mbox{and}\quad  u(x)= \sum_{n=1}^\infty  P[\hat Z_1+\cdots + \hat Z_n = -x],$$
 where $\hat Z_1 = S_{\sigma_{(-\infty,-1]}}$  and $\hat Z_k$ are independent copies of  $\hat Z_1$. 
 
%Thm3
\begin{Thm}\label{thm3}\,  Let $\ga=2-\a$. Then for each $M>1$ uniformly  for $0\leq x <Mc_n$, as $n\to\infty$
\beqn\label{eq_thm3}
Q^n_{\{0\}}(x,y) \sim  \left\{\begin{array}{lr}
{ \displaystyle  \bigg\{ a^\dagger(x)f^0(n) +\frac{U_{{\rm d}}(x)\mathfrak{p}_1(-x_n)}{U_{{\rm d}}(c_n)n} \bigg\}a^\dagger(-y)   } \quad  & (-\La_n(-y) <\!< y_n <\!<1),
\\[4mm]
{\displaystyle   \frac{a^\dagger(x)}{n}\mathfrak{f}^{\,-y_n}(1) + \frac{U_{{\rm d}}(x)K_1(y_n)}{U_{{\rm d}}(c_n)c_n}   }\quad\;\;  & (x_n \to 0, y_n \in [1/M, M]),
\end{array} \right.
\eeqn
where  $K_{t}(\eta) =0 $ for $\eta\leq 0$ and 
\beqn\label{Kt}
K_{t}(\eta) =\lim_{\xi\downarrow 0} \frac1{\xi}\mathfrak{p}^{(-\infty,0]}_{t} (\xi,\eta) \quad \mbox{for} \quad \eta>0.
\eeqn
\end{Thm}
\v2

We know that
 if $\ga = 2 -\a$, then $U_{{\rm d}}$ varies regularly of index 1 [24] and
\beqn\label{eqR(b)}
\mathfrak{f}^{\,x}(t) = xt^{-1}\mathfrak{p}_t(-x)\quad (x>0) \quad \mbox{and}\quad 
\mathfrak{p}^{\{0\}}_{t} (x,y) =\mathfrak{p}^{(-\infty,0]}_{t} (x,y)\;\; (x, y>0)
\eeqn
  (cf., e.g.,  \cite[Corollary 7.3]{Bt}). Note that the first case of (\ref{eq_thm3}) conforms to the second one of (\ref{eq_thm2}) valid for  $x_n \asymp 1$ in view of Proposition \ref{prop1}.
 
By Proposition \ref{prop1} we have $Q^n_{\{0\}}(x,y)\sim a^\dagger(x)f^0(n)a^\dagger(-y) $ for $-M<x_n <\!< \La_n(x)$.   Comparing this with  the   first relation of (\ref{eq_thm3}) shows that  for $x\geq 0$, $x_n <\!< \La_n(x)$ implies $U_{{\rm d}}(x)/U_{{\rm d}}(c_n) <\!< \La_n(x)$. The converse implication is obvious because of the bound  $U_{{\rm d}}(x)/U_{{\rm d}}(c_n) \geq x_n\{1+o(1)\}$  ($x\leq c_n$) (see Lemma \ref{L8.8}(ii)).   It therefore follows that 
 for $0\leq x <Mc_n$,
\beqn\label{La/x}
U_{{\rm d}}(x)/U_{{\rm d}}(c_n) <\!< \La_n(x) \quad\mbox{if and only if} \quad x_n <\!< \La_n(x).
\eeqn

By (\ref{La/x}) it follows that in  both cases of (\ref{eq_thm3}) the second terms are superfluous if    $x_n<\!< \La_n(x)$.
Taking this into account  we reformulate  the result for $\ga=2-\a$ as  the following
%Cor3
\begin{Cor}\label{cor2} Let $\ga=2-\a$. Then for any   $M>1$, 

{\rm (i)} \, uniformly  for $x$, as $n\to\infty$
\beqn\label{eq_cor2}
f^x(n) \sim  \left\{\begin{array}{lr}
a^\dagger(x)f^0(n)  \qquad & (-1<\!< x_n <\!< \La_n(x)),\\
{\displaystyle  a^\dagger(x)f^0(n)  + \frac{U_{{\rm d}}(x)\mathfrak{p}_1(-x_n)}{U_{{\rm d}}(c_n)n} } \quad  & (M^{-1}\La_n(x)<x_n <M),
\end{array}\right.
\eeqn

{\rm (ii)} \, uniformly for $- \La_n(-y)<\!< y_n <\!< 1$ and  $|x_n|<M$
\beqn\label{eq_cor21} Q^n_{\{0\}}(x,y) \sim  f^x(n)a^\dagger(-y).
\eeqn
%and uniformly for $1/M\leq y_n <M$.
%\beqn\label{eq_cor22}
%Q^n_{\{0\}}(x,y) \sim  \left\{\begin{array}{lr}
%{ \displaystyle  a^\dagger(x)f^{-y}(n)   } \quad  & (-1 <\!<x_n <\!<\La_n(x)),
%\\[1mm]
%{\displaystyle   a^\dagger(x)f^{-y}(n) + \frac{U_{{\rm ds}}(x)K_1(y_n)}{U_{{\rm ds}}(c_n)B_n}   }\quad\;\;  & (M^{-1}\La_n(x) <x_n <\!<1).
%\end{array} \right.
%\eeqn
\end{Cor}
\v2

  In view of  (\ref{La/x}) it follows from  (\ref{eq_cor2}) that  as $n\to\infty$
\beqn\label{crs_ov}
f^x(n) \sim \left\{\begin{array}{ll} \k_{\a,\ga}  \, a^\dagger(x)c_n/ n^{2} \quad & (0\leq x_n <\!<  \La_n(x)), \\[2mm]
{\displaystyle \frac{U_{{\rm d}}(x)\mathfrak{p}_1(-x_n)}{U_{{\rm d}}(c_n)n} } \quad & (\La_n(x) <\!< x_n <M),
 \end{array} \right.
 \eeqn
  and  the two terms on the RHS are comparable with each other if $\e \La_n(x)<\!< x_n <\!< \e^{-1}\La_n(x)$ for each $\e>0$.  %On comparing with the result for $\sigma^x_{(-\infty,0]}$ obtained by Doney \cite{D} 
  The same kind of crossover    takes place in (\ref{eq_thm3}): in the both cases of it  the crossover occurs around  $x \asymp c_n\La_n(x)$ as in (\ref{crs_ov}). 
 % On restricting to the case $y_n \in [1/M,M]$  (i.e., the second  case of (\ref{eq_thm3}))
   This may be exhibited in a noticeable  form of the local limit theorem that corresponds to Corollary \ref{loc_clt}.
%for the conditional law of $S^x_n$ conditioned to avoid  $0$ which we state as a corollary. 
Let $q(\xi)$ designate  the density of the stable meander  of length $1$ at time $1$ (see (\ref{Bt})). 
Recall that  $h$ is the density  in  the conditional limit theorem (\ref{Belkin}). We know $h(\xi)>0$ for $\xi\neq 0$, while $q(\xi)= 0$ for $\xi\leq 0$. 
%Cor4
\begin{Cor}\label{cor-1} If $\ga=2-\a$, then
   uniformly for $0 \leq x_n <\!< 1$ and $y\in \Z$,
  \beqn\label{mix}
P[S^x_n = y \,|\, \sigma^x_{\{0\}}>n] =  \frac{\th_{n,x} h(y_n)+ (1-\th_{n,x}) q(y_n) }{c_n} +  o\Big(\frac{r_{n,x,y} +1}{c_n}\Big),
\eeqn
where  
$$\th_{n,x} = \frac{P[\sigma^x_{(-\infty,0]} < n <  \sigma^x_{\{0\}} ]}{ P[n< \sigma^x_{\{0\}} ]}\quad \mbox{ and}\quad  r_{n,x,y} = \1(0\leq y<c_n) \frac{L(c_n)}{L(y)}.$$ 
\end{Cor}

 Plainly   the convex combination of $h$ and $q$ on the RHS of   (\ref{mix}) is  proper  if and only if $ P[\sigma^x_{(-\infty,0]} > n] \asymp  P[\sigma^x_{\{0\}} > n]$, of which an  analytic expression   is provided by $\La_n(x) \asymp x_n$. Although Corollary \ref{cor-1}  is almost a corollary of the proof of Theorem \ref{thm3}   its proof  will be given in  Section 6 (given after Proposition  \ref{prop6.2})   since  we need to apply results proved therein in order to verify the uniformity 
in $y$  outside $[0, c_n]$.

\v2
%Remark 2.2.
\begin{rem} \label{rem2.1}

 (a)
We shall see   (in  Lemma \ref{lem3.1}) that   as $|x|\to\infty$
\beqn \label{a(x)}
\frac{a(x) L(x)}{|x|^{\alpha-1}} \, \longrightarrow\,   \left\{\begin{array}{lll} 0 \quad  &\mbox{if} \quad  \ga x >0, |\ga| = 2-\a \\
\k^a_{\a,\ga,\, {\sg x}}>0\; \quad&\mbox{otherwise} 
\end{array}\right.
 \eeqn
where     $\k^a_{\a,\ga,\, {\sg x}}$ is a constant (depending on $\a, \ga$ and  ${\sg x}$) which is positive if $\ga\,{\sg x} \neq   2-\a$ and equals $1/\Ga(\a)$ if $\ga\,{\sg x} = -  2+\a$. In particular
\v2

{\it if \, $\ga=2-\a$,}
\beqn\label{a(-y)}
\begin{array}{ll}
 \;\; & a(x) \sim \left\{
\begin{array}{lr} [\Ga(\a)]^{-1}|x|^{\a-1}/L(x) \quad &\mbox{as} \;\;  x\to -\infty,\\[1mm]
o(a(-x)) \quad  &\mbox{as} \;\;  x\to \infty,  
\end{array}\right. \;\mbox{and}
\\[5mm]
 & a(-c_n) \sim [\Ga(\a)]^{-1}n/c_n \;\;\mbox{and}  \;\; \La_n(c_n) \to 0 \quad \mbox{as} \;\; n\to\infty.
\end{array}
\eeqn
%In case  $\ga=2-\a$ and $x\to \infty$ the behaviour of $a(x)$ varies according to that of  $P[X<x]$, but we know, in addition to the facts mentioned above as well as   at the beginning     of this section, that $a(x) \sim 1/P[\sigma^0_{(-\infty, -x]}<\sigma^0_{\{0\}}]$, entailing  \lq asymptotic monotonicity' \cite[Corollary 5.2]{Uattrc}.

%\v2
%(b)\,  For $\ga \neq 2-\a$,   as $x\to\infty$ under $x_n\to 0$
%$$f^x(n)\sim  \k_{\a,\ga}  \, a(x) c_n/ n^{2}\sim {\rm const_\cdot} \frac{ \mathfrak{f}^{\, x_n}( 1)L(c_n)}{nL(x)}  $$
 % (see (\ref{cor0}) for the latter equivalence). If  $X$ belongs to the domain of normal attraction,  the two expressions on the RHS of (\ref{eq_thm1}) are asymptotically equivalent to each other in this regime, so that   $f^x(n)  \sim \mathfrak{f}^{x_n}(1)/n$ as $x\to\infty$ under $x_n<M$.

\v2
(b) \,In (\ref{NC}) one may replace 1 by any positive constant $\la$, or what is the same thing, replace $L$ by $\tilde L = L/\la$ and $\psi$ by $\tilde \psi =\la \psi$  so that  
$n c_n^{-\a}\tilde L(c_n) \to 1$ and $S_n/ c_n \stackrel{{\rm law}} \longrightarrow  Y_{\la} $.
This causes merely  the replacements of the functions associated with  $Y_t$ by those with $\tilde Y_{t}=Y_{\la t}$, e.g.,  $\mathfrak{f}^{x_n}(1) \to \la\mathfrak{f}^{x_n}(\la)$, $\mathfrak{p}^{\{0\}}_1(x_n,y_n) \to \mathfrak{p}_{\la}^{\{0\}}(x_n, y_n)$, and so on. (Note that the formula of Theorem \ref{thm1} becomes $f^0(n) \sim \k_{\a,\ga} \la^{1/\a}c_n/n^2$ since $\mathfrak{p}_{\la}(0) =\la^{-1/\a}\mathfrak{p}_1(0)$.)
\end{rem}

%\v2
%(b) \,In (\ref{NC}) one may replace 1 by any positive constant $\la$, or what effects  the same thing, replace $c_n$  by  $\tilde c_n=\la^{-1/\a} c_n$ so that if   $\tilde Y_t = Y_{\la t}$ and $\tilde L = \la L$, then 
%$S_n/\tilde c_n \stackrel{{\rm law}} \longrightarrow \tilde Y_{1} $  and $n\tilde c_n^{-\a}\tilde L(\tilde c_n) \to 1$.
%This causes merely  the replacements of the functions associated with  $Y_t$ by those with $\tilde Y_{t}=Y_{\la t}$ (along with the replacements:\,$c_n \to \tilde c_n$, $x_n \to \tilde x_n= \la^{1/\a} x_n$), e.g.,  $\mathfrak{f}^{x_n}(1) \to \la\mathfrak{f}^{\tilde x_n}(\la)$, $\mathfrak{p}^{\{0\}}_1(x_n,y_n) \to \mathfrak{p}_{\la}^{\{0\}}(\tilde x_n, \tilde y_n)$, and so on. (Note that the formula of Theorem \ref{thm1} becomes $f^0(n) \sim \k_{\a,\ga} \la^{1/\a}\tilde c_n/n^2$ since $\mathfrak{p}_{\la}(0) =\la^{-1/\a}\mathfrak{p}_1(0)$.)

\vskip2mm

Below we provide some consequences on $\mathfrak{p}^{\{0\}}_{1}(\xi, \eta)$ that are drawn from   Theorems \ref{thm2} and \ref{thm3} (and that the present author fails to find explicit statements for  in the existing literature) and also known facts on  $\mathfrak{f}^{\,\xi}(t)$.
\v2

{\sc Some consequences on $\mathfrak{p}^{\{0\}}_{1}(\xi, \eta)$ and $K_t(\xi)$} 
\v2 

(i) \, Let $\gamma = 2-\alpha$. For the present purpose we may suppose $L\equiv 1$ and  $E|\hat Z|<\infty$ so that $c_n\sim n^{1/\a}$ and  $U_{{\rm d}}(x)/U_{{\rm d}}(c_n)\sim x_n$. Formula (\ref{eq_cor21})   implies that  if one takes the successive  limit as  first  $x_n\to \xi >0$, $y_n\to \eta >0$ as well as $n\to\infty$ and then   $\xi \vee \eta \to 0$, then
\[
 \frac{Q^n_{\{0\}}(x,y)} {f^x(n)a(-y)}=   \frac{Q^n_{\{0\}}(x,y) }{ \mathfrak{p}^{\{0\}}_{ n}(x, y)} \frac{\mathfrak{p}^{\{0\}}_{ n}(x, y)}{f^x(n)a(-y)} \to 1. 
\]
Since   the first ratio of the middle member approaches 1 by virtue of the last relation of (\ref{eq_thm3}), it therefore follows from the second formula of (\ref{eq_cor2}) (or  (\ref{Doney})) that   as $\xi \to 0$ and $\eta=y_n \to 0$ 
\beqn\label{eqR(d)}
%\frac{(c_\circ n)^{1-1/\a} x_n \mathfrak{p}_{c_\circ}(- x_n) f^{-y}(n) /c_\circ\k_{\a,\ga}  }
%{p^n_{\{0\}}(x,y)} \sim  
\frac{ \mathfrak{p}^{\{0\}}_{1}(\xi, y_n) n^{1-1/\a}}{\xi \mathfrak{p}_{1}(-\xi) a(-y)} 
\to  1.
\eeqn
%Letting $\xi<0$ and $\eta>0$ in the above we obtain a similar relation.
On noting  $\mathfrak{p}^{\{0\}}_{1}(-\xi, -\eta) = \mathfrak{p}^{\{0\}}_{1}(\eta, \xi)$ and using 
(\ref{a(-y)}) this  shows that  
\beqn\label{R3eq}
\mathfrak{p}^{\{0\}}_{1}(\xi, \eta ) \sim  \frac{\mathfrak{p}_{1}(0)}{ \Ga(\a)} \times \left\{\begin{array}{lr}  \xi \eta^{\a-1} \quad  &(\xi \downarrow 0,  \eta \downarrow 0),\\
 (-\xi)^{\a-1}(-\eta) \quad & (\xi \uparrow 0, \eta \uparrow 0).
\end{array} \right.
\eeqn
(See Lemma \ref{L8.6} for the case when one of  $\xi$ and $\eta$ is fixed.) In a similar way we deduce
$$\mathfrak{p}^{\{0\}}_{1}(\xi, \eta ) \sim \kappa_{\a,\ga}[\Ga(\a)]^{-2}(-\xi\eta)^{\a-1} \quad \quad  (\xi \uparrow 0,  \eta \downarrow 0).
$$
\v2

(ii) \, 
In the same way as above but with
\beqn\label{eqR(d)1}
\frac{ \mathfrak{p}^{\{0\}}_{1}(\xi, y_n) n^{1-1/\a}}{\mathfrak{f}^{\,\xi}(1) a(-y)} 
\to  1
\eeqn
in place of (\ref{eqR(d)})
 we deduce  % from Theorem \ref{thm2} 
  that if  $|\ga| <2-\a$,   
 \beqn\label{xi/eta}
 \begin{array} {rr}
 \mathfrak{p}^{\{0\}}_{1}(\xi, \eta) %\sim n^{-2+1/\a}a(n^{1/\a}\xi)a(-n^{1/\a}\eta) 
\sim  \k'_{\a,\ga} \{\sin[\tst12 \pi(\a + (\sg \xi) \ga)]\}\{\sin[\tst12 \pi(\a - (\sg \eta) \ga)]\} |\xi\eta|^{\a-1} &\\[4mm]
\mbox{as}\quad |\xi|\vee |\eta|\to 0,&
\end{array}
\eeqn
where $\k'_{\a,\ga}=\k_{\a,\ga}[\Ga(1-\a)/\pi]^2$.  Similarly
$\lim_{\xi\to \pm 0}  \mathfrak{p}^{\{0\}}_{1}(\xi, \eta)/\xi^{\a-1} = \k^a_{\a,\ga,\pm} \mathfrak{f}^{-\eta}(1)$ for $\eta >0$.
%These are valid, if restricted to $\xi<0, \eta >0$,  also  for $\ga=2-\a$.
\v2
(iii) \, Let $\gamma = 2-\alpha$.  Then  $\mathfrak{p}_{1}(0) = 1/\alpha \Gamma(1-1/\alpha)$   and by (\ref{R3eq}) 
\beqn\label{R3f1}
K_{1}(\eta) \sim  \frac{\mathfrak{p}_{1}(0)}{ \Ga(\a)} \eta^{\alpha-1} \qquad (\eta \downarrow 0).
\eeqn
The scaling relation of $K_t(\eta)$ is given by $K_t(\eta) =K_1(\eta/t^{1/\a})/t^{2/\a}$.
\v2
\v2

{\sc Asymptotic properties of  $\mathfrak{f}^{\,\xi}(t)$.}
\v2

The density function $\mathfrak{f}^{\,x}(t)$ satisfies the scaling relation
\beqn\label{scl_rl} 
\mathfrak{f}^{\,x}(\la t) = \mathfrak{f}^{\, x/t^{1/\a}}(\la)/t = \mathfrak{f}^1( \la t/|x|^\a)/|x|^\a \qquad (x\neq 0, t >0, \la > 0).
\eeqn 
In case $\ga =|2-\a|$,  expansions of $\mathfrak{f}^{\,x}(t)t = \mathfrak{f}^{\,x/t^{1/\a}}(1)$ into power series of $x/t^{1/\a}$ are known. Indeed, if $\ga = 2 -\a$,  owing to (\ref{eqR(b)})   the series expansion for $x>0$ 
is obtained from that of  $t^{1/\a}\mathfrak{p}_t(-x)$ which is found in   \cite{F},  while  for $x<0$, the  series expansion is  derived by Peskir \cite{P}. In the recent paper  
 \cite[Theorem 3.14]{KKPW} Kuznetsov et al. obtain  a similar  series expansion for all cases. Here we state the leading term for $\mathfrak{f}^1(t)$  and an error estimate (as ($t \to \infty$) that  are deduced from these expansions.

As $t  \to\infty$
\beqn\label{cor0}
  \mathfrak{f}^1(t) =   \left\{ \begin{array}{ll}  [-1/\Ga(-1/\a)]t^{-1-1/\a}\{1+ O(t^{-1/\a})\} \quad &\mbox{if}\quad 
\ga=2-\a, \\[2mm]
 \k^{\mathfrak{f},+}_{\a,\ga} \, t^{-2+1/\a} \{1+ O(t^{1-2/\a})\}\quad &\mbox{if}\quad 
\ga\neq 2-\a, 
\end{array} \right.
\eeqn 
where
\beq
\k^{\mathfrak{f}, +}_{\a,\ga} = \frac{\Ga(2-\a)\sin (\pi/\a) \,\sin [\frac\pi2 (\a +\ga)]}{\a \pi^2 \mathfrak{p}_1(0)}
=\frac{\sin \pi/\a}{ \pi \mathfrak{p}_1(0)}\int_0^\infty u^{1-\a}\mathfrak{p}_1'(-u)du. 
\eeq
[The first expression shows  that  $\k^{\mathfrak{f}, +}_{\a,\ga}$    is positive if  $\ga < 2-\a$ and zero if $\ga =2-\a$; see Lemma \ref{lem7.1} for the second expression.]

In  \cite[Theorem 3.14]{KKPW} an asymptotic expansion as  $t\to 0$  (see Lemma \ref{lem7.01}) is also obtained which  entails 
$$ \mathfrak{f}^1(t) =  [\a^3\Ga(\a)/\Ga(2-\a)]\kappa^{\mathfrak{f},+}_{\a,\ga}\,  t^{1/\a} + O(t^{1+1/\a}) \quad \mbox{as}\quad t\to 0.$$

We observe that the leading term in the second formula of (\ref{cor0}) is deduced from Corollary \ref{cor1}. To this end  we may let 
$L\equiv 1$, for which  if $\ga\neq 2-\a$,
$$\lim_{\xi\downarrow 0} \frac{\mathfrak{f}^\xi(1)}{\xi^{\a-1}} = \lim_{\xi\downarrow 0} \lim_{x_n \to \xi}\frac{nf^x(n)}{\xi^{\a-1}} 
= \lim_{\xi\downarrow 0}  \lim_{x_n\to\xi} \frac{a(x)\k_{\a,\ga}}{\xi^{\a-1}n^{1-1/\a}}  
=  \k^{a}_{\a,\ga,+}\kappa_{\a,\ga} = \k^{\mathfrak{f},+}_{\a,\ga}
$$
which  is the same  as what is to be observed   because of the scaling relation  (\ref{scl_rl}). 
%$$ \frac{\a \Ga(2-\a)\,\sin (\pi/\a)\sin [\frac\pi2 (\a +\ga)]}{\Ga(\frac1{\a}) \sin[ \frac{ \pi}{2}(\a - \ga)/\a]} $$

\v2\v2
 %2.3 [Case  \ga=2-\a xy<0]
{\bf 2.3. Case  $\ga=2-\a$ with $xy<0$ and  upper bounds outside the principal regime.}
\v2\v2
 From Theorem \ref{thm3} or Proposition \ref{prop1} 
  the regime  $x_n>\e \La_n(x), \, y_n <-\e \La_n(-y)$ is excluded for any $\e>0$, where     there arises  a difficulty in  estimating  $p^n_{\{0\}}(x,y)$ in general; below we give a result under an extra assumption on the tail  as $t\to-\infty$ of the distribution function
$$F(t) := P[X\leq t].$$
Suppose that $\ga=2-\a$ so that $F(-t) =o(1-F(t))$ ($t\to\infty$).   In  \cite[Theorem 2(iii)]{Uladd} a criterion for the limit
 $$C^+ :=  \lim_{x\to  +\infty}  a(x) \leq \infty$$
(which exists) to be  finite  is obtained.   Under the present assumption on $F$ it says that
 \beqn\label{a_bdd}
 \sum_{y=1}^\infty F(-y) [a(-y)]^2 <\infty \quad \mbox{and} \quad F(-2)>0
 \eeqn
is  necessary and sufficient  for   $0<C^+<\infty$.  
  \v2
%Thm4
\begin{Thm}\label{thm4}\, Suppose that $\ga=2-\a$ and  (\ref{a_bdd}) holds. Then, given  $M>1$, uniformly for $-M\leq y_n<0\leq x_n <M$,
% along with  $x_n\wedge y_n \to 0$
\v2
{\rm (i)} \quad $Q^n_{\{0\}}(x,y) \sim %\left\{\begin{array}{lr} 
{\displaystyle 
a^\dagger(x)a^\dagger(-y)f^0(n) +\frac{ a^\dagger(x)|y_n|\mathfrak{p}_{1}(y_n) +  a^\dagger(-y)x_n\mathfrak{p}_{1}(-x_n)} { n} }$ \\
\qquad \qquad \qquad \qquad \qquad\qquad\qquad\qquad\qquad \qquad\qquad \qquad\qquad \qquad\qquad $(x_n\wedge (-y_n) \to 0),$
\v2
{\rm (ii) }\quad 
$Q^n_{\{0\}}(x,y) \sim  {\displaystyle  C^+\frac{(x_n-y_n)\mathfrak{p}_1(y_n -x_n)}{ n}=\frac{C^+  \mathfrak{f}^{\, x_n-y_n}(1) }n }
 \qquad\quad   (x_n\wedge (-y_n) >  1/M) $,
\v2\n
  as $n\to\infty$.
\end{Thm}

An application  of  Theorem \ref{thm4} leads to the next result which  exhibits  a way  the condition   $C^+<\infty$  is reflected in   the behaviour of   the walk $S^x$, $x>0$: conditioned on   $S^x_n=y$, $y<0$  it enters $(-\infty,-1]$ without visiting the origin  \lq continuously' or by a very long jump for large $x, -y$  according as $C^+$ is finite or not.  Exactly the same behaviour of the pinned walk is observed in \cite{U1dm_f} in the case $\sigma^2<\infty$  but with the condition (\ref{a_bdd}) replaced by $E[|X|^3; X<0] <\infty$ which is equivalent to $\lim_{x\to +\infty} [a(-x)-x/\sigma^2] <\infty$.

%Prop2.2
  \begin{Prop}\label{prop2.2}\,  Let  $\ga=2-\a$. Then  for each  $M\geq 1$, as $n\to\infty$   under the constraint   $-Mc_n< y <0< x< Mc_n$
\begin{eqnarray}\label{jump}
&& P[S^x_{\sigma_{(-\infty,0]}}<- R \,|\,\sigma^x_{\{0\}}>n, S^x_n =y]   \nonumber
\\[2mm]
&&\quad  \longrightarrow  \left\{ 
\begin{array}{ll}
0 \quad\mbox{ as}  \;\;  R \to\infty \;\; \quad  \mbox{uniformly for} \;\; x, \,y\; &\mbox{ if}\quad C^+ <\infty, \\[1mm]
1 \quad\mbox{ as}  \;\;   x\wedge |y| \to\infty \quad  \mbox{ for each}\; R>0\;\; &\mbox{ if}\quad C^+ =\infty.
\end{array} \right.
\end{eqnarray}
\end{Prop}
\v2

The conditional probability on the LHS of (\ref{jump})  is  of interest  even for $y>0$  in  case  $0\leq x_n < M\La_n(x)$ when the walk $S^x$  conditioned  as above    visits the negative half line with a positive probability---as made explicit in   Corollary \ref{cor-1}.  The same proof will apply in this case to show that  the  formula  (\ref{jump}) remains valid for  $0\leq x_n < M\La_n(x)$, $0<y< Mc_n$  if  its LHS is replaced by
$$P[S^x_{\sigma_{(-\infty,0]}}<- R \,|\,\sigma^x_{\{0\}}>n> \sigma^x_{(-\infty,0]}, S^x_n =y].$$

%We state  the  following upper bound as a proposition, a unified but partly reduced version obtained  by combining theorems above and the results  in  Section 6.  Write for $x>0$
%$$T(x)= F(-x) +1-F(x) \quad \mbox{and}\quad \bar a(x) =\frac12 [a(x)+a(-x)]$$
%Our assumption implies  $F(x) \sim BL(x)/x^\a$ for some constant $B>0$ (see Appendix (A) for related matters) and by Lemma \ref{lem3.1} $\bar a(x) \sim$

In  case   $|x_n|\to\infty$ upper bounds are provided by the following proposition, where we include a reduced version of that  for the case  $|x_n|<1$ given above.

%would be sometimes  useful. 
%Proposition 2.3
\begin{Prop}\label{prop2.4} \,
{\rm (i)} \, For all  admissible $\ga$ and $M>1$, there exists a constant $C_M$ such that for all $n\geq 1$ and $x\in \Z$, 
\[ 
Q^n_{\{0\}}(x,y)  \leq C_M \bigg(\frac{|x|^{\alpha-1} \vee 1}{L(x) n^2/c_n} \wedge \frac{L(x)}{|x|^{\alpha} \vee 1} \bigg) \frac{ |y|^{\alpha -1}\vee1}{L(y)} \quad (|y_n|<M); 
\]
%For $|\ga| <2-\a$, the restriction  $|y_n|<M$ may be removed.

 {\rm (ii)}  If $\ga =2-\a$, there exists a constant $C$ such that for all $x, y\in \Z$,
\beqn\label{eqL2.1}
Q^n_{\{0\}}(x,y) \leq  
C\bigg[\frac{a^\dagger(x)a^\dagger(-y)}{ n^{2}/c_n} +\frac{a^\dagger(-y)\{(x_n)_+\wedge 1\}  + a^\dagger (x)\{(y_n)_- \wedge 1\}}{n} \bigg]. 
\eeqn
 \end{Prop}
\v2
For $\ga=2-\a$, we shall provide a lower bound   which entails that
 {\it if either $C^+<\infty$ or  $F(x)$ varies regularly as  $x\to-\infty$,  then  (\ref{eqL2.1}) restricted to $|x|\vee|y|<Mc_n$ is exact}
(see Remark  \ref{r_P6.1}(i)).
In view of  $f^x(n) = Q^n(x,0)$ it follows from Proposition \ref{prop2.4}(i) that for all  $\gamma$,
$$f^x(n) \leq C  P[\, |X|> |x|\,] \qquad (n\geq 1,  x\in \Z).$$

(i) of Proposition \ref{prop2.4}  follows from Lemma \ref{L6.1} for $|x|<c_n$ and from  Lemma \ref{L6.2} for $|x|\geq c_n$. %Use  Theorem \ref{thm3}  in addition if $|\ga|=2-\a$. 
(ii)    follows from Theorems \ref{thm2} and \ref{thm3} in case   $ |x|\vee |y| \leq c_n$ with $xy\geq 0$, 
from Proposition \ref{prop6.2}(i)  in case $y\leq 0 \leq x$, from  Lemmas \ref{L6.1}(i)  and its dual (see Proposition \ref{prop6.2}(ii))  in case $0< (-x)\wedge y\leq c_n$,  from Lemmas \ref{L6.1} and \ref{L6.2} in case $|x|\wedge  |y|\leq c_n \leq |x|\vee |y|$ with $xy\geq 0$  ((\ref{eqL2.1}) is much weaker in this case) and from  the bound  $p^n(x) \leq C/c_n$  in case $|x|\wedge  |y|\geq c_n$ with $xy\geq 0$ as well as in case $(-x)\wedge y\geq c_n$. 
\v2

 %One may reasonably  conjecture that   (\ref{eqL2.1}) restricted to $|x|\vee(y|<Mc_n$ is exact, namely for some positive constant $C$
%$$C^{-1} Q^n_{\{0\}}(x,y) \leq  \frac{a^\dagger(x)a^\dagger(-y)}{ n^{2}/c_n} +\frac{a^\dagger(-y)(x_n)_ +  + a^\dagger (x) (y_n)_-}{n}  \leq C Q^n_{\{0\}}(x,y),$$
%this being true if either $C^+<\infty$ or  $F(x)$ varies regularly as  $x\to-\infty$ (see Remark \ref{r_P6.1}).
 \v2
  %Rem2.3
\begin{rem}\label{rem2.2}\, When  the condition 2) (strong aperiodicity) is not assumed, the results stated above must be subjected to  some modifications.
Let $\nu\geq 1$ denote  the  period of the walk, which amounts to assume  (in addition to (\ref{f_hyp})) that  $p^{\nu n}(0) >0$  and $p^{\nu n+j }(0) =0$ ($1\leq j < \nu)$ for all sufficiently large  $n$.
 Then  the process  $ \tilde S_n:= S_{\nu n}/\nu$, $n=0, 1,\ldots$ is a strongly aperiodic walk   on $\Z$ such that  
 $1-E[e^{i \th  \tilde S_1}] = 1-[\phi(\th/\nu)]^\nu \sim \nu^{1-\a} L(1/\th)\psi(\th)$,
 % where $\tilde L(x) = \nu L(x)$ and $\tilde \psi(\th)=\psi(\th/\nu) = E[e^{i\th Y_{1/\la^\a}}]$, 
 and $S^x_k\neq 0$ a.s.  for all combinations of  $x\in \nu \Z$ and $k\notin \nu\Z$; hence     the results on $Q^n_{\{0\}}$ restricted to $( \nu\Z)\times ( \nu\Z)$ follow  from  results of the aperiodic walks and the extension to $\Z$ is then readily performed by using $E[a(S_1^x)]=a^\dagger(x)$ if it concerns Theorems \ref{thm1} and \ref{thm2} and Proposition \ref{prop1};  this in particular results in
 $$f^0(n) = \nu \1( n\in \nu \Z)  f_*(n)\{1+o(1)\} \qquad (n\to\infty),$$
  where $ f_*(n) =  \k_{\a,\ga} c_n/n^2 $
 (with  $c_n$ satisfying (\ref{NC})),   and  for $x, y$ such that $p^n(y-x)>0$, 
$$Q^n_{\{0\}}(x,y) \sim \nu a^\dagger(x) f_*(n) a^\dagger(-y)   \quad \mbox{under} \quad |x_n|<\!< \La_n(x)  \; \mbox{and} \; |y_n| <\!< \La_n(-y)$$
 (note that $\tilde a(x/\nu) =a(x)$, $x\in \nu \Z$  for the last relation), and analogously for the other ranges of $x, y$.
As for the extension of  Theorems \ref{thm3} and \ref{thm4} we must replace the renewal function $U_{{\rm d}}$ by  $\tilde U_{{\rm d}}$ say,  the corresponding one  for the walk $\tilde S$, but  the formulae therein are kept unchanged under this replacement except for obvious modifications as made above, since
 \beqn\label{U*} 
 \tilde U_{{\rm d}}(x) /\tilde U_{{\rm d}}(c_n) \sim U_{{\rm d}}(x) /U_{{\rm d}}(c_n) \quad \mbox{as}\quad 
 n\wedge x\to\infty
 \eeqn
(see Lemma \ref{L8.9}).    
  \end{rem}
 \v2\v2
 %2.4 Extension
{\bf 2.4. Extension to  an arbitrary finite set.}
\v2
Let $A$ be a finite subset of $\Z$. Suppose for simplicity that for some $M>1$
\beqn\label{HA}
G_A(x,y)>0\quad \mbox{if} \quad  |x|\wedge |y| >M, 
\eeqn
where for a non-empty  $B\subset \Z$, $G_B$ denotes the Green function for the walk killed on $B$:
\beqn\label{green}
G_B(x,y) = \sum_{n=0}^\infty Q_B^n(x,y).
\eeqn
Since   $p$ is carried by an unbounded set, (\ref{HA}) implies $G_A(x,y)>0$ for all $x, y$.   There exists 
 \beqn\label{def_u}
u_A(x) =\lim_{|y|\to\infty} G_A(x,y)
\eeqn
(since  $\sigma^2=\infty$) \cite[Theorem 30.1]{S}.   $u_A$ is positive and  harmonic for the killed walk:   $u_A(x)= \sum_{z\notin A}p(z-x)u_A(z)>0$ for all $x\in \Z$. 
Put  
\[
f^x_A(n) = P[\sigma^x_A=n].
\]
In order to obtain the asymptotic form of $Q^n_A(x,y)$ and $f^x_A(n)$  we may  simply replace   $a^\dagger(x)$ by $u_A(x)$ and  $a^\dagger(-y)$ by $u_{-A}(-y)$  on the RHS of formulae given  in  Theorems \ref{thm2} to \ref{thm4} and Corollaries \ref{cor1} and \ref{cor2}---and accordingly replace 
 $C^+$ by  $C_A^+ = \lim_{x\to\infty} u_A(x)$ (positive under (\ref{HA})) in Theorem \ref{thm4}(ii)---with the resulting formulae  valid  in the same range of variables. Thus by Corollary \ref{cor1}  we have
\beqn\label{f/A}
f_A^x(n) \sim\left\{ \begin{array}{ll}  u_A(x)f^0(n)\quad &\mbox{for} \quad |x_n| <\!< \La_n(x),\\
{\displaystyle  \mathfrak{f}^{\, x_n}(1)  } \quad &\mbox{uniformly for} \quad 1/M <|x_n| <M.  \;\,  \end{array}\right.
\eeqn
Similarly corresponding to  Theorem \ref{thm3} we have:

\begin{Thm}\label{thm5}
If $\ga=2-\a$,  for each $M>1$ uniformly  for $|x| <Mc_n$, as $n\to\infty$
\[ %beqn\label{eq_thm5}
Q^n_{\{0\}}(x,y) \sim  \left\{\begin{array}{lr}
{ \displaystyle  \bigg\{ u_A(x)f^0(n) +\frac{U_{{\rm d}}(x_+)\mathfrak{p}_1(-x_n)}{U_{{\rm d}}(c_n)n} \bigg\}u_{-A}(-y)   } \quad  & (-\La_n(-y) <\!< y_n <\!<1),
\\[4mm]
{\displaystyle   \frac{u_A(x)}{n}\mathfrak{f}^{\,-y_n}(1) + \frac{U_{{\rm d}}(x_+)K_1(y_n)}{U_{{\rm d}}(c_n)c_n}   }\quad\;\;  & (x_n \to 0, y_n \in [1/M, M]).
\end{array} \right.
\]
[For $-M<x_n<\!< \La_n(x)$ (resp. $\La_n(x)<\!< x_n<M$), the second (first) terms on the RHS are negligible as compared to the first (second).]
\end{Thm}
% If $|\ga|<2-\a$, then  for any $M>1$, uniformly for $|x_n| \vee |y_n|<M$, as $n\to\infty$
%\beqn\label{Q/A}
%Q^n_{A}(x,y) \sim \left\{\begin{array}{lr}
%{\displaystyle   u_A(x)f^{0}(n) u_{-A}(-y) } \quad &( |x_n| \vee |y_n| \to 0),\\[2mm]
%{\displaystyle \mathfrak{f}^{\,x_n}(1)u_{-A}(-y)/n} &  (y_n \to 0, |x_n|>1/M),
%\\[2mm]
%{\displaystyle  u_A(x)\, \mathfrak{f}^{\,-y_n}(1)/n } \quad &( x_n \to 0, |y_n|>1/M),\\[2mm]
%\mathfrak{p}^{\{0\}}_{1} (x_n,y_n)/c_n & (|x_n| \wedge |y_n|\geq 1/M),
%\end{array} \right.
%\eeqn
%and  if $\ga= 2-\a$, this holds under the constraint on $x, y$ given in Proposition \ref{prop1}.

\v2
Corollaries \ref{loc_clt} and \ref{cor-1} with $\sigma^x_A$ in place of  $\sigma_{\{0\}}^x$ remains in force without any modification. The next corollary would be  also worthy of note (see Lemma \ref{L7.3} for proof),  providing an example that exhibits the significance of the condition $C^+<\infty$. In case $\sigma^2<\infty$  a similar result holds when and only when $xy<0$ as $|x|\wedge |y|\to\infty$ \cite[(1.20]{U1dm_f}.
%the walk makes occasional very  long  jumps to the right and drifts towards the left  between them.  The result given below may be understood to  reflect   this tendency of drifting to the left that becomes notably  strong  when  $C^+<\infty$.

%Cor5
\begin{Cor}\label{cor5}\,  Let  $\ga = -2+\a$. Then for any finite subset  $B\subset \Z$ such that $A\subset B$,  as $n\wedge x \wedge y\to\infty$ under $y_n<M$,  $x_n < M\La_n(x)$
\[
P[\sigma^x_B<n\,|\, S_n^x=y, \sigma_{A}^x>n] = 
 \frac{\lim_{x\to\infty}E[u_A(S^x_{\sigma_B}); S^x_{\sigma_B}\notin A]}{C_A^+ + \mathfrak{p}_1(0)x_+n/c^2_n}\{1+o(1)\}. 
\]
\end{Cor}
\v2
  
  %Write $H^\infty_{A}\{a\} =  \lim_{|x|\to\infty}E[a(S^x_{\sigma_A})]$ and
  On taking   $y\in A$  in  the formula of Theorem \ref{thm5} one obtains  
  %Cor6
%\begin{Cor}\label{cor3}  \,
for all admissible  $\ga$
$$Q^n_A(x,y)=  P[\sigma_A^x =n, S^x_n =y]  \sim f_A^x(n)u_{-A}(-y)  \quad (y\in A)$$
%\end{Cor} 
 uniformly for $|x_n|<M$, as $n\to\infty$. 
From the definition of $G_A(x,y)$ it follows  that $u_A$ is the probability distribution of  the hitting place of $A$ by the dual walk \lq started at infinity'. 
By  (\ref{f/A})  it therefore  follows that
 \beqn\label{Kstn}
 \sum_{x\in A} \sum_{y\in A} Q^n_{A}(x,y) =\sum_{x\in A}f_A^x(n) \sim f^0(n).
 \eeqn

\v2
%Rem4
\begin{rem}\label{rem4C}\, As was mentioned in Introduction Kesten \cite{K} 
%and Kesten and Spitzer \ref{KS} 
obtained  asymptotic formulae of  $Q_A^n(x,y)$ with $x, y$ fixed for a large class of random walks on multidimensional lattices $\Z^d$: in particular Theorem 6a of \cite{K}   specialized to  one-dimensional  walk may  read in the present notation
$$\lim_{n\to\infty}Q_A^n(x,y) \Big/ \sum_{z\in A} \sum_{w\in A} Q^n_{A}(z,w)  = u_A(x)u_{-A}(-y) \quad (y\notin A),$$ 
valid whenever  the walk is  recurrent, strongly aperiodic and   having infinite variance.   When $\ga=0$ and $L\equiv 1$ relation (\ref{Kstn})  is also observed in \cite{K}.
\end{rem}

The rest of the paper is organized as follows. In Section 3 we  prove some preliminary lemmas that are fundamental for our proofs of  Theorems \ref{thm1} and \ref{thm2} and Proposition \ref{prop1} that are given in Section 4. 
The proof  of Theorem \ref{thm3}   is given  in Section 5.  In Section 6 some estimations of $Q^n_{\{0\}}(x,y)$ are made in case $xy <0$ and, for this purpose,  beyond the regime $|x|\vee |y|=O(n^{1/\a})$.  Propositions \ref{prop6.1} and \ref{prop6.2} given there  provide a lower and upper bound, respectively and Theorem \ref{thm4} and Proposition \ref{prop2.2} are   proved after them.   In Section 7 the results  are extended to those for an arbitrary finite set instead of the single point set $\{0\}$. In Section 8  we deal with   the limit stable process and present some properties of $\mathfrak{f}^{\, \xi}(t)$, $\mathfrak{p}_t^{\{0\}}(x,y)$ and $U_{{\rm d}}(x)$.
In  the Appendix  we provide
 (A) condition (\ref{f_hyp})  expressed  in terms of  the tails of $F$ and some related facts,
 (B)   the representation of   $a^\dagger(x), x\geq 0$ as an integral of  $a(x), x<0$ in case $\ga =2-\a$,   and  (C) 
\lq some bounds of escape probabilities' from the origin.

%Section3
\section{Preliminary lemmas. }

In this section we show that the conditional  law of   $S^x_n/c_n$   given $\sigma^x_{\{0\}}>n$ converges 
as $n\to\infty$  uniformly for  $|x_n| <\!< \La_n(x)$---extending  Belkin's result (\ref{Belkin})  for the special case $x=0$---
 under a certain assumption if $\ga=|2-\a|$ that will be  removed in the next section.
%The following result is given in \cite{B2} (see (\ref{eq8.2}), (\ref{kapp}) for the identification of the expression of the constant  on the RHS of the formula of it). Its proof,  omitted in 
%\cite{B2}, is found in \cite{Uattrc}. 
%Lem3.1
\begin{lem}\label{lem3.1}   
\, Put  $\kappa^a_{\a,\ga,\pm} = -\Ga(1-\a) \pi^{-1} \sin [\tst12\pi (\a\pm\ga)]$. Then 
$\kappa^a_{\a,\ga,\pm} >0$ if $|\ga| <2-\a$,  and  $\kappa^a_{\a,\ga,\pm}  = 0$ or $1/\Ga(\a)$  according as $\pm \ga >0$ or $<0$  if $|\ga| =2-\a$; and
\v2
{\rm (i)}\qquad \;
${\displaystyle    \lim_{x\to \pm\infty} \frac{ a(x) L(x)}{|x|^{\a-1}} = \kappa^a_{\a,\ga,\pm},}$
\v2
{\rm (ii)} \qquad ${\displaystyle   \lim_{x\to \pm\infty}  \{a(x+1)-a(x)\}|x|^{2-\a} L(x) = \pm (\a-1)\kappa^a_{\a,\ga,\pm}. }$
\v2\n
% in particular
%if $\ga=2-\a$, then 
%\[a(x) = \left\{ 
%\begin{array} {ll} \Ga(\a)^{-1} |x|^{\a-1}/L(x) \{1+o(1)\} \quad  & \mbox{as}\quad  x\to -\infty, 
%\\[1mm]
%o(x^{\a-1}/L(x))     &\mbox{as}\quad  x\to +\infty.
%\end{array} \right. \]
\end{lem}
\n
\pf\, (i)  is given in \cite{B2} (without proof). Here we give a proof, which is partly used in the proof of (ii).    The proof is based on 
\beqn\label{eq3.00} 
\int_0^\infty \left\{\begin{array}{c}1-\cos u\\
\sin u \end{array}\right\} \frac{du}{u^\a}
 =  \left\{\begin{array}{l}  -\Ga(1-\a) \sin\frac12 \pi\a, \\[1mm]
 \Ga(1-\a) \cos\frac12\pi\a,
 \end{array}\right.
 \eeqn 
 (cf. \cite[pp.10, 68]{E}, \cite[p.260]{WW}).
 In the representation 
 $a(x)= \frac1{2\pi}\int_{-\pi}^\pi (1-e^{ix \th})(1-\phi(\th))^{-1}d\th$ (valid at least if $EX=0$) we replace $1-\phi(\th)$ by $\psi(\th)L(1/\th)$,  its principal part about zero,  
 and compute the resulting integral. 
 Changing a variable we have
\beqn\label{L3.1}
\int_{-\pi}^\pi\frac{1-e^{ix\th}}{\psi(\th) L(1/\th)} d\th = |x|^{\a-1}\int_{-\pi|x|}^{\pi |x|}
\frac{1-\cos u \mp i\sin u}{\cos \frac12 \pi\ga +i(u/|u|) \sin \frac12 \pi\ga} \cdot \frac{|u|^{-\a}du}{L(x/u)}\quad (\pm = x/|x|),
\eeqn
which  an easy computation with the help of  (\ref{eq3.00}) shows  to be asymptotically equivalent   to
$$%\frac1{2\pi}\int_{\pi}^\pi\frac{1-e^{ix\th}}{c_\circ \psi(\th)}=
 \frac{-2\Ga(1-\a)  |x|^{\a-1}}{L(x)}\Big[\cos\frac{\pi\ga}{2} \sin \frac{\pi\a}{2} \pm \sin\frac{\pi\ga}{2} \cos\frac{\pi\a}{2} \Big] $$
 as $|x|\to \infty$. 
The combination of the  sine's and cosine's in the square brackets being equal to  $\sin[\frac12\pi(\a\pm \ga)]$ we find the equality  (i), provided  that the replacement mentioned at the beginning causes only a negligible term of the magnitude  $o(|x|^{\a-1})/L(x)$, but this   is
assured from the way of computation carried out above since  the integrand of the RHS integral in (\ref{L3.1}) is summable on $\R$.  

   For the proof of (ii) it suffices to show that
\beqn\label{L3.11}
\int_{-\pi}^\pi e^{ix \th}(1-e^{i\th})\bigg[\frac1{1-\phi(\th)}- \frac1{\psi(\th)L(1/\th)}\bigg]d\th 
=o\bigg(\frac{|x|^{\a-2}}{L(x)}\bigg),
\eeqn
for  $a(x+1)-a(x)= \frac1{2\pi}\int_{-\pi}^\pi e^{ix \th}(1-e^{i\th})(1-\phi(\th))^{-1}d\th$ and  this integral with $1-\phi(\th)$ replaced by $\psi(\th)L(1/|\th|)$   is asymptotically equivalent to $ \pm (\a-1)\kappa^a_{\a,\ga,\pm} |x|^{\a-2}/L(x)$ as one sees by looking at the increment of the RHS of  (\ref{L3.1}).  Recall $L$ is so chosen  as  to be smooth enough. Because of  the fact that if $\psi(\th) L(1/\th) = \{1-\phi(\th)\}(1+\de(\th))$ then  $\de'(\th) \th \to 0$ ($\th \to 0$) (cf. (\ref{F/psi})), the relation 
(\ref{L3.11}) then can be shown in a standard way. 
\qed

\v2
 %Lem3.2
\begin{lem}\label{lem01}  \,{\rm (i)} There exists a constant $C$ such that 
$$|p^n(0) - p^n(x)|\leq C|x|/c_n^2 \qquad  (x\in \Z, n\geq 1).$$

{\rm (ii)}  For any $\e>0$,\,
${\displaystyle \frac1{a^\dagger(x) }\sum_{k>\e n}|p^k(0) - p^k(x)| \to 0 \quad (n\to\infty) \;\; \mbox{under} \;\; |x_n|  <\!<\La_n(x).}$
\end{lem}

Note that if the walk is left-continuous (right-continuous),   the assertion (ii) of Lemma \ref{lem01}
 is void for $x> 0$ ($x< 0$). The same remark applies to the subsequent lemmas.
  
\v2\n
\pf\, By (\ref{f_hyp}) $- \log |\phi(\th)| \sim (\cos \frac12 \ga\pi)|\th|^\a L(1/\th)$ ($\th \to 0$) and by the Fourier representation
\beqn \label{pL01_1}
p^n(0) - p^n(x) = \frac1{2\pi}\int_{-\pi}^\pi [\phi(\th)]^n(1-e^{-ix\th})d\th
\eeqn
as well as the assumed strong aperiodicity that entails $|\phi(\th)|<1$ ($0< |\th|\leq \pi$)
it follows that for any  positive $\e>0$ small enough 
\beqn\label{pL01_2}
|p^n(0) - p^n(x)| \leq  \frac{|x|}{\pi}\bigg[\int_{0}^\e e^{-\la n \th^\a L(1/\th)}  \th d\th + O(e^{-Cn})\bigg]
\eeqn
with  $\la := \frac12  \cos \frac12 \ga\pi$ and a positive constant $C=C(\e)$. Bringing in the variable $t= \th c_n$ we have
$$n\,\th^\a L(1/\th) =  t^\a L(c_n/t)/L(c_n) \{1+o(1)\}$$
because of  (\ref{NC}).  
Hence, on  choosing $\e >0$ so that  $u^{-\a/2}L(u)$ is decreasing on $u \geq 1/\e$, 
$$t^\a L(c_n/t)/L(c_n)  = t^{\a/2} [(t/c_n)^{\a/2} L(c_n/t)]/[ c_n^{-\a/2} L(c_n)] > t^{\a/2} \quad \mbox{for}\quad 1\leq  t\leq \e c_n.$$
 Hence
$$\int_{0}^\e e^{-\la n \th^\a L(1/\th)} \th  d\th \leq \int_{0}^{1/c_n} \th d\th  + \frac1{c_n^2}\int_{1}^{(\e c_n)\vee 1} e^{- t^{\a/2} \{\la+o(1)\}} tdt \leq \frac{C'}{c_n^2},$$
which combined with (\ref{pL01_2}) shows the first assertion (i).  The second formula (ii) follows from the first by an easy computation.
 \qed
\v2
 
%Lem3.3
\begin{lem}\label{lem021}  If $\ga\neq 2-\a$, then there exists a constant  $C$ such that
$$\sum_{k=0}^\infty  |p^k(0) - p^k(-x)|\leq Ca(x)\quad \mbox{for}\quad x\geq 0.$$
\end{lem}
\n\pf\, Note that by Lemma \ref{lem01}  the infinite series on the LHS  is convergent, so that we may consider only large values of $x$.   For each  $x$ large
enough let $j(x)$ be the smallest integer such that $c_{j(x)} \geq x$, which entails 
\beqn\label{pL021}
j(x) \sim  x^\a/L(x)
\eeqn
in view of the definition of $c_n$. Since  $\mathfrak{p}_{1}(-x/c_k)< \frac12 \mathfrak{p}_{1}(0)$ if $k/j(x)$ is small enough and $x$ large enough, we can choose, owing to the local limit theorem,
 positive constants $\de$ and  $N$ that do not depend on  $x$ such that  
\beqn\label{pL0211}
p^k(0)-p^k(-x)>0 \quad \mbox{for}\quad N\leq k\leq \de j(x).
\eeqn
Using Lemma \ref{lem01} we deduce that
$$\sum_{k>\de j(x)}|p^k(0) - p^k(-x)|  \leq C\de^{1-2/\a}xj(x)/c^2_{j(x)}.$$
Here $c_{j(x)}$ may be replaced by $x$  and because of (\ref{pL021}) the RHS is dominated by 
a constant multiple of $x^{\a-1}/L(x)$, hence of  $a(x)$  if $\ga\neq 2-\a$,  owing to Lemma \ref{lem3.1}. Because of (\ref{pL0211}) this concludes the proof. \qed

%Lem3.4
\begin{lem}\label{lem02}  Let  $\k= \k_{\a,\ga}/(1- \frac1{\a})$. If $\ga\neq 2-\a$, then 
%as $n\to\infty$
\v2
 {\rm (i)}\quad $P[\sigma^x_{\{0\}}>n] \sim  \k a^\dagger (x)c_n/n \quad(n\to\infty) $ \quad  under \;\; $0\leq x_n<\!< \La_n(x),$
\v2
 {\rm (ii)}\quad $\exists \; C>0$, \;   $P[\sigma^x_{\{0\}}>n] 
 \leq C(a^\dagger(x)P[\sigma^0_{\{0\}}>n] + x_n)$ \quad for \; $x\geq 0, \, n\geq 1$.
\end{lem}

Recall that for $\ga\neq 2-\a$, the condition  $0\leq x_n<\!< \La_n(x)$ is the same as $0\leq x_n<\!< 1$.

\v2\n
\pf\, For $x=0$ the assertion (i) reduces to  (\ref{S_f}). %proved by Belkin  \cite[Lemma 2.1]{B1} 
 For the proof of  (i) it therefore suffices to show 
\beqn\label{pL02}
P[\sigma^x_{\{0\}}>n] \sim a^\dagger(x) P[\sigma^0_{\{0\}}>n] \qquad (0< x_n <\!<\La_n(x)).
\eeqn
 We have the identity
$$p^n(-x) =f^x(n) + \sum_{k=1}^{n-1}p^k(-x)f^0(n-k) \qquad (n\geq 1, x\in \Z).$$
Noting  that  for $x=0$, the RHS can be written as  $\sum_{k=0}^{n-1}p^k(0)f^0(n-k)$,  we find 
$$f^x(n)= p^n(-x)-p^n(0) + \sum_{k=0}^{n-1}[p^k(0)-p^k(-x)]f^0(n-k) \qquad (n\geq 1, x\neq 0),$$
and by the usual method of generating function or by  direct computation   one  can easily derive
 \beqn\label{pL0212}
 P[\sigma^x_{\{0\}}\geq n]= \sum_{k=0}^{n-1} [p^k(0)-p^k(-x)]P[\sigma^0_{\{0\}}\geq n-k] \quad\;\; (n\geq 1, x\neq 0).
 \eeqn
We split the sum on the RHS at $k=\e n$  with a small $\e>0$.
Now suppose  $\ga\neq 2-\a$  and  apply   Lemmas \ref{lem01} (ii) and  \ref{lem021}  as well as  the  result for $x=0$ to see that  
the sum restricted to  $0\leq k\leq \e n$ is written  as
 \beqn\label{pL0213}
 \sum_{k\leq \e n} =  a^\dagger(x)P[\sigma^0_{\{0\}}\geq n] \{1 +o_\e(1)\}
 \eeqn
where $o_\e(1)$ is bounded on $x\geq 0$ and  tends to zero as $n\to\infty$ and $\e\to 0$ in this order uniformly for  $0\leq x_n <\!< \La_n(x)$, 
whereas for the other sum 
\beqn\label{pL0214}
\sum_{\e n< k < n} \leq \sum_{\e n< k < n} \frac{Cx}{c_k^2}\cdot\frac{c_{n-k}}{n-k} 
\sim \frac{Cx}{\e^{1/\a}c_n}\int_{\e}^1\frac{1}{u^{2/\a}(1-u)^{1-1/\a}}.
\eeqn
These  estimates together conclude (\ref{pL02})  as well as (ii). \qed
%On noting  that $xc_n^{-1}/ \{a^\dagger(x)P[\sigma^0_{\{0\}}>n]\} < C'x_n/\La_n(x)$, 

%\frac{c_\circ^{1/\a}\sin \pi/\a}{\pi \mathfrak{p}_1(0)}
\v2
%By Lemma \ref{lem01} 
%\[\sum_{k\geq \e n} |p^k(0)-p^k(-x)| \leq C \frac{|x|n}{B^2_n} \leq C''\frac{|x|}{\La_n(x)}\cdot a^\dagger(x),\]
 The same (rather simplified) proof as above verifies that even for $\ga =2-\a$ the statements (i) and (ii) of Lemma \ref{lem02} hold if  $x$ is confined in an arbitrarily  fixed finite interval. 
 %We also notify that there exists a constant $C$ such that  for all $\ga$
 %\beqn\label{b_HT}
 %P[\sigma^x_{\{0\}} >n] \leq C|x_n|\qquad   (x \neq 0, n\geq 1),
 %\eeqn
%which is immediate from (\ref{pL0212}) with the help of the estimate $P[\sigma^0_{\{0\}}>n] =O(c_n/n)$ and   Lemma \ref{lem01}(i).

We shall prove a local version of  formula  (\ref{S_f}) in the next section. Combined with it---with the above remark taken into account---the next lemma will remove the restriction $\ga\neq 2-\a$ from Lemma \ref{lem02}.  
%Lem3.5
\begin{lem}\label{lem04}   Suppose   $f^0(n) \leq C c_n/n^2$ for some constant $C$.  Then both (i) and (ii) of Lemma \ref{lem02} hold for all admissible $\ga$.
\end{lem}

\n
\pf\, Put $a(x,n) = \sum_{k=n}^\infty [p^k(0) -p^k(-x)]$. In the identity  (\ref {pL0212})  we split the sum on its RHS at  $\e n>0$,  which we suppose to be an integer for simplicity. Then summation by parts  transforms the sum  restricted on $k\leq \e n$ into
$$  a(x)P[\sigma^0_{\{0\}}\geq n]  - a(x,\e n+1) P[\sigma^0_{\{0\}}\geq (1-\e)n]+ \sum_{k=1}^{\e n} f^0(n-k)a(x,k)\qquad (x\neq 0).$$
By Lemma \ref{lem01}(i) we have  $a(x,k) \leq C|x|k/c_k^2$, which together with (\ref{S_f}) and the assumption of the lemma shows that the second and the third terms  are  dominated in absolute value by  a constant multiple of  
$|x_n| <\!< a(x)P[\sigma^0_{\{0\}}\geq n]$
under $|x_n|<\!<\La_n(x)$.  This combined with the bound (\ref{pL0214}) which  is valid also for $\ga=2-\a$ concludes the proof. \qed
 
\v2

Now we can prove  an extension of the aforementioned result (\ref{Belkin})  of \cite{B1}.
%Lem3.6
\begin{lem} \label{lem05}
If  $\ga\neq 2-\a$, then for any interval  $I \subset \R$, as  $n\to\infty$ 
 under $0\leq x <\!< c_n$ \beqn\label{eL04}
 P[ S^x_n/c_n \in I \,|\, \sigma^x_{\{0\}} > n] \, \longrightarrow \,  \int_I h(\xi)d\xi.
 \eeqn
%where $h$ is a probability density on $\R$ which is bounded and  continuous and whose characteristic function, denoted by  $\hat h(\th)$, is given by
%\beqn\label{eqL03}\hat h(\th) = 1- \psi(\th)\int_0^1 \xi^{\frac1{\a}-1}e^{  \psi(\th)(\xi -1)}d\xi.
%\eeqn
\end{lem}  
\n
\pf\,  This proof  will reduce the general case of $x$ to the case  $x=0$ that is treated in \cite{B1}.     As in  \cite{B1},  observe  the identity 
$$Q^n_{\{0\}} (x,y) = p^n(y-x)-\sum_{k=1}^{n-1}f^x(k)p^{n-k}(y)\qquad (y\in \Z),$$
and on putting $r^x_n= P[\sigma^x_{\{0\}}>n]$   carry out  summation by parts, which results in      
$$Q^n_{\{0\}} (x,y)=\sum_{k=0}^{n-1}r^x_k[p^{n-k}(y)-p^{n-k-1}(y)] + r^x_{n-1}\, p^0(y) +p^n(y-x)-p^n(y).$$
Then we  can  write the Fourier series 
$\sum_{y\in\Z} Q^n_{\{0\}} (x,y) e^{i \th y}$  as 
$$- \sum_{k=0}^{n-1}r^x_{k}[\phi^{n-k-1}(\th)-\phi^{n-k}(\th)]+r^x_{n-1}+(e^{i\th x}-1)\phi^n(\th) \quad ( -\pi<\th \leq \pi).$$
Writing  $\th_n =\th/c_n$ 
we accordingly   obtain 
\beq 
E[e^{i\th S^x_n/c_n}\,|\, \sigma^x_{\{0\}}>n] &=&\frac1{r^x_n} \sum_{y\in\Z} Q^n_{\{0\}} (x,y) e^{i \th_n y} \\
&=&  \frac{r^x_{n-1}}{r^x_n}-  \sum_{k=0}^{n-1}\frac{r^x_k}{r^x_n}[1-\phi(\th_n)] \phi^{n-k-1}(\th_n)
 +\frac{e^{i\th_n x}-1}{r^x_n} \phi^n(\th_n).
\eeq
Now let $\ga <2-\a$.  Then using Lemma \ref{lem02}(i) we observe that for each small  $\e>0$
$$\frac{r^x_k}{r^x_n} = \frac{r^0_k}{r^0_n}\{1+o(1)\}\;\;  (k>\e n)\;\;  \mbox{and}\;\; \frac{|e^{i\th_n x}-1|}{r^x_n} = |\th|\times o(1) \qquad  \mbox{for} \;\;  0\leq x_n <\!< \La_n(x),
$$
of which  the former relation shows that  $x$ can be replaced by  $0$  for $k>\e n$ in the sum on the right-most member of the above equalities and the latter that the last term  of it tends to zero. By Lemma \ref{lem02}(i) it also follows that  $r^x_{n-1}/r^x_n\to 1$ uniformly for $|x_n|<\!< \La_n(x)$. 
The reduction to the case $x=0$ follows if we can  show that  the sum over  $k\leq \e n$ tends to zero as $n\to \infty$ and $\e\to 0$ in turn. 
To this end observe  that by (i)
 and (ii) of Lemma \ref{lem02}
 $$\sum_{k\leq \e n} \frac{r^x_k}{r^x_n} \leq C\sum_{k\leq \e n} \frac{r^0_k}{r^0_n}+ C\sum_{k\leq \e n}  
\frac{|x|/c_k}{a^\dagger(x)c_n/n},$$
of   which  the first sum on the RHS  is bounded by  $C'\e^{1/\a}n$ and the second one is at most a constant multiple of 
$|x| n^2/[a^\dagger(x)c_n^2] =n |x_n|/\La_n(x)= o(n)$. 
Thus  the desired reduction is achieved,  for $1-\phi(\th_n) =O(1/n)$.  \qed

\v2
%Rem3.1
\begin{rem}\label{remL3.5}
If $x \geq 0$ is fixed,  (\ref{eL04}) holds for all  $\ga$ unless the walk is  left-continuous. 
In case $\ga=2-\a$,   if  $f^0(n) =O(c_n/n)$  is assumed in addition, then  (\ref{eL04})   holds  under $0 \leq x_n<\!< \La_n(x)$. (The same proof applies; see what is remarked before Lemma \ref{lem04}.)
\end{rem}

\v2

 In the sequel  we let $\k= \k_{\a,\ga}/(1- \frac1{\a})$   as in Lemma \ref{lem02}.

%Lem3.7
\begin{lem}\label{lem06} \,
\beqn\label{h/f}
h(\xi) =   \mathfrak{f}^{\, -\xi} (1)/\kappa \qquad (\xi\in \R).
\eeqn
\end{lem}

The identity above would be comprehended  in view of the identity $Q^n_{\{0\}}(0,x) = f^{-x}(n)$, it being expected that  if $x_n\to\xi \neq 0$, 
$f^{-x}(n) \sim
\mathfrak{f}^{-\xi}(1)/n$   and $Q^n_{\{0\}}(0,x) \sim  h(\xi)P[\sigma_{(-\infty,0]}^0>n]/c_n\sim
\kappa h(\xi)/n$,  and hence (\ref{h/f}). % (these are verified shortly but (\ref{h/f}) is used for the verification).  
 
    For the proof of (\ref{h/f}) one may suppose that $L\equiv 1$ and  the analysis under this assumption that  is carried out in \cite{Uattrc} yields  $f^{-x}(n) \sim \mathfrak{f}^{-\xi}(1)/n$ (as $x_n\to \xi$), hence (\ref{h/f}). We shall  see $Q^n_{\{0\}}(0,x) \sim \kappa h(\xi)/n$ shortly (see (\ref{a/h/f})) and this together with the functional limit theorem also leads to  (\ref{h/f}).
    We shall provide another more direct proof in Section 8  (Lemma \ref{L8.10}) by computing the Fourier transform of $\mathfrak{f}^{-\xi}(1)$. % by using an integral representation of $\mathfrak{f}^{-\xi}$. 

%Section4
  \section{ Proof of Theorems \ref{thm1}, \ref{thm2} and Proposition \ref{prop1}.} 
\v2
The proof is made by showing several lemmas and is based on Lemma \ref{lem05} that asserts the  convergence---with some uniformity in $x$---of the conditional law of $S^x_n/c_n $ given  
$\sigma^x_{\{0\}}>n$. Employing   Lemma \ref{lem02}(i) we can rephrase  this convergence  result  as follows:  

{\it If  $\ga\neq 2-\a$ and $\fa$ is a continuous function on $\R$ and $I$ is  a finite interval of the real line,  then   uniformly for  $|y_n|<M$ and $0\leq x_n <\!<  1$, as $n\to\infty$}
 \beqn\label{crL03}
\sum_{w: w/c_n\in I} Q^{n}_{\{0\}}(x,y+w)\fa(w/c_n) \sim \kappa a^\dagger(x)\frac{c_n}{n}\int_I h(y_n+\xi)\fa(\xi)d\xi.
 \eeqn  
Here and in the sequel  $M$ denotes  an arbitrarily fixed constant larger than  1. 
%Lem4.1
\begin{lem}\label{lem4.1} \, If  $\ga\neq 2-\a$, then  uniformly for  $ |y_n| \in [1/M, M]$, as $n\to\infty$
\beqn\label{eL07}
Q^n_{\{0\}}(x,y) \sim  a^\dagger(x)  \mathfrak{f}^{\, -y_n}(1)/n \qquad \mbox{under} \;\; 0\leq x_n <\!< \La_n(x).
\eeqn
The comment given to (\ref{eL04}) in Remark \ref{remL3.5} is applicable to this  relation; in particular
\beqn\label{eL071}
f^x(n) =Q^n_{\{0\}}(0,-x)\sim  \mathfrak{f}^{\, x_n}(1)/n
\eeqn
uniformly for $|x_n| \in [1/M, M]$ for every admissible $\ga$.
\end{lem}

As is noted previously the first assertion  of the lemma is void for $x>0$ if the walk is left-continuous when the
asymptotics is described  by means of   the stable meander. % similarly for the right continuous walk.
\v2
\n
\pf\,  (\ref{eL07}) is essentially  the local limit theorem corresponding to  (\ref{crL03}) in view of Lemma \ref{lem06}.  We are going to derive  the former  from the integral one 
with the help of  Gnedenko's theorem.
 The idea of the following proof is
borrowed from \cite{DW} (the proof of  Theorem 5 in it).  Taking  $m=\lfloor \e^{2\a} n\rfloor$ with a small $\e>0$  we decompose
$$Q^n_{\{0\}}(x,y) = \sum_{z\in \Z\setminus \{0\}}  Q^{n-m}_{\{0\}}(x,z) Q^m_{\{0\}}(z,y).$$
 Note  $c_m/c_{n-m}\sim \e^2/ (1-\e^{2\a})^{1/\a}$.  We apply  (\ref{crL03})   in the form  
 \beqn\label{pL3.4}
 n\sum_{|u|<\e c_{n-m}} Q^{n-m}_{\{0\}}(x,y-u)\frac{\fa(u/c_m)}{c_m} \sim \frac{\kappa a^\dagger(x)}{1-\e^{2\a}} \int_{|\xi|<\e} h(y_{n- m}-\xi) \frac{\fa(\xi/\tilde \e^2)}{\tilde\e^2}d\xi
 \eeqn
valid for each  $\e>0$ fixed, where $\tilde \e =\e/(1-\e^{2\a})^{1/2\a}\sim \e$ ($\e\to 0$). % and $C(\e)= (1-\e^{2\a})^{\frac1{\a}-1}$.

%The following estimates hold as $n\to\infty$ and $\e\to 0$:
Let $|y_n| \in [1/M, M]$.  It is easy to see  
$$\sup_{z: |z-y|<\e c_n} |Q^m_{\{0\}}(z,y) - p^m(y-z)| = p^m(y-z)|\times o_\e(1),$$
Since  $\e c_n/c_m \sim 1/\e$ and  $\mathfrak{p}_1(\pm 1/\e) =O(\e^{\a+1}) $ according to \cite{Zol} (cf. \cite[ Eq(14.34-35)]{Sk}),  we also have  for all sufficiently large $n$ 
$$\sup_{z: |z-y| \geq  \e c_{n-m}}Q^m_{\{0\}}(z,y) <  C \e^{\a+1}/c_m$$
which combined with the preceding bound yields
$$\bigg|Q^n_{\{0\}}(x,y) - \sum_{|z-y| <\e c_{n-m}} Q_{\{0\}}^{n-m}(x,z)p^m(y-z)\{1+o_\e(1)\} \bigg|\leq  C'\frac{P[\sigma^x_{\{0\}}>n-m] \e^{\a+1}}{c_m}.$$
On the LHS $p^m(y-z)$ may be replaced by $\mathfrak{p}_{1}((y-z)/c_m)/c_m$   (whenever $\e$ is fixed) whereas  the RHS  is dominated by a constant multiple of $a^\dagger(x)[c_n /nc_m] \e^{\a+1}\sim a^\dagger(x) \e^{\a-1}/n$. Hence,  after a change of variable 
\[
nQ^n_{\{0\}}(x,y) =  n  \sum_{|u| <\e c_{n-m}} Q_{\{0\}}^{n-m}(x,y-u) \frac{\mathfrak{p}_{1}(u/c_m)}{c_m} \{1+o_\e(1)\} + a^\dagger(x) \times o_\e(1).
\]
Now on letting $n\to\infty$ and $\e\to0$ in this order  (\ref{pL3.4})  shows that   the RHS  can be written as  
\beqn\label{a/h/f}
\kappa a^\dagger(x)\bigg\{\int_{|\xi|<\e} h(y_n-\xi)  \frac{\mathfrak{p}_{1}(\xi /\e^2)}{\e^2}d\xi +o_\e(1)\bigg\} =a^\dagger(x) \{ \kappa h(y_n) +o_\e(1)\}
\eeqn
and   we find the formula (\ref{eL07}) owing to  Lemma \ref{lem06}. \qed

\v2 

%Lem4.2
\begin{lem}\label{lem4.2}  Uniformly for $ |x_n| \in [1/M, M]$ and $ |y_n| \in [1/M, M]$,
$$Q^n_{\{0\}}(x,y) =\frac1{c_n}\{ \mathfrak{p}^{\{0\}}_{1}(x_n,y_n) +o(1)\}.$$
[ Note that if $\ga=|2-\a|$, then $\mathfrak{p}^{\{0\}}_{1}(x_n,y_n) =0$ for $\ga x>0$ and $\ga y<0$.] 
\end{lem}
\n
\pf\,  Let $|x_n|, |y_n| \in [1/M,M]$.  By (\ref{eL071}) $f^x(k) \sim  \mathfrak{f}^{\, x_k}(1)/k $ for $k \asymp c_n$ and we see
\beq
Q^n_{\{0\}}(x,y) &=& p^n(y-x) -\sum_{k=0}^{n-1} f^x(k) p^{n-k}(y)
\\
&=&\frac{ \mathfrak{p}_{1}(y_n-x_n) }{c_n}\{1+o(1)\}-  \sum_{k=\e n}^{(1-\e)n} \frac{\mathfrak{f}^{\, x_k}(1)}{k}\cdot \frac{\mathfrak{p}_{1}(y_{n-k})}{c_{n-k}}\{1+o(1)\} \\
&& \, + R_1(\e) + R_2(\e),
\eeq
where $R_1(\e)$ and  $R_2(\e)$ denote the sums of $f^x(k) p^{n-k}(y)$ over $1\leq k < \e n$ and $(1-\e)n < k\leq n$, respectively. Since $c_k/c_n \sim (k/n)^{1/\a}$ uniformly for $\e n <k<(1-\e)n$, on employing  scaling property  $\mathfrak{p}_t(\xi) = \mathfrak{p}_1(\xi/t^{1/\a})/t^{1/\a} $ and similar one of $\mathfrak{f}^{\,\xi}(t)$  the last sum above  is asymptotically equivalent to
$$\frac1{nc_n} \sum_{k=\e n}^{(1-\e)n} \mathfrak{f}^{\, x_n}(k/n)\mathfrak{p}_{1-k/n}(y_{n}) \sim \frac1{c_n}\int_\e^{1-\e} \mathfrak{f}^{\, x_n}(t)\mathfrak{p}_{1-t}(y_{n})dt. $$
By the local limit theorem (\ref{llt})  and (\ref{eL071}) we see 
$$R_1(\e) \leq C\frac1{c_n}P[\sigma_{\{0\}}^x <\e n] \leq C\frac1{c_n}P[\sigma_{(-\infty,0]}^x <\e n] = o_\e(1)\times \frac1{c_n}, \quad \mbox{and}$$
$$ R_2(\e) \leq \bigg(\sup_{(1-\e)n<k\leq n} f^x(k)\bigg) \sum_{k=0}^{\e n}p^k(y) \leq C\frac1{n} \sum_{k=1}^{\e n} \frac1{c_k} < C' \e^{1-1/\a}/c_n.$$ 
and similarly  $(\int_0^\e+ \int_{1-\e}^1) \mathfrak{f}^{\, x_n}(t)\mathfrak{p}_{1-t}(y_{n})dt \to 0$ as $\e\to 0$. Collecting these estimates we obtain the relation of the lemma since $\mathfrak{p}^{\{0\}}_{1}(x_n,y_n) = \mathfrak{p}_{1}(y_{n}-x_n) - \int_0^1 \mathfrak{f}^{\, x_n}(t)\mathfrak{p}_{1-t}(y_{n})dt$. \qed  

%Lem4.3
\begin{lem}\label{lem4.3} \, If $\ga\neq 2-\a$, then uniformly for $0\leq x_n <\!< \La_n(x)$ and $ y_n \in [1/M, M]$,   as $n\to\infty$ and $\e\downarrow 0$ 
$$\frac{\sup \{ Q^n_{\{0\}}(x,z) :\,|z_n|<\e  \;\mbox{or} \; |z_n|>1/\e\}}{Q^n_{\{0\}}(x,y)} \;\longrightarrow\; 0.$$
The same comment as given to (\ref{eL04}) in Remark \ref{remL3.5} applies  to this relation.
\end{lem}
\n
\pf\, 
 For simplicity   suppose $m:=n/2$ to be an integer and decompose 
 \beqn\label{pL09}
  Q^n_{\{0\}}(x,z) =  \sum_w Q^m_{\{0\}}(x,w)Q^m_{\{0\}}(w,z).
 \eeqn 
Observe that  $\sup_{|w|\leq c_n/2\e}Q^m_{\{0\}}(w,z) \leq c_m^{-1} [\sup_{|u|>1/2\e}\mathfrak{p}_1(u) +o(1)]$ \,  uniformly for  $|z_n|>1/\e$. If   $\ga\neq 2-\a$, then    $\sum_{|w| > c_n/2\e}Q^m_{\{0\}}(x,w)  /P[\sigma^x_{\{0\}}>n] \to 0$  as $n\to\infty$ and $\e\downarrow 0$ uniformly for $0\leq x_n <\!< \La_n(x)$ in view of Lemma \ref{lem05}, and hence by  Lemma \ref{lem02} (ii)
   \beqn\label{pL091}
  \sup_{|z_n|>1/\e}  Q^n_{\{0\}}(x,z)  \leq P[\sigma^x_{\{0\}}>n] c_n^{-1}\times o_\e(1)
= a^\dagger(x)n^{-1}\times o_\e(1),
\eeqn
where $o_\e(1)\to 0$ as $n\to \infty$ and $\e\downarrow 0$ interchangeably.      As for  the case $|z_n|<\e$, further making decomposition
$Q^m_{\{0\}}(w,z)
 =\sum_{w'\neq 0} Q^{m/2}_{\{0\}}(w,w')Q^{m/2}_{\{0\}}(w',z)$   we deduce,
 $$
 \sup_w Q^m_{\{0\}}(w,z) \leq \sup_{w,w'} Q^{m/2}_{\{0\}}(w,w')\sum_{u\neq 0}Q^{m/2}_{\{0\}}(u,z) \leq CP[\sigma^{-z}_{\{0\}}> m/2]/c_n = c_n^{-1}\times o_\e(1)$$
 so that 
 $Q^n_{\{0\}}(x,z)  \leq C' P[\sigma^x_{\{0\}} >m]c_n^{-1}\times o_\e(1) =[ a^\dagger(x)/ n]\times o_\e(1)$, which 
 together with (\ref{pL091}) concludes the proof,  since $ Q^n_{\{0\}}(x,y) \asymp a^\dagger(x)/n$ for $|y_n| \in [1/M, M]$ owing to Lemma \ref{lem4.1}. 
 \qed

%Lem4.4
\begin{lem}\label{lem4.4}   Uniformly under $|x_n| <\!< \La_n(x)$ and $|y_n|<\!< \La_n(-y)$, as $n\to\infty$ 
$$Q^n_{\{0\}}(x,y) \sim a^\dagger(x) f^0(n) a(-y).$$
%If $\ga=2-\a$,  then this holds under $0\leq x <\!< \La_n$ and $|y|\<\!< \La_n(y)$. 
 \end{lem}
\n
\pf\,  Pick a constant  $\la$ such that $0<\la<1$ and put $m =\lfloor \la n\rfloor$.   Consider the decomposition
$$Q^n_{\{0\}}(x,y)  = \sum_{z\neq 0} Q^m_{\{0\}}(x,z) Q^{n-m}_{\{0\}}(z,y). $$

Let $\ga\neq |2-\a|$. Then by Lemma \ref{lem05}
$$\frac1{P[\sigma^x_{\{0\}} >m]}\bigg(\sum_{|z_n|< \e } +\sum_{|z_n| >\e^{-1}}\bigg)Q^m_{\{0\}}(x,z) \, \longrightarrow\, \bigg(\int_{|\xi|<\e/ \la^{1/\a}}+\int_{|\xi|> \e^{-1}/\la^{1/\a}}\bigg)h(\xi)d\xi, 
$$
which together with Lemma \ref{lem4.3} applied to  $Q^{n-m}_{\{0\}}(z,y)= Q^{n-m}_{\{0\}}(-y,-z)$ shows  that the sum over $z$ subject to $|z_n|<\e$ or $|z_n|>\e^{-1}$ is negligible as $\e\to 0$. As for the  sum over $z$: $\e\leq |z_n|\leq \e^{-1}$ we see that its summands are expressed as 
$$a^\dagger(x)a^\dagger(-y)\mathfrak{f}^{\,-z_m}(1)\mathfrak{f}^{\, z_{n-m}}(1)/[m(n-m)]\{1+o(1)\}$$ 
owing to
Lemma \ref{lem4.1}, and, on noting that $z_m$ and $z_{n-m}$ can be replaced by  $z_n/\la^{1/\a}$ and
$z_n/(1-\la)^{1/\a}$, respectively and  that $\mathfrak{f}^{\, \xi/t^{1/\a}}(1) = \mathfrak{f}^{\,\xi}(t)t$, the sum itself is given as
$$\sum_{z:\e\leq  |z_n| \leq \e^{-1}} = \frac{a^\dagger(x)a^\dagger(-y) \la(1-\la)c_n}{m(n-m)}\int_{\e <|\xi|<1/\e} \mathfrak{f}^{\, -\xi}(\la)\mathfrak{f}^{\, \xi}(1-\la)d\xi \{1+o(1)\}.$$
Since $\e$ may be arbitrarily small, we conclude
\beqn\label{Q}
Q^n_{\{0\}}(x,y)  \sim \frac{a^\dagger(x)a^\dagger(-y) c_n}{n^2}\int_{-\infty}^{\infty}  \mathfrak{f}^{\, -\xi}(\la)\mathfrak{f}^{\, \xi}(1-\la)d\xi
\eeqn
(under  $|x_n|<\!<\La_n(x), |y_n| <\!< \La_n(-y)$).  Taking $x=y=0$ this becomes
\beqn\label{f*f}
f^0(n) \sim  \frac{c_n}{n^2} \int_{-\infty}^{\infty}  \mathfrak{f}^{\, -\xi}(\la)\mathfrak{f}^{\, \xi}(1-\la)d\xi.
\eeqn
Hence the formula of the lemma follows.

In view of Remark \ref{remL3.5} (and the comments given at the end of the succeeding lemmas)  (\ref{Q}) with   $x, y$ fixed---in particular  (\ref{f*f})---holds for all  $\ga$, which in turn 
shows that  the assertion  of the lemma is also valid for $\ga =|2-\a|$. \qed
 
\v2
{\it Proof of Theorem \ref{thm1}.} Comparing (\ref{f*f}) with  $P[\sigma^0_{\{0\}}>n] \sim  \kappa\, c_n/n$  ($\kappa= \kappa_{\ga,\a}/(1-\frac1{\a})$) one obtains the identity
\beqn
\int_{-\infty}^{\infty}  \mathfrak{f}^{\, -\xi}(\la)\mathfrak{f}^{\, \xi}(1-\la)d\xi =\kappa_{\a, \ga},
\eeqn
which together with (\ref{f*f}) shows Theorem \ref{thm1}.

\begin{figure}[t]
 \begin{center}
\begin{picture}(405,205)(0,-10)

%(Left)
\put(85,40){{\footnotesize $\gamma=2-\alpha$}}

\put(25,120){\vector (1,0){150}}
\put(167,125){{\footnotesize  {\sc x}}}
\put(100,60){\vector (0,1){126}}
\put(104,178){{\footnotesize {\sc y}}}

\thicklines
\put(40,145){\line (1,0){118}}
\put(40,95){\line (1,0){118}}
\put(75,70){\line(0,1){100}}
\put(125,70){\line(0,1){100}}

\put(130,99){{\footnotesize $ y_n = -1/M$}}
\put(128,72){{\footnotesize $x_n = 1/M$}}

%vertical dashed line

\put(106, 70){\line(0,1){4}}
\put(106, 76){\line(0,1){4}}
\put(106, 82){\line(0,1){4}}
\put(106, 88){\line(0,1){4}}
\put(106, 94){\line(0,1){4}}
\put(106, 100){\line(0,1){4}}
\put(106, 106){\line(0,1){4}}
\put(106, 112){\line(0,1){4}}
\put(106, 118){\line(0,1){4}}
\put(106, 124){\line(0,1){4}}
\put(106, 130){\line(0,1){4}}
\put(106, 135){\line(0,1){3}}
\put(106, 140){\line(0,1){3}}
\put(106, 146){\line(0,1){4}}
\put(106, 152){\line(0,1){4}}
\put(106, 157){\line(0,1){3}}
\put(106, 163){\line(0,1){3}}
\put(106, 168){\line(0,1){2}}

%horizontal dashes line

\put(40,112){\line(1,0){2}}
\put(44,112){\line(1,0){4}}
\put(50,112){\line(1,0){4}}
\put(56,112){\line(1,0){4}}
\put(62,112){\line(1,0){4}}
\put(68,112){\line(1,0){4}}
\put(74,112){\line(1,0){4}}
\put(80,112){\line(1,0){4}}
\put(86,112){\line(1,0){4}}
\put(92,112){\line(1,0){4}}
\put(98,112){\line(1,0){4}}
\put(104,112){\line(1,0){4}}
\put(110,112){\line(1,0){4}}
\put(116,112){\line(1,0){4}}
\put(122,112){\line(1,0){4}}
\put(128,112){\line(1,0){4}}
\put(134,112){\line(1,0){4}}
\put(140,112){\line(1,0){4}}
\put(146,112){\line(1,0){4}}
\put(152,112){\line(1,0){4}}

\thinlines
%(1)
\put(101,112){\line(0,1){56}}
\put(102,112){\line(0,1){56}}
\put(103,112){\line(0,1){56}}
\put(104,112){\line(0,1){56}}
\put(105,112){\line(0,1){56}}

%(2)
\put(101,95){\line(0,-1){24}}
\put(102,95){\line(0,-1){24}}
\put(103,95){\line(0,-1){24}}
\put(104,95){\line(0,-1){24}}
\put(105,95){\line(0,-1){24}}

\put(40,114){$....................$}
\put(42,117){$...................$}
\put(42,122){$...................$}
\put(40,125){$....................$}
\put(42,128){$...................$}
\put(40,131){$....................$}
\put(40,134){$....................$}
\put(42,137){$...................$}
\put(40,140){$....................$}
\put(42,143){$...................$}
\put(42,146){$...................$}
\put(42,149){$...................$}
\put(40,152){$....................$}
\put(42,155){$...................$}
\put(40,158){$....................$}
\put(42,161){$...................$}
\put(40,164){$....................$}
\put(42,167){$...................$}

\put(127,114){$..........$}
\put(125,117){$.........$}
\put(127,122){$.........$}
\put(125,125){$..........$}
\put(127,128){$.........$}
\put(127,131){$.........$}
\put(125,134){$..........$}
\put(127,137){$.........$}
\put(127,140){$.........$}
\put(125,143){$..........$}
\put(127,146){$.........$}
\put(127,149){$.........$}
\put(125,152){$..........$}
\put(127,155){$.........$}
\put(125,158){$..........$}
\put(127,161){$.........$}
\put(125,164){$..........$}
\put(127,167){$.........$}

\put(40,92){$....................$}
\put(42,89){$...................$}
\put(40,86){$....................$}
\put(40,83){$....................$}
\put(40,80){$....................$}
\put(42,77){$...................$}
\put(40,74){$....................$}
\put(40,71){$....................$}

%%%%
%%%

%(Right)
\put(274,40){{\footnotesize  $\gamma= -2+\alpha$}}

\put(225,120){\vector (1,0){152}}
\put(369,110){{\footnotesize {\sc x}}}
\put(300,60){\vector (0,1){126}}
\put(304,178){{\footnotesize {\sc y}}}
\put(330,136){{\footnotesize $y_n = 1/M$}}
\put(220,162){{\footnotesize $x_n = - 1/M$}}

\thicklines
\put(240,145){\line (1,0){118}}
\put(240,95){\line (1,0){118}}
\put(275,70){\line(0,1){100}}
\put(325,70){\line(0,1){100}}

%vertical dashed line
%\thinlines
%\put(294, 64){\line(0,1){2}}
\put(294, 70){\line(0,1){4}}
\put(294, 76){\line(0,1){4}}
\put(294, 82){\line(0,1){4}}
\put(294, 88){\line(0,1){4}}
\put(294, 94){\line(0,1){4}}
\put(294, 100){\line(0,1){4}}
\put(294, 106){\line(0,1){4}}
\put(294, 112){\line(0,1){4}}
\put(294, 118){\line(0,1){4}}
\put(294, 124){\line(0,1){4}}
\put(294, 130){\line(0,1){4}}
\put(294, 136){\line(0,1){4}}
\put(294, 142){\line(0,1){4}}
\put(294, 148){\line(0,1){4}}
\put(294, 154){\line(0,1){4}}
\put(294, 160){\line(0,1){4}}
\put(294, 166){\line(0,1){4}}

%horisontal dashed line

\put(240,128){\line(1,0){2}}
\put(244,128){\line(1,0){4}}
\put(250,128){\line(1,0){4}}
\put(256,128){\line(1,0){4}}
\put(262,128){\line(1,0){4}}
\put(268,128){\line(1,0){4}}
\put(274,128){\line(1,0){4}}
\put(280,128){\line(1,0){4}}
\put(286,128){\line(1,0){4}}
\put(292,128){\line(1,0){4}}
\put(298,128){\line(1,0){4}}
\put(304,128){\line(1,0){4}}
\put(310,128){\line(1,0){4}}
\put(316,128){\line(1,0){4}}
\put(322,128){\line(1,0){4}}
\put(328,128){\line(1,0){4}}
\put(334,128){\line(1,0){4}}
\put(340,128){\line(1,0){4}}
\put(346,128){\line(1,0){4}}
\put(352,128){\line(1,0){4}}

\thinlines
%(1)
\put(302,145){\line(0,1){20}}
\put(304,145){\line(0,1){20}}
\put(306,145){\line(0,1){20}}
\put(308,145){\line(0,1){20}}
\put(310,145){\line(0,1){20}}
\put(312,145){\line(0,1){20}}
\put(314,145){\line(0,1){20}}
\put(316,145){\line(0,1){20}}
\put(318,145){\line(0,1){20}}
\put(320,145){\line(0,1){20}}
\put(322,145){\line(0,1){20}}
\put(324,145){\line(0,1){20}}

%(2)
\put(302,70){\line(0,1){24}}
\put(304,70){\line(0,1){24}}
\put(306,70){\line(0,1){24}}
\put(308,70){\line(0,1){24}}
\put(310,70){\line(0,1){24}}
\put(312,70){\line(0,1){24}}
\put(314,70){\line(0,1){24}}
\put(316,70){\line(0,1){24}}
\put(318,70){\line(0,1){24}}
\put(320,70){\line(0,1){24}}
\put(322,70){\line(0,1){24}}
\put(324,70){\line(0,1){24}}

%(3)
\put(327,120){\line(0,-1){24}}
\put(329,120){\line(0,-1){24}}
\put(331,120){\line(0,-1){24}}
\put(333,120){\line(0,-1){24}}
\put(335,120){\line(0,-1){24}}
\put(337,120){\line(0,-1){24}}
\put(339,120){\line(0,-1){24}}
\put(341,120){\line(0,-1){24}}
\put(343,120){\line(0,-1){24}}
\put(345,120){\line(0,-1){24}}
\put(347,120){\line(0,-1){24}}
\put(349,120){\line(0,-1){24}}
\put(351,120){\line(0,-1){24}}
\put(353,120){\line(0,-1){24}}
\put(355,120){\line(0,-1){24}}
\put(357,120){\line(0,-1){24}}
\put(359,120){\line(0,-1){24}}

%(4)
\put(273,95){\line(0,-1){24}}
\put(271,95){\line(0,-1){24}}
\put(269,95){\line(0,-1){24}}
\put(267,95){\line(0,-1){24}}
\put(265,95){\line(0,-1){24}}
\put(263,95){\line(0,-1){24}}
\put(261,95){\line(0,-1){24}}
\put(259,95){\line(0,-1){24}}
\put(257,95){\line(0,-1){24}}
\put(255,95){\line(0,-1){24}}
\put(253,95){\line(0,-1){24}}
\put(251,95){\line(0,-1){24}}
\put(249,95){\line(0,-1){24}}
\put(247,95){\line(0,-1){24}}
\put(245,95){\line(0,-1){24}}
\put(243,95){\line(0,-1){24}}
\put(241,95){\line(0,-1){24}}

%\put(295,116){$....................$}
%\put(295, 69){$....................$}
\put(295, 72){$....................$}
\put(296, 75){$...................$}
\put(295, 78){$....................$}
\put(295, 81){$....................$}
\put(296, 84){$...................$}
\put(295, 87){$....................$}
\put(295, 90){$....................$}
\put(296, 93){$...................$}
\put(295, 96){$....................$}
\put(296, 99){$...................$}
\put(295, 102){$....................$}
\put(295, 105){$....................$}
\put(295, 108){$....................$}
\put(296,111){$...................$}
\put(295,114){$....................$}
\put(296,117){$...................$}
\put(295,122){$....................$}
\put(296,125){$...................$}

%Overlapping(1)
%\put(304, 99){$........$}
%\put(302, 90){$........$}
%\put(297,107){$...................$}
%\put(295,112){$....................$}
%\put(297,117){$...................$}

%\put(301, 75){$.......$}
%\put(303, 76){$.......$}
%\put(301, 80){$.......$}
%\put(303, 81){$.......$}
%\put(301, 85){$.......$}
%\put(303, 86){$.......$}

\put(241,125){$..........$}
\put(239,122){$..........$}
\put(239,118){$..........$}
\put(241,115){$..........$}
\put(239,112){$..........$}
\put(241,109){$..........$}
\put(239,106){$..........$}
\put(241,103){$..........$}
\put(241,100){$..........$}
\put(239,97){$..........$}
\put(241,94){$..........$}
\put(241,91){$..........$}
\put(239,88){$..........$}
\put(239,85){$..........$}
\put(241,82){$..........$}
\put(239,79){$..........$}
\put(239,76){$..........$}
\put(241,73){$..........$}

%\put(295,169){$....................$}
\put(295,165){$....................$}
\put(296,162){$...................$}
\put(295,159){$....................$}
\put(295,156){$....................$}
\put(296,153){$...................$}
\put(295,150){$....................$}
\put(295,147){$....................$}

\put(12,15){{\footnotesize {\sc Figure 2}:  The dotted regions  indicate the range of validity for (\ref{eq_thm2})  in the extreme  }}
\put(12,3){{\footnotesize  cases $\ga =\pm 2-\a$: they are symmetrical to each other about the diagonal $x-y =0$ }}
\put(12,-9){{\footnotesize while
each one is symmetrical to itself about the diagonal $x+y=0$}}
\end{picture}
\end{center}
% \caption{{\bf u} is  the  point of intersection where the half line  $\{(\eta,a): \eta>0\}$ meets
%  the  line passing through $\x$ and   tangential to  the circle   $\partial U(a)$.}

\vspace*{0cm}
\end{figure}
\v2\v2

{\it Proofs of Theorem \ref{thm2} and Proposition \ref{prop1}.}  What are asserted are virtually  involved  in Lemmas \ref{lem4.1}, \ref{lem4.2} and \ref{lem4.4} because of duality. This is immediate for Theorem \ref{thm2}. Proposition \ref{prop1} is verified as follows. Denote by $D_+$ and $D_-$ the diagonals $x-y=0$ and $x+y=0$, respectively.   The range of $(x,y)$ for validity of the formulae in  (\ref{eq_thm2})  is symmetrical about $D_-$ to itself and  about $D_+$ to that for the dual walk.
 If $\ga=-2+\a$, then by Lemma \ref{lem4.1} (\ref{eq_thm2}) holds for  $0\leq x <\!< c_n, |y_n|\in [1/M,M]$  and  their reflections about $D_-$  (the darkened region in  {\sc Figure 2} right). To this region  Lemmas \ref{lem4.2} and \ref{lem4.4} add the regions 
$x, |y_n|\in [1/M,M]$ and $0\leq x, -y <\!< c_n $, respectively. As a consequence it follows that
for $\ga=-2+\a$,  the range of validity contains 
\beqn\label{pT2}
 0\leq x_n < M, \quad y_n\in [-M, 0] \cup [1/M,M].
 \eeqn
Now let $\ga=2-\a$.  First note that  the assumption on $f^0(n)$ in Remark \ref{remL3.5} is valid (since Theorem \ref{thm1} has been shown) so that  the formula of \ref{lem4.1} is available.    Then by Lemma \ref{lem4.4} and \ref{lem4.1}  (\ref{eq_thm2}) holds 
for $0\leq x_n <\!< \La_n(x)$ along with  $y$ subject to  $ y_n\in [-M, -1/M]$ or $-\La_n(y)<\!< y_n <M$ (the darkened region in {\sc Figure 2} left). We may add to this  the reflection of the region (\ref{pT2}) about $D_+$, and finally find the symmetrization of  the resulting region   about $D_-$  agrees with  the one described in Proposition \ref{prop1}.  \qed

%\section5{ Case $\Ga=2-\a$: Proof of Theorem \ref{thm3}}
\section{ Proof of Theorem \ref{thm3}}

The proof is based on the  result of Doney \cite[Proposition 11]{D}, given below, for the walk killed not at 0 but on the negative half line. 
 Let  $V_{{\rm a}}$ denote the renewal function of weakly ascending  ladder height process of the walk $S$.  Recall that   $q(\eta)$ and $\hat q(\eta)$, $\eta \geq 0$ denote  the densities of the stable meander of length $1$ at time $1$ for $Y$ and $-Y$, respectively (see (\ref{Bt}) for the definition). 
Doney \cite{D} obtains  an elegant asymptotic formulae of the probability $P[S_n =x-y, \sigma^0_{[x+1,\infty)}>n] $ ($x\geq 0$),
  which under the present assumption and with our notation may be rewritten (by using the duality relation) as 
\beqn\label{Doney}
 Q^n_{(-\infty,-1]}(x,y)   \sim \left\{\begin{array}{lr}
{\displaystyle \frac{U_{{\rm d}}(x)V_{{\rm a}}(y) \mathfrak{p}_{1}(0)}{ nc_n} } \quad & (x_n \downarrow 0, y_n \downarrow 0), \\[3mm]
{\displaystyle  
\frac{V_{{\rm a}}(y)  P[\sigma^0_{[0, +\infty)}>n]  \hat q(x_n)}{ c_n} } \quad  & (y_n \downarrow 0, x_n >1/M ),\\[3mm]
{\displaystyle  
\frac{ U_{{\rm d}}(x) P[\sigma^0_{(-\infty, -1]}>n] q(y_n)}{ c_n} } \quad  & (x_n \downarrow 0, y_n>1/M),\\[3mm]
\mathfrak{p}^{(-\infty,0]}_{1} (x_n,y_n)/c_n \quad & \;\;   (x_n \wedge y_n\geq 1/M) \end{array} \right.
\eeqn 
valid for all admissible $\ga$ uniformly for $x\vee y <Mc_n$. 
 
 Let $\ga=2-\a$ throughout  the rest of this section. 
 It holds  \cite{VW} that  as $n\to\infty$
 \beqn\label{D0}
  P[\sigma^0_{(-\infty, -1]}>n]U_{{\rm d}}(c_n)\to\k_\circ>0
  \eeqn 
  (with  $\k_\circ = 1/\Ga(1-\frac1\a)$  by Lemma \ref{L8.8}(i)) and $U_{{\rm d}}(\xi c_n)/U_{{\rm d}}(c_n) \to \xi$.  Combining these with  the last two equivalences in (\ref{Doney}) we see that $K_1(\eta)=\lim_{\xi\downarrow 0} \mathfrak{p}^{(-\infty,0]}_{1} (\xi,\eta)/\xi= \k_\circ q(\eta)$, so that  the third case in (\ref{Doney}) is equivalently given as
 \beqn\label{D1}
  Q^n_{(-\infty,-1]}(x,y)   \sim  \frac{ U_{{\rm d}}(x) K_1(y_n)}{ U_{{\rm d}}(c_n)c_n}  \quad   (x_n \downarrow 0, y_n\in [1/M, M]).
  \eeqn
Putting $m= \lfloor \de n\rfloor$ with a small number $\de>0$   we decompose
 $Q^n_{\{0\}}(x,y) = I+II$, where
$$ I= \sum_{z=1}^\infty Q^m_{(-\infty,0]}(x,z)   Q^{n-m}_{\{0\}}(z,y) 
\quad\mbox{and}\quad
 II= \sum_{k=0}^{m}\sum_{z=1}^\infty Q^k_{(-\infty,0]}(x,-z)   Q^{n-k}_{\{0\}}(-z,y).
 $$
 [Note that here appears not  $Q^k_{(-\infty,-1]}(x,\pm z)$ but $Q^k_{(-\infty,0]}(x,\pm z)$.]
Asymptotic forms of  the second factors of the  summands are essentially given in Proposition \ref{prop1} (as is observed shortly), by which,  together with (\ref{Doney}), we can compute the double sum with an appropriate accuracy.
Bringing in the function
$$W_{n,x} =  \frac{ U_{{\rm d}}(x)}{U_{{\rm d}}(c_n)c_n}.$$
 our main task then  is performed by  showing that
 \beqn\label{I/II}
 \begin{array}{ll}
(a)\, &\;  I = Q^n_{(-\infty,0]}(x,y)  + o(W_{n,x})\\[2mm]
(b)  &  II = {\displaystyle a(x)\frac{\mathfrak{f}^{-y_n}(1)}{n}\{1+o(1)\}  + W_{n,x}\times o_\de(1) }
 \end{array}
\quad \mbox{  uniformly   for}\;\;\left\{ \begin{array}{ll}   0 < x_n<\!< 1,\\[1mm]
 y_n\in [1/M, M]
 \end{array}\right.
   %\qquad\qquad\qquad\qquad\qquad  \qquad\qquad\qquad\qquad\qquad
  \eeqn
  \vskip1mm\n
  ($o_\de(1)\to0$ as $n\to\infty$ and $\de\downarrow 0$)
 by which the formula of Theorem \ref{thm3} follows if restricted to  this same regime  since $Q_{(-\infty,0]}(x,y)  \asymp W_{n,x}$ ($x_n\downarrow 0, y_n \asymp 1$).    

\v2
{\it Proof of  (\ref{I/II}a). }\, Let $y_n\in [1/M,M]$. Since for any $\e>0$, 
$Q^n_{\{0\}}(z,y) \sim Q^n_{(-\infty,0]}(z,y) $ for $z_n>\e$ in view of Proposition \ref{prop1} and since $Q^n_{(-\infty,0]}(z,y) \leq Q^n_{\{0\}}(z,y) \leq C/c_n$,  it suffices to show that
\beqn\label{5.3a} \frac1{c_n}\sum_{1\leq z <\e c_n} Q^m_{(-\infty,0]}(x,z)   \leq C \e^\a W_{n,x}
\eeqn
in which $\e$ can be made arbitrarily small. 
By (\ref{Doney}) the LHS is dominated by 
$$\frac{U_{{\rm d}}(x)}{nc_n^2}\sum_{z< \e  c_n} V_{{\rm a}}(z) \sim \a \e^\a \frac{U_{{\rm d}}(x) V_{{\rm a}}(c_n)}{nc_n},$$ 
verifying (\ref{5.3a}) since $V_{{\rm a}}(c_n) \sim C_1 n/ U_{{\rm d}}(c_n)$ (cf. Lemma \ref{L8.8}(ii)). \qed
%\begin{lem}\label{T3L0} $ I= Q^n_{(-\infty,0]}(x,y)  + o(W_{n,x})$.
%\end{lem}

\v2
The proof of (\ref{I/II}b) is given  after showing a few  lemmas. Since $\ga=2-\a$ entails that  the ascending ladder height has infinite expectation,  according to  \cite[Corollary 1]{Uladd} 
we have
\beqn\label{a_Rep} 
a^\dagger(x) = E[a(S^x_{\sigma_{(-\infty,0]}})]= \sum_{k=0}^{\infty}\sum_{z=1}^\infty Q^k_{(-\infty,0]}(x,-z) a(-z)\quad x\in \Z
\eeqn
(see (\ref{H/a})).  Put
$R(w) = \sum_{z=1}^\infty p(-z-w)a(-z)$. Using   $G_{B}$   defined in  (\ref{green})
we rephrase  identity (\ref{a_Rep})  as 
\beqn\label{a/R}
a^\dagger(x)= \sum_{w=1}^\infty G_{(-\infty,0]}(x,w)R(w).
\eeqn 
%%%%
%Lem5.1 
\begin{lem}\label{T3L1} There exists a constant  $C$ such that for $y\geq 2x\geq 1$,
\beqn\label{eqL5.1}
\sum_{w=y}^\infty G_{(-\infty,0]}(x,w)R(w) \leq C \frac{U_{{\rm d}}(x)a(y/2)}{U_{{\rm d}}(y)}.
\eeqn
 \end{lem}
 \n
 \pf\,  In this proof  we use Lemma \ref{lem3.1}(ii) which verifies  that  $a(-x)$ is increasing for all sufficiently large $x$---the use of only  (i) of Lemma \ref{lem3.1} suffices but necessitates  certain modifications  involving  a circuitous argument. 
 Summing by parts   we have 
 $$R(w) =  \sum_{j=w}^\infty r(j), \quad \mbox{where} 
 \quad r(j) =\sum_{z=1}^\infty p(-z-j)[a(-z) -a(-z+1)].$$
%By what is mentioned above $r(j)\geq 0$ for $j$ large enough. 
Put $v(x) = V_{{\rm a}}(x)- V_{{\rm a}}(x-1)$, $x\geq 1$ and $v(0)=V_{{\rm a}}(0)$. Then
\beqn\label{r_G}
G_{(-\infty,0]}(x,w) = \sum_{k=1}^{x\wedge w} u(x-k)v(w-k)
\eeqn
(cf.   \cite[Proposition 19.3]{S}).
On summing  by parts again
\beqn\label{pL5.1}
\sum_{w=x+1}^\infty v(w-k)  R(w) = -V_{{\rm a}}(x-k) R(x+1)+  \sum_{w=x+1}^\infty V_{{\rm a}}(w-k) r(w).
\eeqn
 In particular   for $y\geq x+1$,
\beqn
\label{eLT3L1}
\sum_{w=y}^\infty G_{(-\infty,0]}(x,w) R(w)  \leq  U_{{\rm d}}(x)\sum_{w=y}^\infty V_{{\rm a}}(w)r(w).
\eeqn
Make the decomposition  $a(x) =I+II$, where
$$I= \sum_{w=1}^x \sum_{k=1}^w u(x-k)  v(w-k)R(w), \quad II= \sum_{k=1}^x u(x-k) \sum_{w=x+1}^\infty v(w-k)R(w).$$
For all values of  $x$ large enough we have  $R(x)\leq R(w)$ for $1\leq w\leq x$ (because of the monotonicity  of  $a(-z)$ ($z>\!>1$) mentioned at the beginning of the proof) and, accordingly, 
$$I=  \sum_{k=1}^x u(x-k) \sum_{w=k}^xv(w-k)R(w) \geq  \sum_{k=1}^x u(x-k) V_{{\rm a}} (x-k) R(x),$$
 which  together with (\ref{pL5.1}) as well as the regularity of $V_{{\rm a}}$ shows   that 
$$I+II \geq  \sum_{k=1}^x u(x-k) \sum_{w=x+1}^\infty V_{{\rm a}}(w-k) r(w) \geq c\, U_{{\rm d}}(x) \sum_{w=2x}^\infty V_{{\rm a}}(w)r(w)$$ 
with some positive constant $c$.
Rewriting this result as
 $$\sum_{w=2x}^\infty V_{{\rm a}}(w)r(w) \leq C' a(x)/U_{{\rm d}}(x)$$
valid for all  $x\geq 1$ and combining this with (\ref{eLT3L1}) we conclude  that for $y\geq 2x$,
 \[
 \qquad  \qquad  \sum_{w=y}^\infty G_{(-\infty,0]}(x,w)R(w)   \leq C' U_{{\rm d}}(x)a(y/2)/U_{{\rm d}}(y). \qquad  \qquad  \qquad   \qed
  \]

 %Lem5.2
\begin{lem}\label{L5.2} \,Uniformly for $1/M<y_n<M$ and $0\leq x_n <M$,   as $n\to\infty$ and $\e\downarrow 0$ 
$$\frac{\sup \{ Q^n_{(-\infty,0]}(x,z) :\,1\leq z_n<\e  \;\mbox{or} \; z_n>1/\e\}}{Q^n_{(-\infty,0]}(x,y)} \;\longrightarrow\; 0.$$
\end{lem}
\n
\pf\,  For $0\leq x_n <M$,  the laws of  $S^x_n/c_n$ conditioned on   $\sigma_{(-\infty,0]}^x>n$ 
constitute a tight family, 
and the  proof of Lemma \ref{lem4.3} is available with obvious modifications.\qed 
  %Lem5.3 
\begin{lem}\label{T3L2}\,
 Uniformly for $1\leq x < c_n$,  as $y\wedge n\to\infty$  under $y <c_n$
 $$\sum_{k=m}^\infty \sum_{w=1}^y Q^k_{(-\infty,0]}(x,w)R(w) = \frac{U_{{\rm d}}(x)V_{{\rm a}}(y)}{c_m} \times o(1) \qquad (m=\lfloor \de n\rfloor).
$$
\end{lem}
\n
\pf\,
We split the sum defining $R(w)$ at $z=w$.  Recall  $F(t) =P[X\leq t], t\in \R$.  Performing summation by parts we then deduce 
$$\sum_{z=1}^w p(-z-w)a(-z) \leq  C\int_0^w F(-t-w)\frac{t^{\a-2}}{L(t)}dt =
 \frac{L(w)}{w^\a}\int_0^w \frac{t^{\a-2}}{ L(t)}dt\times o(1) = o(1/w)$$
 as $w\to\infty$, where  we have $o(1)$ since $F(-z)= o(L(z)/z^{\a})$.  The other sum is evaluated in a similar way: as a result we obtain
 $$R(w) =o(1/w).$$
Now by (\ref{Doney})  
$Q^k_{(-\infty,0]}(x,w) \leq C  U_{{\rm d}}(x)V_{{\rm a}}(w)/kc_k$ for all $1\leq w \leq c_n$ and $ k \geq m$ and the required estimate follows immediately. \qed

  %Lem5.4
\begin{lem}\label{T3L3}\, Uniformly for $1\leq x\leq y$, as $y\to\infty$
$$\sum_{z= y}^\infty\sum_{w=1}^y G_{(-\infty,0]}(x,w)p(-z-w)a(-z) = \frac{U_{{\rm d}}(x)V_{{\rm a}}(y)}{y} \times o(1).
$$
\end{lem}
\n
\pf\, Since $F(-z)/F(z) \to 0$ as $z\to\infty$, uniformly for $1\leq w \leq y$,
 $$\sum_{z=y}^\infty p(-z-w)a(-z)  \leq F(-y-w)a(-y) + C\sum_{z=y}^\infty F(-z-w)\frac{z^{\a-2}}{L(z)} =o(1/y)$$
 as $y\to\infty$,  whereas
  $\sum_{w=1}^y G_{(-\infty,0]}(x,w) \leq C  U_{{\rm d}}(x)V_{{\rm a}}(y)$, showing  the asserted equality.
 \qed
\v2

{\it Proof of  (\ref{I/II}b).} \; Since   by Proposition \ref{prop1}
$$Q^{n-k}_{\{0\}}(-z,y) = n^{-1}a(-z) \mathfrak{f}^{-y_n}(1)\{1+o_\de (1)\}$$ 
uniformly for $1\leq k\leq m$, $1\leq z < c_m \, (\sim \de^{1/\a}c_n)$,
by virtue of (\ref{a_Rep}) and  (\ref{a/R}) it suffices to show 
$\sum_{k=0}^{m}\sum_{z=1}^\infty Q^k_{(-\infty,0]}(x,-z)  a(-z) = a^\dagger(x)+ o(nW_{n,x}),$
which follows if we prove
that for each $\de>0$  fixed, $$\sum_{k=m}^\infty \sum_{w=1}^\infty Q^k_{(-\infty,0\}}(x,w)R(w)= o(nW_{n,x}); \mbox and$$
$$\sum_{z> c_m}\sum_{w=1}^\infty G_{(-\infty,0]}(x,w)p(-z-w)a(-z) = o(nW_{n,x}).$$
The former one follows from Lemmas \ref{T3L1} and \ref{T3L2}, while the latter from Lemmas \ref{T3L1} and \ref{T3L3}, for
$$a(c_m/2) = o(a(-c_m)) = o(n/c_n ), %\quad V_{{\rm a}}(c_m) \asymp \frac{n}{U_{{\rm d}}(c_n)},   
\quad \frac{V_{{\rm a}}(c_m)}{c_m} \sim C_\de\frac{n}{U_{{\rm d}}(c_n)c_n}. \qquad \qed$$

\v2
{\it Proof of Theorem \ref{thm3}.} \;We have already shown (\ref{I/II}a) 
and,  accordingly,  the second  formula of (\ref{eq_thm3}) restricted to $0<x_n <\!< 1$.  The cases  $x_n>1/M$ and  $-\La_n(x)<\!< x_n \leq 0$ being included in Proposition \ref{prop1},
 up to now  we have shown
\beqn\label{eq_T3}
Q^n_{\{0\}}(x,y) \sim    \frac{a^\dagger(x)}{n}\mathfrak{f}^{\,-y_n}(1) + \frac{U_{{\rm d}}(x)K_1(y_n)}{U_{{\rm d}}(c_n)c_n}   \quad\;\;   (-\La_n(x)<\!< x_n < M,  y_n \in [1/M, M]).
\eeqn

The first formula of  (\ref{eq_thm3}) is derived from (\ref{eq_T3})  as in the proof of Lemma \ref{lem4.4}. Indeed by Proposition \ref{prop1} we have  $Q^n_{\{0\}}(z,y) \sim \mathfrak{f}^{\, z_n}(1)a^{\dagger}(-y)/n$ for $z_n \in [M^{-1}, M]$, $-\La_n(-y)<\!<y_n<\!< 1$, and  on putting $m=\lfloor \la n\rfloor$ and applying (\ref{eq_T3}) to   $Q^m_{\{0\}}(x,z)$
$$\sum_{z=1}^{Mc_n} \frac{U_{{\rm d}}(x)K_1(z_m)}{U_{{\rm d}}(c_m)c_m} Q^{n-m}_{\{0\}}(z,y) 
\sim \frac{U_{{\rm d}}(x)a^\dagger(-y)}{nU_{{\rm d}}(c_n)} \sum_{z=1}^{Mc_n} \frac{K_1(z_n\la^{-1/\a})
\mathfrak{f}^{\, z_n/(1-\la)^{1/\a}}(1)}{\la^{2/\a}(1-\la)c_n}.$$
Using  the scaling relations of $K_t(\eta)$ and  $\mathfrak{f}^{\eta}(t)$ we see that the  above sum approaches 
$$\int_0^\infty K_\la(\eta) \mathfrak{f}^{\, \eta}(1-\la)d\eta = \lim_{\xi \downarrow 0} \frac{1}{\xi} \mathfrak{f}^{\, \xi}(1)= \mathfrak{p}_1(0),$$
and recalling (\ref{f*f}) as well as  the computation leading to (\ref{Q}) we have  the first formula of  (\ref{eq_thm3}). This finishes the proof. \qed

\v2
%Rem5.1
\begin{rem}\label{R_D}\, Let  $\ga=2-\a$ and put  $\ell^*(x) =\int_0^x P[-\hat Z>t]dt$.  Here we observe that  as  $n\wedge y\to\infty$ under $\ell^*(y)/\ell^*(c_n) \to 1$ and $ y_n <\!< 1$
 \beqn\label{DT1}
  \begin{array}{lr}
 P[S^x_n=y,   \sigma^x_{(-\infty,-1]} <n \leq  \sigma^x_{\{0\}}] \sim   f^0(n) a^\dagger(x)a^\dagger(-y)  & (x_n \downarrow 0),\\[2mm]
Q^n_{(-\infty,0]}(x,y) \sim {\displaystyle \frac{U_{{\rm d}}(x-1)\mathfrak{p}_1(-x_n)}{U_{{\rm d}}(c_n)n}a(-y) } & (0 < x_n <M).
 \end{array} 
\eeqn
The   result corresponding to the first equivalence  for the case $y_n\in [1/M,M]$ that may read
 \beqn\label{DT}
 P[S^x_n=y,   \sigma^x_{(-\infty,-1]} < n < \sigma^x_{\{0\}}] \sim  {\displaystyle   \frac{a^\dagger(x)}{n}\mathfrak{f}^{\,-y_n}(1)}
 \quad  (x_n \downarrow 0,  y_n \in [1/M, M])
\eeqn
 follows immediately from the manner we have derived (\ref{eq_T3}). By the same token the first equivalence of (\ref{DT1})  holds for  $0 \leq x_n <\!< \La_n(x)$ and $0\leq y_n <M$ when the RHS of it is the leading term of $Q^n_{\{0\}}(x,y)$.  From  these fact  we infer that in most subcases of $0\leq x <\!< c_n$ the sums in the  formulae on the RHS of  (\ref{eq_thm3}) (in  Theorem \ref{thm3}) correspond to the obvious decomposition of $Q^n_{\{0\}}(x,y)$ as the
 sum of the probabilities on the LHS of  (\ref{DT1});  it in particular follows that the dominant contribution to $Q^n_{\{0\}}(x,y)$ comes  from the walk trajectories  that enter the negative half line before  $n$ if $x_n<\!< \La_n(x)$ and those that do not if   $ \La_n(x) <\!< x_n $. 
 
 For verification of (\ref{DT1}) first note  that
 the function $\ell^*(x) $ varies slowly at infinity and $U_{{\rm d}}(x) \sim x/\ell^*(x)$). (Cf. \cite[Lemma 12]{VW}.) 
 %and $u(x) \sim 1/\ell^*(x), \cite[Lemma 7.1]{Uladd}.) 
 It also follows that  $V_{{\rm a}}(y) \sim y^{\a-1}\ell^*(y)/[L(y)\Ga(\a)]$ (cf. Lemma \ref{L8.8}). The second relation then is checked by looking at the first formula in (\ref{Doney}) if 
 $x_n \to0$. In case $x_n>1/M$ one can resort to the analogue of  formula (\ref{D1}) which reads
 \beqn\label{D1V}
  Q^n_{(-\infty,-1]}(x,y)   \sim  \frac{ V_{{\rm a}}(y) x_n\mathfrak{p}_1(-x_n)}{ \Ga(\a)V_{{\rm a}}(c_n)c_n}  \quad   (y_n \downarrow 0, x_n\in [1/M, M]):
  \eeqn
for the derivation   use  $P[\sigma^0_{[0, \infty)}>n]\hat q(\xi) = \xi\mathfrak{p}_1(\xi)/[\Ga(\a)V_{{\rm a}}(c_n)]$
(cf. Lemmas \ref{L8.6} and \ref{L8.8}). The first formula of (\ref{DT1}) follows from the second of it together with the first case of  (\ref{eq_thm3}). Thus  (\ref{DT1}) has been verified. 
\end{rem}

%(e) \,  comparing (\ref{eq_thm3})   with the asymptotic form  of $Q_{(-\infty,0]}^n(x,y)$ in \cite{D} (see Section 5 of the present paper) we then  infer that    uniformly for $0\leq x,  y < Mn^{1/\a}$, as $n\wedge y \to\infty$ 
%\beqn\label{comp}
%Q_{\{0\}}^n(x,y) \sim Q_{(-\infty,0]}^n(x,y) + a^\dagger(x) f^0(n) a(-y)\eeqn
%(this follows also from the proof of Theorem \ref{thm3}) and  for $ M^{-1}  \La_n(x) < x_n < M \La_n(x)$  the two terms on the RHS are comparable with each other, of which the first  represents the  contribution of the trajectories of the walk  staying in $[1,\infty)$    and the second the contribution  of those  entering into $(-\infty,0]$ (both until $n$).
% (Note that  the first term is superfluous for $x_n <\!< \La_n(x)$ and the second for $x_n >\!> \La_n(x)$.)

%Section6
%section{Esitmation of $p^n_{\{0\}}(x,y)$ in case $xy< 0$}
\section{Upper and lower bounds of $Q^n_{\{0\}}(x,y)$  and proof of Theorem \ref{thm4}}

Here we derive estimates of $Q^n_{\{0\}}(x,y)$ for $x, y$ not necessarily confined in $ |x_n|\vee |y_n|<M$, that lead to Proposition  \ref{prop2.4} and are useful  for the proof of   Theorem \ref{thm4}. Throughout this section
 we  assume $\ga =2-\a$ unless stated otherwise  explicitly.  %the case $xy<0$ for $|\ga|<2-\a$ being included in Theorem \ref{thm2}.    
  Sometimes we suppose $E[-\hat Z] <\infty$, which entails $\ga=2-\a$.
%[In case  $E|\hat Z|=\infty$ and  $\ga=2-\a$ there arises  a troublesome question caused by the obscure nature of $L(x):= U_{{\rm d}}(x)/x$ (cf. \cite{VW}; see also Remark \ref{Rprop1}(a)).]  
 
  By (\ref{La/x}) (and by $x_n \leq U_{{\rm d}}(x)/U_{{\rm d}}(c_n)\{1+o(1)\}$ ($0<x < c_n$)) it follows that
 \beqn\label{U/x}
 \frac{a^\dagger (x)c_n}{n}+\frac{U_{{\rm d}}(x)}{U_{{\rm d}}(c_n)} 
 \asymp \frac{a^\dagger(x)c_n}{n} + x_n \qquad (0<x <Mc_n),
 \eeqn 
by which  the  two expressions on the above  two sides will be  interchangeable in most places of this section.

%Prop6.1
\begin{Prop}\label{prop6.1}   Let $\ga=2-\a$. For each $M> 1$, 
 there exists a positive constant  $c$ such that for  $-M< y_n < 0<  x_n<M$, 
\beqn\label{eqP1}
Q^n_{\{0\}}(x,y) \geq c[ D_n(x,y) \vee D_n(-y,-x)],
\eeqn
where
$$D_n(x,y) =  \bigg( \sum_{z=2}^{x} p(-z)zV_{{\rm a}}(z)a(-z)\bigg)\frac{U_{{\rm d}}(x)}{x} 
\bigg[  \frac{a^\dagger(-y)}{n^{2}/c_n} + \frac{|y_n|}{n} \bigg].$$
\end{Prop}
\n\pf\,  The walk is supposed to be not left-continuous, otherwise the result being trivial.
Let  $j(x)$  denote the smallest integer $j$ such that $c_j \geq  x$.  This proof employs   the obvious lower bound 
   $$Q^n_{\{0\}}(x,y) \geq  \sum_{\de j(x) \leq k \leq n/2}\; \sum_{1\leq w \leq c_{k/\de}}\;\sum_{z=1}^{x}Q^{k-1}_{(-\infty,0]}(x,w)p(-z-w) Q^{n-k}_{\{0\}}(-z, y)
$$
valid for any constant $\de>0$.   $\de$  needs  to be  chosen so small  that for all $n, x$ large enough, $c_{\de j(x)}< \eta c_n$ for some $\eta<1$.  This is fulfilled  with $\eta =1/2$ by taking
$\de = 1/(3 M)^\a$,
for $c_{\de j(x)}/c_n\sim \de^{1/\a}x_n <\de^{1/\a}M =1/3$. Let $0<x_n, -y_n<M$.

 For $k, w, z$ taken from the range of summation above, we have  by Theorem \ref{thm3} (see also Corollary \ref{cor2})
$$Q^{n-k}_{\{0\}}(-z,y) =Q^{n-k}_{\{0\}}(-y,z) \asymp f^{-y}(n-k)a(-z)\asymp  \bigg\{\frac{a^\dagger(-y)}{n^{2}/c_n}  +\frac{U_{{\rm d}}(-y)}{U_{{\rm d}}(c_n)n}\bigg\}a(-z),
$$
and, noting that for  $k\geq \delta j(x)$ and  $x\leq c_{k/\delta} \sim c_k\de ^{-1/\a}$, we apply (\ref{Doney}) to obtain 
$$Q^{k-1}_{(-\infty,0]}(x,w) \asymp  U_{{\rm d}}(x)V_{{\rm a}}(w) /kc_k.$$
  Hence,    putting
\beqn\label{s_m}
m(x) =  \sum_{z=1}^x\;\sum_{w=1}^x p(-z-w)V_{{\rm a}}(w)a(-z)
\eeqn 
we have
$$Q^n_{\{0\}}(x, y) \geq c' m(x)U_{{\rm d}}(x)\bigg\{\frac{a^\dagger(-y)}{n^{2}/c_n}  +
\frac{|y_n|}{n}\bigg\}\sum_{\de j(x)\leq k \leq n/2}\;\frac1{ kc_k},$$
where (\ref{U/x}) is used. 
Since  $x\sim c_{j(x)}$ and,  by our choice of $\de$,  $c_{\de j(x)} <c_n/2\sim 2^{1/\a -1}c_{n/2}$,  
the last sum is bounded below by  a positive multiple of $1/x$. 
 In the double sum in (\ref{s_m}) restricting its range of  summation to $z+w\leq x$ and  making change of variables we see
\beqn\label{m_lb}
 m(x)\geq  
 \sum_{k=2}^xp(-k) \mathop{{\sum}^{*}}_{z=-k+1}^{k-1\;} V_{{\rm a}}\Big(\frac{k-z}{2}\Big) a\Big(\!\!-\frac{k+z}{2}\Big),
 \eeqn
where under the symbol $\mathop{{\sum}^*}$ the summation is restricted to $z$ such that $k-z$ is even.  It is easy to see that $m(x)\geq c'' \sum_{k=2}^x p(-k)kV_{{\rm a}}(k)a(-k)$, showing
 the required
lower bound of (i) in view of duality.  \qed

\v2
%Remark6.1
\begin{rem}\label{r_P6.1}   (i)  If $F(x)$ is regularly varying as $x\to -\infty$ (of index $ -\beta$ necessarily  $\beta\geq \a$), then   $\sum_{z=2}^x p(-z)zV_{{\rm a}}(z)a(-z) \sim C_1 \sum_{z=2}^x F(-z)V_{{\rm a}}(z)a(-z) \sim C_2 a(x)x/U_{{\rm d}}(x)$ with positive  constants $C_1, C_2$  (the first  equivalence  is due to Karamata's theorem and  the second to \cite[Proposition 6.2]{Upot}),   
so that $D_n(x,y) \asymp a(x)[a^\dagger(-y)c_n/n^2 + |y_n|/n]$ and in view of Proposition \ref{prop6.2}(i) (given shortly) the lower bound (\ref{eqP1})  is exact so that  for some positive constant $C$,
\[
\begin{array}{lr}
&{\displaystyle C^{-1} Q^n_{\{0\}}(x,y) \leq  \frac{a^\dagger(x)a^\dagger(-y)}{ n^{2}/c_n} +\frac{a^\dagger(-y)x_n + a^\dagger (x) y_n}{n}  \leq C Q^n_{\{0\}}(x,y)}\\[4mm]
&\quad  \mbox{for} \quad -M <y_n\leq 0\leq x_n<M.
\end{array}
\]

(ii)\, Suppose that $E[-\hat Z]<\infty$ and  $\lim_{x\to\infty} a(x) =\infty$. Since  $U_{{\rm d}}(x)\sim x/E[-\hat Z]$ and the latter condition implies (in fact is equivalent to)   $\sum_{z=2}^{\infty} p(-z)zV_{{\rm a}}(z)a(-z) =\infty$, it then follows from  Proposition \ref{prop6.1}  that
as  $n\wedge x\wedge (-y)\to \infty$ under $x\wedge (-y) <Mc_n$
$$\frac{Q^n_{\{0\}}(x,y)}{(\La_n(x)\vee x_n \vee \La_n(-y) \vee |y_n|)/n}\,  \longrightarrow \, \infty.$$ 
 \end{rem} 

%Let   $E|\hat Z|<\infty$ so that   $a(x) \asymp \sum_{w=1}^{x}\sum_{z=1}^\infty  p(-w-z)[ V_{{\rm a}}(z)]^2$ (\cite[Theorem 2(i), (iii)]{Uladd}) and $U_{{\rm d}}(x)\sim x/E|\hat Z|$. Suppose  that $F(x)$ is regularly varying as $x\to -\infty$ of index $ -\beta$ (necessarily  $\beta\geq \a$).  %If  $\beta = \a$, then  (\ref{eqP21}) reduces to (\ref{eqP1}) (whether   $E|\hat Z|$ is finite or not).   
%Then we deduce  $ \sum_{w=1}^{x} p(-w)w^{2\a-1} \asymp  a(x)$ 
%so that in view of Proposition \ref{prop6.2}(i) (given shortly) the lower bound (\ref{eqP1})  is exact. 

 \v2
% Here we state the following consequence of Theorems \ref{thm5} for convenience of citation.  

%Lem6.1
\begin{lem}\label{L6.1}\,  {\rm (i)}  If $\ga=2-\a$, then  for any $M>1$  there exists a constant $C$ such that 
$$ Q_{\{0\}}^n(x,y) \leq  C\bigg[\frac{a^\dagger(x)}{n/c_n}  
+ (x_n)_+ \bigg]
\times \left\{
 \begin{array}{lr}
{\displaystyle \bigg(\frac{a^\dagger(-y)}{n} \wedge  \frac{L(y)n/c_n}{y^\a\vee 1}\bigg)} 
  &\mbox{for} \;\;  |x_n|\leq M, \, y \geq 0,\\[4mm]
{\displaystyle 
o\bigg( \frac{L(y)n/c_n}{|y|^\a}\bigg) }  &\mbox{for} \;\;  |x_n|\leq M, \, y_n \leq -1/M.
\end{array}\right.
$$

{\rm (ii)}   If  $|\ga| <2-\a$,  
$${\displaystyle Q_{\{0\}}^n(x,y) \leq C \frac{a^\dagger(x)}{n/c_n}  
 \bigg(\frac{a^\dagger(-y)}{n} \wedge  \frac{L(y)n/c_n}{|y|^\a\vee 1}\bigg)} \quad \mbox{for} \quad  |x_n|\leq M.
 $$
\end{lem}
\n
\pf\,  Let $ |x_n|<M$.   By Theorem \ref{thm3} we have as before 
 $$
 Q_{\{0\}}^n(x,y) \sim
  f^x(n)a^\dagger(-y) \asymp  \bigg[\frac{a^\dagger(x) c_n}{n} +
  \frac{U_{{\rm d}}(x) {\bf 1}(x\geq 0) }{U_{{\rm d}}(c_n)}\bigg]\frac{a^\dagger(-y)}{n} \quad \mbox{for }\;\; 0\leq y_n \leq 4M.$$
Because of  (\ref{U/x}), in case  $y\geq 0$   it therefore suffices to show
\beqn\label{y>>n}
  Q_{\{0\}}^n(x,y) \leq C_M\bigg[\frac{a^\dagger(x)c_n}{n}  
+ (x_n)_+\bigg] \frac{n/c_n}{a(-y)y} \quad\mbox{for }\;\; y_n>4M.
\eeqn
  Putting $R=\lceil y/3 \rceil $, $N= \lfloor n/2 \rfloor$ we make the decomposition
\begin{eqnarray}\label{eq_L6.1}
Q_{\{0\}}^n(x,y) &=&
\sum_{k=1}^{N} \sum_{z\geq R}  P[S_{k}^x=z,  \sigma^x_{[R,\infty)} =k >\sigma^x_{\{0\}}]Q^{n-k}_{\{0\}}(z,y) \nonumber\\
&&+ \sum_{z < R}P[  \sigma^x_{[R,\infty)}\wedge \sigma^x_{\{0\}} > N, S^x_{N} =z]Q_{\{0\}}^{n- N}(z,y)  \nonumber\\
&=& J_1+ J_2\quad \mbox{(say)}. 
\end{eqnarray}  %cf. Lemma \ref{lem6.6})
By a  local  large deviation bound  (cf. \cite[Theorem 1]{D2}, \cite[Theorem 2.3]{Brg}) it follows that
\beqn\label{(E)}
P[S_n = z] \leq C_0 nL(z)|z|^{-\a}/c_n  \quad (z\neq0)
\eeqn 
which   entails 
\beqn\label{(B)}
\sup_{z\leq 2R}\,\sup_{ k\leq N} Q^{n-k}_{\{0\}}(z,y) \leq C'_0\frac{n L(R)}{c_nR^\a}.
\eeqn
Rewrite $J_1$ as
$$J_1 = \sum_{k=1}^{N}  \sum_{w<R, w\neq 0} \,\sum_{z = R}^\infty Q^{k-1}_{\{0\}\cup [R,\infty)}(x,w)p(z-w) Q^{n-k}_{\{0\}}(z,y). % \quad\mbox{where} \quad D_k(w) = \sum_{z=R}^\infty .
$$ 
 We split the inner most summation  at $z=2R$.      By  (\ref{(B)})  it follows that 
\beqn\label{(C)}
\sum_{k=1}^{N} \sum_{0\neq w<R}\sum_{z= R}^{2R}  Q^{k-1}_{\{0\}\cup[R,\infty)}(x,w)p(z-w) Q^{n-k}_{\{0\}}(z,y) \leq C_1\frac{n L(R)}{c_nR^\a} P[\sigma^x_{[R,\infty)}< \sigma^x_{\{0\}}]. 
\eeqn
 The contribution from $z>2R$ is at most a constant multiple of
  \beqn\label{(D)}
 \sum_{z>2R}p(z-R)/c_n \leq C  L(R)/c_n R^\a \asymp   L(y)/c_n y^\a.
  \eeqn 
 In Appendix (C) (Lemma \ref{lem9.4}) we prove
 \beqn\label{(A)}
P[ \sigma^x_{[R,\infty)} <\sigma^x_{\{0\}}] \leq C'[a^\dagger(x)/a(-R)  + x_+ R^{-1}].
\eeqn
Collecting  these bounds we can conclude 
$$J_1 \leq C' \bigg[\frac{a^\dagger(x)}{a(-y)}  + \frac{x_+}{ y} + \frac1{n} \bigg]
\frac{n/c_n}{a(-y)y},$$
On the other  hand   on employing  the bound (\ref{(B)}) again
%$p^n(x)\leq C |x|^{-\a}L(x)n/c_n$ (\cite[Theorem 1]{D2}, \cite[Theorem2.3]{Brg})%Cn^{1-1/\alpha}/|x|^\alpha$ (see Lemma \ref{lem6.6})
\beq
J_2 \leq P[\sigma^x_{\{0\}} \geq \tst 12 n]\sup_{ z\leq R} Q_{\{0\}}^{n-N}(z,y) 
&\leq& C [nf^x(n)] nL(y)/y^\a c_n \\
&\leq& 
 C' \bigg[\frac{a^\dagger(x)c_n}{n} +  (x_n)_+\bigg] \frac{n/c_n}{a(-y)y},
 \eeq
 where  (\ref{U/x}) is used for the last inequality.
Noting that   $n/c_n a(-y)$ is bounded for $y>c_n$    and comparing terms we find (\ref{y>>n})  obtained.  Thus the assertion  has been proved in case $y\geq0$. If  $\ga<|2-\a|$, the same proof also shows the result  for $y<0$. 

Let  $\ga=2-\a$ and $y<- c_n/M$.   Because of $F(y)= o(L(y)/|y|^\a)$ we have, instead of  (\ref{(E)}), 
 $p^n(y) = o(nL(y)/|y|^\a c_n)$ ($y \to -\infty$)
according to \cite{Brg}.  Hence, on putting $R= \lfloor -y/3\rfloor$,
$$
\sup_{z\geq -2R}\,\sup_{ k\leq N} Q^{n-k}_{\{0\}}(z,y) \leq o(n L(R)/c_nR^\a).
$$
On the other hand by Lemma 5.5 of \cite{Upot}
$$P[ \sigma^x_{(-\infty, -R]} <\sigma^x_{\{0\}}] \leq [a^\dagger(x)/a(-R)]\{1+o(1)\}.$$
With these bounds we can follow the same lines as  above with $\sigma^x_{(-\infty,-R]}$
in place of $\sigma^x_{[R,\infty)}$ to obtain the asserted estimate of the lemma.
\qed

\v2 For the convenience of later citations  we write down two simplified and (partly reduced) version of the bounds of Lemma \ref{L6.1}.
In the next one  $\gamma$ may be any admissible constant.
%Lem6.2 
\begin{lem}\label{L6.2}\, There exists a constant $C$ such that for all  $n\geq 1$,
\[
Q^n_{\{0\}}(x,y)  \leq C\frac{ [|y|^{\a-1}/L(y)]\wedge (n/c_n)}{|x|^{\a}/L(x)}
\quad\mbox{if}\quad 1\leq   |y|<\frac12 |x|. \]
%where $o(1)$ is bounded and tends to zero as $|x_n|\to \infty$.
\end{lem}
\n
\pf\, For $|y_n| \geq M$ the asserted bound follows from the large deviation estimate (\ref{(E)}). For  $|y_n|< M$,  it is rephrased in  the dual
form which  is given as
\beqn\label{eqL5.4}
Q^n_{\{0\}}(x,y) \leq C  |x|^{\a-1}L(y)/ (|y|^{\a}L(x)) \quad  \mbox{for} \quad 0<|x_n|<M,  |y_n| >1/M,
\eeqn
and hence follows from Lemma \ref{L6.1}. \qed

%Lem6.3
\begin{lem}\label{L6.3} \, Let $\ga=2-\a$. Then
\beqn\label{eqL6.3}
Q^{n}_{\{0\}}(x,y)  \leq C[a^\dagger(-y)\wedge a(-c_n)]\frac{a^\dagger(x)c_n/n +(x_n \wedge 1)}{n}  \quad  (x\geq 0, y\geq 0).
\eeqn
\end{lem}
\n
\pf\, This follows from  Lemma \ref{L6.1}(i) for $x_n<1$  and  Lemma \ref{L6.2}  for $x_n\geq 1$. \qed

%Lem6.4
\v2
\begin{lem}\label   {L6.4}  Suppose  $\ga = \a-2$ and define $\omega_{n,x,y}$ for $x\neq 0$ and $y>0$ via
\beqn \label{iv}
 Q_{\{0\}}^n(x,y)= a(-y)f^x(n)\om_{n,x,y}.
 \eeqn
Then,  $\omega_{n,x,y}$ is dominated by a constant multiple of  $1 \wedge  y^{-1}$ (in particular uniformly bounded), and tends to unity as $y_n \to 0$   and $n\to\infty$ uniformly for  $0<x_n < M$  for each  $M>1$. 
\end{lem}
\n
\pf\, The convergence of $\omega_{n,x,y}$  to unity follows from Theorems \ref{thm2} and
 \ref{thm3}  and  the asserted bound of $\om_{n,x,y}$ follows from Lemma \ref{L6.1}.
 \qed

%Prop6.2
\begin{Prop}\label{prop6.2}\, Suppose $\ga=2-\a$. Then for some constant $C$
\v2
{\rm (i)}\quad  ${\displaystyle  Q^n_{\{0\}}(x,y) \leq  
C\bigg[\frac{a^\dagger(x)a^\dagger(-y)}{ n^{2}/c_n} +\frac{  a^\dagger(-y)(x_n\wedge1)+ a^\dagger (x)(|y_n|\wedge 1)} { n} \bigg]  } \quad  (x\geq 0, y\leq 0) $,
\v2
{\rm (ii)}\quad 
 ${\displaystyle  Q^n_{\{0\}}(x,y) \leq  \frac{C}{c_n}
\bigg[\frac{a^\dagger (-y)}{n/c_n} + (y_n)_-\bigg]\bigg (\frac{a(x)}{n/c_n}\wedge \frac{n/c_n}{a(x)|x_n|}\bigg)  }\qquad  (x \leq -1, |y|\leq Mc_n). $
 %\begin{flushright}
%\end{flushright}
%$ \qquad\qquad\qquad\qquad\qquad\qquad\qquad\qquad\qquad\qquad\qquad\qquad\qquad\quad $,\\
%where $H^x_B\{a\} =\sum_{z\in B} H^x_B(z)a(z)$. 
%[If $E|\hat Z|=\infty$, then $H^z_{[0,\infty)}\{a\} =a(z)$ for $z<0$ so that the RHS of (ii) has the same form as that of (i), whereas  if $E|\hat Z| <\infty$,   $H^z_{[0,\infty)}\{a\} =o(a(z))$ as $z\to -\infty$.]
\end{Prop}
\n
\pf\,  
For the proof of (i)  we apply Lemma \ref{L6.3}  to have
\beqn\label{P2_1}
Q^{n-k}_{\{0\}}(z,y) =Q^{n-k}_{\{0\}}(-y,-z) \leq Ca(z)\frac{a^\dagger(-y)+(|y_n|\wedge 1)n/c_n}{n^2/c_n}  \quad  (z<0, k\leq n/2).
\eeqn
 We have  $ES_{[0,\infty)}=\infty$ so that  for $x\geq 0$,  $Ea(S^x_{\sigma_{(-\infty,0]}})= a^\dagger(x)$ (by (\ref{H/a})), and  from (\ref{P2_1}) we deduce 
\beq
P[\sigma^x_{(-\infty,0]} \leq n/2, \sigma^x_{\{0\}} >n, S^x_n= y] 
&\leq&
\sum_{z<0}P[S^x_{\sigma_{(-\infty,0]}}=z]\sup_{k\leq  n/2}Q^{n-k}_{\{0\}}(z,y) \\
&\leq& Ca^\dagger(x)\frac{a^\dagger(-y)+(|y_n|\wedge 1)n/c_n}{n^2/c_n }. 
\eeq

Let $\hat S^x$ and $\hat \sigma_B^x$ denote the dual walk and its hitting time, respectively. Then 
\[
Q^n_{\{0\}}(x,y) -P[\sigma^x_{[0,\infty)} \leq n/2, \sigma^x_{\{0\}} >n, S^x_n= y]  
 \; \leq \;
 P[\hat \sigma^y_{(-\infty,0]} \leq n/2, \hat \sigma^y_{\{0\}} >n, \hat S^y_n= x].
\]
By duality relation the probability on the RHS is the same as what we have just estimated but with $x$ and $y$ replaced by $-y$ and $-x$, respectively, and hence we can  conclude  (i).

By duality (ii) is immediate from Lemma \ref{L6.1}(i). \qed
\v2

%The decomposition  (\ref{mix})  is almost obvious from our proof of  Theorem \ref{thm3}  (see 
%Remark \ref{R_D}) (except for the uniformity in $y$ outside $[0, c_n]$) 

\v2
%Prf of Cor 4
{\bf Proof of Corollary \ref{cor-1}.}
The decomposition  (\ref{mix}) follows from what is mentioned in 
Remark \ref{R_D}. Here we observe that  it  is actually  involved in  Theorem \ref{thm3} itself.
Indeed, we have $K_1(\xi) =\a \mathfrak{p}_1(0)q(\xi)$ (Lemma \ref{L8.6})  and $h(\xi) =\mathfrak{f}^{-\xi}(1)/\kappa$, which together with the second case of   (\ref{eq_thm3}) shows
that 
for $0\leq x <\!< c_n$,  the conditional probability on the LHS of (\ref{mix}) 
is written as
\beq
&&\left\{\begin{array}{lr}
[\th'_{n,x}h(y_n)+\th''_{n,x}q(y_n)]/c_n  +o(1/c_n) &\quad (y_n\in [1/M,M])\\[2mm]
[(1-\frac{1}{\a})\th'_{n,x} + \a^{-1} \th''_{n,x}]a^\dagger(-y)/n  +o(1/c_n)&\quad (0\leq y_n <\!<1)
\end{array}\right. \\[3mm] 
&&\;\; = \th'_{n,x}\frac{h(y_n)}{c_n}+\th''_{n,x}\frac{q(y_n)}{c_n} +\bigg(1+\frac{L(c_n)}{L(y)} \bigg)\times o\bigg(\frac1{c_n}\bigg),
\eeq
 where
$$\th'_{n,x} = \frac{(1-\frac1{\a})^{-1}\k_{\a,\ga} \, \La_n(x)}{P[\sigma^x_{\{0\}}>n]},   \quad \th''_{n,x} = \frac{\a \mathfrak{p}_1(0) U_{{\rm d}}(x-1) /U_{{\rm d}}(c_n)}{P[\sigma^x_{\{0\}}>n]}.  $$
On the other hand according to \cite{D}
%$P[\sigma^x_{(-\infty,0]} > n] \sim $, which together  shows 
 $$(1-\th_{n,x})P[\sigma^x_{\{0\}}>n]  =P[\sigma^x_{(-\infty,0]} > n] \sim 
U_{{\rm d}}(x-1)/U_{{\rm d}}(c_n)\Ga(1-1/{\a}) \quad (x_n\downarrow 0),
 $$ 
which, with the help of  $\a\mathfrak{p}_1(0)\Ga(1-\frac1{\a})=1$,  yields $1-\th_{n,x} \sim \th''_{n,x}$. By (\ref{eq_cor2}) we deduce
$$P[\sigma^x_{\{0\}}>n] \sim \kappa \La_n(x) +\a\mathfrak{p}_1(0)U_{{\rm d}}(x)/U_{{\rm d}}(c_n)$$
so that $\th'_{n,x}+\th''_{n,x}
\to 1$. Hence $\th'_{n,x} = \th_{n,x} +o(1)$, showing that  (\ref{mix}) holds uniformly for $0\leq x <\!< c_n$ and $y\in [0,Mc_n]$ for each  $M>1$.

 From Proposition \ref{prop6.2}(i) and Lemma \ref{L6.1}(i)   it follows that for $0\leq x <Mc_n$, 
$$Q^n_{\{0\}}(x,y)/ [\La_n(x) + x_n] \leq C[\La_n(-y)   + |y_n|\, ]/c_n  = o(1/c_n) \quad (-1 <\!<y_n< 0) 
$$
and 
$$Q^n_{\{0\}}(x,y)/ [\La_n(x) + x_n] = o(1/c_n) \quad (y_n  < -1/M \;\; \mbox{or} \;\; y_n\to\infty) , $$
respectively. These together show that   (\ref{mix})  holds also uniformly for $y\notin [0, Mc_n]$, since  $P[\sigma^x_{\{0\}}>n] \asymp \La_n(x)+x_n$ ($0\leq x <Mc_n$). 
\qed
\v2
In the sequel it is convenient to  bring in some notation. For $B\subset \Z$, $H^x_B(y)$ denote the  mass function of hitting distribution to $B$: 
  $H_B^x(y) = P[S^x_{\sigma_B}=y] $. It follows that
\beqn\label{H_B}
H_B^x(y) =  \sum_{z\notin B}G_B(x,z) p(y-z) \quad \mbox{for}\quad y\in B.
\eeqn
We write $H^x_B\{\fa\} =\sum_{y\in B} H_B^x(y) \fa(y)$ for a function $\fa$ on $B$.  
For the special case  $B= (-\infty,0]$ denote by  $h^x(n,y)$  the corresponding space-time mass function: 
\beqn\label{D_h}
h^x(n,y) = P[\sigma^x_{(-\infty,0]} =n, S^x_n=y]\qquad (y\leq 0),
\eeqn
which is the restriction of $Q_{(-\infty,0]}(x,y)$ to $y\leq 0$.  Write also $H_{(-\infty,0]}^{+\infty}(y)= \lim_{x\to\infty}H_{(-\infty,0]}^{x}(y)$;  $H_{(-\infty,0]}^{+\infty}$ is a probability on $(-\infty,0]$ if and only if $E[-\hat Z]<\infty$ (\cite{S}).

%Lem6.5
\begin{lem}\label{L6.5} Suppose $E[-\hat Z]<\infty$. Then,

{\rm (i)} for $M>1$ and $\e>0$,  uniformly for $0 < x_n<M$ and $y\leq 0$
$$h^x(n, y) = \frac{x_n \mathfrak{p}_{1}(-x_n)}{n} \Big[ H^{+\infty}_{(-\infty,0]}(y)\{1+o_\e(1)\} + r(n,y)\Big]$$
where $o_\e(1)$ is bounded and tend to zero  as $n\to\infty$ and $\e\to 0$ in this order and
$$|r(n,y)| \leq  C_M\sum_{z>\e c_n} a(-z)p(y-z)$$
for a constant $C_M$   depending  only on  $M$ and $F$; and

 {\rm (ii)} there exists a constant $C$ such that for all $x\geq 1, y <0$ and  $n\geq 1$,
$$h^x(n, y) \leq C  \bigg(\frac1{n}\wedge \frac{L(x)}{x^\a}\bigg)x_n  H^{+\infty}_{(-\infty,0]}(y).$$
 \end{lem}
\n
\pf\,   Suppose $E|\hat Z|<\infty$, so that  $U_{{\rm d}}(x)\sim x/E|\hat Z|$.   Let $\e>0$ and  in the expression  
\beqn\label{++}
h^x(n+1,y) =\sum_{z=1}^\infty Q^{n}_{(-\infty,0]}(x,z)p(y-z)
\eeqn
 we divide the sum into two parts,  the sum on  $z< \e c_n$  and 
the remainder which  are denoted by $\Sigma_{<\e c_n}$ and    $\Sigma_{\geq\e c_n}$, respectively.  %We know that $U_{{\rm ds}}(x)V_{{\rm as}}(z) \sim xz^{\alpha-1}/c_\circ \Gamma(\alpha)$ ($x\wedga z \to\infty$) 
By Doney's result (\ref{Doney}) (cf. also  (\ref{D1V}))  together with   $\Ga(\a) V_{{\rm a}}(c_n) \sim n/U_{{\rm d}}(c_n)$  it follows that 
$$Q^{n}_{(-\infty,0]}(x,z) = \frac{xV_{{\rm a}}(z-1) \mathfrak{p}_{1}(-x_n)}{E|\hat Z| nc_n}\{1+o_\e(1)\}
\quad \mbox{uniformly for}  \;\; z\leq \e c_n,$$
%(note $\lim u(x) = 1/E|\hat Z|$ in  (\ref{V/a}))
 and substituting  this and  using
\beqn\label{bd_H}
      \frac1{E|\hat Z|} \sum_{z=1}^\infty V_{{\rm a}}(z-1)p(y-z) = H^{+\infty}_{(-\infty,0]}(y)  \quad (y\leq 0),
 \eeqn
(recall (\ref{H_B}) and (\ref{r_G}))  we deduce that  
$$\Sigma_{<\e c_n} 
= \frac{x_n \mathfrak{p}_{1}(-x_n)}{n}\bigg[H^{+\infty}_{(-\infty,0]}(y) - \frac1{E|\hat Z|}  \sum_{z>\e c_n} V_{{\rm a}}(z)p(y-z)  \bigg] \{1+o_\e(1)\}. $$
By Lemma \ref{L5.2}
$$\Sigma_{\geq\e c_n} \leq  C\frac{x_n }{n} \sum_{z>\e c_n} a(-z) p(y-z), $$ 
and hence the assertion  (i)  follows. It in particular follows that 
\beqn\label{bd_h1}
h^x(n, y) \leq C   n^{-1}x_n  H^{+\infty}_{(-\infty,0]}(y)\quad \mbox{for}\quad  0\leq x_n\leq M,  y\leq 0,
\eeqn 
since   $|r(n,y)| \leq C_M'H_{(-\infty,0]}^{+\infty}(y)$ in view of (\ref{bd_H}) and  $a(-z) < CV_{{\rm a}}(z)$ $(z\geq1)$. %By  (\ref{++}) and (\ref{bd_H})   this  bound is readily verified for each  $n$, thus for all $n\geq 1$.

For the proof of  (ii)  we claim that for $y\leq 0$ and $x\geq 2c_n$,
\beqn\label{upb-h2}
\sum_{z=1}^\infty Q_{(-\infty,0]}^{n}(x,z)p(y-z)\leq
%C\bigg(\frac{\mathfrak{p}_{c_\circ n}(-\frac12 x)}{n^{1-1/\a}}+
\frac{CL(x)}{x^{\a}}H_{(-\infty,0]}^{+\infty}(y) +\frac{C}{c_n}F(y-{\textstyle \frac12}x).
\eeqn
%Note that $h^x(n,y)$ is not larger than the sum on  the LHS.  
For verification of (\ref{upb-h2}) we
break the range of summation  into three parts $0<z\leq c_n$,  $c_n <z\leq x/2$ and $z>x/2$,
and denote the corresponding sums by $J_1$, $J_2$ and $J_3$, respectively.
It is immediate  from Lemma \ref{L6.2} and (\ref{bd_H})  that  $J_1 \leq CL(x)x^{-\alpha}H_{(-\infty,0]}^{+\infty}(y)$.  %By the local limit theorem estimate (\ref{eqLLT})  combined with 
% [2\mathfrak{p}_{n}(-x/2)]\wedge %
By   bound (\ref{(E)}), 
 $Q_{(-\infty, 0]}^{n}(x,z)\leq p^n(z-x)\leq  CnL(x)/c_n x^\a$   ($z<x/2$),  while  on using  (\ref{bd_H})
\beqn\label{eq4.-1} \sum_{z> c_n}p(y-z) \leq C_1\sum_{z=1}^\infty \frac{a(-z)}{n/c_n}p(y-z) \leq \frac{2C_1E|\hat Z|}{ n/c_n}\bigg[\sup_{z\geq 1}\frac{a(-z)}{V_{{\rm a}}(z)}\bigg] H_{(-\infty,0]}^{+\infty}(y).
\eeqn
Hence $J_2= C'H_{(-\infty,0]}^{+\infty}(y)L(x)/x^\a$. Finally 
$J_3 \leq C F(y-x/2)/c_n$.
 These
estimates together verify  (\ref{upb-h2}). 
As in (\ref{eq4.-1})  we derive $F(y-{\textstyle \frac12}x)\leq C H^{+\infty}_{(-\infty,0]}(y)/a^\dagger(-x)$.  %= [C_1x_nH^{+\infty}_{(-\infty,0]}(y)/x^\a]n^{1/\a}$. 
%We also have $n^{-1+1/\a}\mathfrak{p}_{c_\circ n}(x/2)\leq C_1/x^\a$ since $\mathfrak{p}_1(w) =O(e^w), w<0$. 
Hence
$$h^x(n,y) \leq C H_{(-\infty,0]}^{+\infty}(y)x_n L(x)/x^\a\quad \mbox{for}\quad x> 2c_n,\,  y\leq 0,\, n\geq 1, $$
which  combined with (\ref{bd_h1})   shows the bound in  (ii).  The proof of Lemma \ref{L6.5} is complete.
\qed
\v2

%%%%%%
%{\it Proof of Theorem \ref{thm4}.}
{\bf Proof of Theorem \ref{thm4}.} If either $x$ or $y$ remains in a bounded set,  the formula  (i) of Theorem \ref{thm4} agrees with  that of Theorem \ref{thm3}, so that
we may and do suppose both $x$ and $-y$ tend to infinity. Note that  the second ratio on the RHS of (i)  is  then asymptotically equivalent to the ratio in (ii), hence  (i) and (ii)  of Theorem \ref{thm4} is written as a single formula. Put
$$\Phi(t;\xi) = t^{-1}\xi \mathfrak{p}_t(-\xi).$$
Then  what  is to be shown may be stated as follows:  as  $n\to\infty$ and $x \vee (-y)\to\infty$
\beqn \label{pT6_1}
Q_{\{0\}}^n(x,y) \sim  \k_{\a,\ga} a(x)a(-y)c_n/n^2 + C^+ \Phi (n; x-y)  
\eeqn 
 uniformly for $-M< y_n< 0< x_n <M$, provided $0< C^+= \lim_{z\to -\infty} a(z)<\infty$. 

We follow the proof in \cite{U1dm} given to the corresponding result for case $\sigma^2<\infty$.   We employ the representation
\beqn\label{eqT5}
Q_{\{0\}}^n(x,y)=\sum_{k=1}^n \sum_{z<0}h^x(k,z)Q_{\{0\}}^{n-k}(z,y).
\eeqn
Break  the RHS into three parts by partitioning the range of the first summation as follows
\beqn\label{eq4.5}
1\leq  k< \e n; ~~ \e n \leq k \leq (1-\e)n;~~ (1-\e)n<k\leq n
\eeqn
and call the corresponding sums $I,~II $ and $I\!I\!I$, respectively. Here $\e$ is a positive constant that will be chosen small. 
The proof is divided into two cases corresponding to (i)  and (ii). Suppose (\ref{a_bdd})  to hold (so that $0<C^+<\infty$). This entails  $E[-\hat Z]<\infty$, which in fact holds under the weaker condition  $\sum_{x=0}^\infty F(-x) [a(x)+a(-x)]<\infty$  \cite[Corollary 4.1] {Unote}.   
Hence $U_{{\rm d}}(x) \sim C_1x$ and $V_{{\rm a}}(x)  
 \sim C_2 a(-x)$ ($x\to\infty$) with some  positive constants $C_1, C_2$.
\v2

{\sc Proof of {\rm (i)} (case  $x_n\wedge |y_n| \to 0$)}:
By duality one may suppose that $x_n\to 0$. From  $ES_{\sigma_{[1,\infty)}}=\infty$ and (\ref{a_bdd})  it follows  that
\beqn\label{Ha}
H_{(-\infty,0]}^x\{a\} = a(x) \quad\mbox{and}\quad C^+=  H_{(-\infty,0]}^{+\infty}\{a\} <\infty,
\eeqn
respectively (see Appendix (B) for the former equality).
From the latter bound above and Lemma \ref{L6.5} (or (\ref{bd_h1}))  one deduces, 
\beqn\label{hh}
 \sum_{k\geq \e n} \sum_{z<0} h^x(k,z)a(z)\leq    C_{\e} x_n
\eeqn
with a constant $C_\e$ depending on $\e$.
As the dual of  (\ref{iv}) of Lemma  \ref{L6.4} we have 
%$$p_{\{0\}}^n(z,y)  \leq  M_\e a(z) f_{-y}(n) \quad \mbox{and}   \;\;  p_{\{0\}}^n(z,y)  \sim a(z) f_{-y}(n) \quad (z_n\to 0, -M<y_n<0) $$
\beqn\label{pf}
Q_{\{0\}}^n(z,y) = a(z)f^{-y}(n)\{1+r_{n,z,y}\} \quad  (z<0, -Mc_n <y < 0)
\eeqn
where $r_{n,z,y}$ is uniformly bounded and 
tends to zero as $z/c_n \to 0$ and $n\to \infty$  uniformly for  $y$,
which together with   (\ref{hh}) shows 
$$II \leq C_{\e,M} x_n f^{-y}(n).$$

%==
%from the bound $H_{(-\infty,0]}^x(y) \leq CH_{(-\infty,0]}^\infty (y)$ (see (\ref{17})) 
%==

Similarly on using  (\ref{pf})  above
$$I=\sum_{1\le k< \e n}\,\sum_{z=-\infty}^{-1}h^x(k,z) a(z)f^{-y}(n-k)\{1+r_{n-k,z,y}\}.$$
For the evaluation of the last double sum we may replace $f^{-y}(n-k)$ by
 $f^{-y}(n) (1+O(\e))$, and observe that  the contribution  of $r_{n-k,z,y}$ to the sum  is negligible since  $\sum_{z> N} H^x_{(-\infty,0]}(z) a(z) \to 0$ ($N\to\infty$)  uniformly in $x$ in view of the second relation of (\ref{Ha}). By (\ref{hh}) the summation over $z$ may be extended to the whole half line $k\geq 1$. Now  applying  the first relation of (\ref{Ha})  we find
$$ I= a(x) f^{-y}(n)\{1+O(\e) +o(1)\}.$$
%and then substitute the expression of $h^x(k,z)$ given  in Lemma \ref{lem5.2}(i) and you find that 

As for $ I\!I\!I$   first observe that by (\ref{pf}) and Theorem \ref{thm3}
 $$\sum_{k=1}^{\e n} Q_{\{0\}}^k(z,y)=G_{\{0\}}(z,y)-r_n\leq C(a(z)\wedge a(y)) \quad \mbox{with}\quad 0\leq r_n\leq C_\e a(z)f^{-y}(n)n$$
 ($y, z<0$).
 If $y_n$ is bounded away from zero so that  $x/y\to 0$, then  $ I\!I\!I=O(x_n/n)=o(y_n/n)$. On the other hand, applying Lemma \ref{L6.5} we  see that if $y_n\to 0$,
 $$ I\!I\!I= x_n\mathfrak{p}_{1}(x_n) n^{-1}\sum_{z<0}H_{(-\infty,0]}^{+\infty}(z)G_{\{0\}}(z,y)(1+O(\e))+ O(x_n f^{-y}(n) ),$$
whereas  by (\ref{Ha})   we infer that  $\sum_{z\leq 0} H_{(-\infty,0]}^{+\infty}(z)G_{\{0\}}(z,y) \to  C^+$ as $y\to -\infty$ (for by subadditivity $|a(-y)-a(z-y)|\leq a(-z)$  so that the dominated convergence is applicable). Hence
% is $O(a(z))$ and tends to zero   uniformly for  $y< z<0$ as  $y\to-\infty$ and $C^+ = \sum_z H_{(-\infty,0]}^\infty(z)$, we conclude 
  $$ I\!I\!I= x_n\mathfrak{p}_{1}(x_n)n^{-1}( C^+ +o(1) +O(\e))+  O(x_n f^{-y}(n))$$
  (in case $-y_n\asymp 1$ the first term on the RHS is absorbed in the second).
Adding these expressions of $I$, $II$ and $I\!I\!I$ yields the desired formula, because of  arbitrariness of $\e$ as well as  the identity $x_n\mathfrak{p}_{1}(x_n)/n =\Phi(n;x)$. 
\v2

{\sc Proof of {\rm(ii)}  (case $x_n\wedge (-y_n) \geq 1/M$).}
 By  Lemma \ref{L6.5}(ii) and Corollary \ref{cor2}(i) it follows that   in this regime 
$$I\leq  \sum_{1\leq k <\e n}\frac {CL(x)} {c_k\,x^{\a-1}}\sum_{z<0}  H_{(-\infty,0]}^\infty(z) a(z)f^{-y}(n) \leq C' \frac{\e^{1-1/\a}}n.$$

 For evaluation of $I\!I\!I$  we change the variable $k$ into $n-k$ and  apply Lemma \ref{L6.2} to $Q^k_{\{0\}}(-y,-z)$  to see that   for any $\de>0$ 
\beq
\sum_{k=1}^{\e n} Q^k_{\{0\}}(z,y) &\leq& C \sum_{k\leq \de j(z)} \frac1{c_k} + C_\de \sum_{\de j(z) <k <\e n} |z|^{\a-1}L(y)/|y|^{\a}L(z)\\
&\leq&
C\de^{1-1/\a}a(z) + C_\de (\e n) a(z)L(y)/|y|^{\a},
\eeq
where $j(z)$ is any function such that $c_{j(z)}\sim z$, or what amounts to the same,  $j(z)/c_{j(z)}\sim a(z)$ as $z\to-\infty$ and   $C_\de$ may depend on $\de$ but $C$ does not. Then by Lemma \ref{L6.5}(ii)
$$I\!I\!I \leq C' n^{-1} \{C\de^{1-1/\a} + C_\de \e\} H_{(-\infty,0]}^{+\infty}\{a\}
 \leq C'' [C\de^{1-1/\a} + C_\de \e] /n, $$
hence for any $\e'>0$ we can choose $\e>0$ and $\de>0$ so that $I\!I\!I  \leq \e'/n$.

 By  Lemma \ref{L6.5}(i), (\ref{pf})  and (\ref{Ha})
$$II=\sum_{ \e n\leq k \le (1-\e)n}\frac{ x_k \mathfrak{p}_{1}(-x_k)}{k}\,\sum_{z=-x}^{-1}H_{(-\infty,0]}^\infty(z) Q_{\{0\}}^{n-k}(z,y)(1+o_\e(1))+ \frac{o(f^{-y}(n))}{\e^{1/\a}}.$$
Here (and  in the rest of the proof) the estimate indicated by $o_\e(1)$ may depend on $\e$ but is  uniform in the passage to the limit under  consideration once $\e$ is fixed.  
Since $-y_n$ is bounded away from zero as well as  infinity,
  we may replace  $Q_{\{0\}}^{n-k}(z,y)$ by $a(z)y_{n-k}\mathfrak{p}_1(y_{n-k})/ (n-k)$ to  see that
\[
II=\sum_{ \e n\leq k \le (1-\e)n}\frac{x_k|y_{n-k}|\mathfrak{p}_{1}(-x_k)\mathfrak{p}_{1}(y_{n-k})}{ k(n-k)}\,\sum_{z=-x}^{-1}H_{(-\infty,0]}^\infty(z)a(z)(1+o_\e(1))+\frac{o(1/n)}{\e^{1/\a}} \nonumber.
\]
On  noting  $x_k\mathfrak{p}_1(-x_k) =x_n\mathfrak{p}_{k/n}(-x_n) = \Phi(k/n;x_n)k/n$ and similarly for $y_{n-k}\mathfrak{p}_1(y_{n-k})$ 
$$\sum_{ \e n\leq k \le (1-\e)n}\frac{x_k|y_{n-k}|\mathfrak{p}_{1}(-x_k)\mathfrak{p}_{1}(y_{n-k})}{ k(n-k)}\,=\frac{1+o(1)}{n}\int_0^1 \Phi(t;x_n)\Phi(1-t; y_n)dt + O\bigg(\frac{\e}{n}\bigg).$$
Here we have used the fact that $\int_0^\e\mathfrak{p}_t(\xi) \xi dt/t= \int_{\xi/\e^{1/\a}}^\infty = O(\e/\xi^\a)$.
 Since  for $\xi>0$, $\Phi(t; \xi)dt$ is the distribution  of the  hitting-time  to zero by the process $\xi+Y_\cdot$, we have
$$\int_0^1 \Phi(t; x_n)\Phi(1-t; -y_n)dt=\Phi(1; x_n - y_n).$$
Hence
\beqn\label{II}
II= \frac1{n} C^+ \Phi(1; x_n - y_n) \{1+o(1)\}+ O\bigg(\frac{\e}{n}\bigg)+ \frac{o(1/n)}{\e^{1/\a}}.
\eeqn
  (as well as  $nI+n I\!I\!I \to 0$)  as $n\to\infty$ and $\e\to 0$ in this order.  Thus  (\ref{pT6_1}) is obtained, the first term on the RHS of it being negligible in the present regime.
    \qed
  
\v2
%Proof of Prop2.2
{\bf Proof of Proposition \ref{prop2.2}.} %We consider only the case $xy<0$, the case $xy>0$ being dealt with in the same way.   
The case $C^+=0$ is trivial.  Let $-Mc_n<y <0<x <Mc_n$.   If    $0<C^+ <\infty $,  on noting that   Theorem \ref{thm4} and Lemma \ref{L6.3} (in the dual form (\ref{P2_1})) together yield
$$
\frac{Q_{\{0\}}^{n-k}(z,y)}{Q^n_{\{0\}}(x,y)} \leq C \frac{a(z)[1 + |y|n/c_n^2]}{1 + |y|n/c_n^2 + xn/c_n^2} \leq Ca(z)  \quad (z<0, k<n/2)$$
 and  that $H_{(-\infty,0]}^x(z)\leq (E|\hat Z|) H_{(-\infty,0]}^\infty(z)$, we deduce  that   the conditional probability
 $$P[S^x_{\sigma_{(-\infty,0]}}<-R, \sigma^x_{(-\infty,0]} <n/2\,|\, \sigma_{\{0\}}>n, S^x_n=y]   =\frac{\sum_{k<n/2}\sum_{z<-R} h^x(k,z)Q^{n-k}_{\{0\}}(z,y)}{Q^n_{\{0\}}(x,y)}$$
  is at most a constant multiple of $\sum_{z<-R}H_{(-\infty,0]}^\infty(z)a(z)$ which approaches zero   as $R\to\infty$.  For the sum over $n/2\leq k \leq n$,  one  uses the bound  $\sum_{n/2\leq k\leq n}Q^{n-k}_{\{0\}}(z,y) \leq G_{\{0\}}(z,y) \leq G_{\{0\}}(z,z) \leq  Ca(z)$ as well as Lemma \ref{L6.5}(ii) to obtain 
 the same bound in a similar way. These together verify the first half of the asserted formula.

The second half   obviously follows  if $E|\hat Z|=\infty$ so that $H^{+\infty}_{(-\infty,0]}$ vanishes.  %in view of the duality as applied in  the proof of Proposition \ref{prop5.2}(i). 
Let $E|\hat Z|<\infty$.  We can then  apply  Lemma \ref{L6.5}(ii) as well as Theorem \ref{thm3} (in a dual form)  %(\ref{eqC2})) 
to see that  the contribution to the sum in (\ref{eqT5}) from   $- R \leq z < 0$ is dominated by a constant  multiple of 
\[
\sum_{-R\leq z  < 0} \sup_{k<n/2}\Big[H_{(-\infty,0]}^x(z) Q^{n-k}_{\{0\}}(z,y) + h^x(n-k,z) G_{\{0\}}(z,y)\Big]  \leq Ca(-R)\bigg[ \frac{a(-y)}{n^2/c_n} +  \frac{x\vee |y|}{nc_n}\bigg]
\]
which is negligible (as $x\vee |y|\to\infty$) as compared with the lower bound of $Q^{n}_{\{0\}}(x,y)$ given by Proposition \ref{prop6.1}, provided that  $C^+=\infty$ 
 (see Remark \ref{r_P6.1}(ii)). 
%  or,  equivalently,  $\sum_{w\geq 1} p(-w)w[a(-w)]^2 =\infty$.      
\qed

%\section{Extension to an arbitrary finite set}
\section{Extension to general finite sets}

Let $A$ be a finite non-empty subset of $\Z$. The function $u_A(x), x\in \Z$  defined  in  (\ref{def_u}) may be  given by 
\beqn\label{u/g}
u_A(x) = G_A(x,y)\1(y\notin A) +\1(x=y\in A)+ a(x-y) - E[a(S^x_{\sigma_A} -y)]
\eeqn
(whether (\ref{HA}) is assumed or not), for the RHS is independent of $y\in \Z$ (cf. \cite[Lemma 3.1]{U1dm_f}, \cite{PS}) and  the difference of the last two terms in it tends to zero as $|y|\to \infty$.
Taking an arbitrary $w_0\in A$ for $y$ it in particular follows that
\beqn\label{uA/a}
u_A(x) = a^\dagger(x-w_0) - E[a(S^x_{\sigma_A} -w_0)].
\eeqn
Hence $u_A(x)\sim a(x)$ as $x\to +\infty$    if $C^+= \lim_{x\to +\infty} a(x)=\infty$; and similarly for the case  $x\to-\infty$.
 If $C^+<\infty$, then  there exists
\beqn\label{C+}
C^+_A:=\lim_{x \to +\infty} u_A(x) = C^+ - H^{\infty}_A\{a(\cdot -w_0)\},
\eeqn
where $H^{\infty}_A(z) := \lim_{|x| \to \infty} P[S^x_{\sigma_A} =z]$ which exists for every $z\in A$ (\cite[Theorem 30.1]{S}).  (Note incidentally  (\ref{uA/a})  shows that   $\lim E[a(S^x_{\sigma_A} -w_0)]$ does not depend  on the choice of $w_0$ since  $a(x-w_0) -a(x-w)\to 0$ for any $w$.) 
The function $u_A$ is harmonic for the  walk killed on  $A$ as noted previously, and  $u_A(S^x_{n})\1( n <  \sigma^x_A)$ is accordingly a martingale for each $x\in \Z$.  It holds that
for $x> m:=\max A$
 \beqn\label{u/A}
 u_A(x) = E[u_A(S^x_{\sigma_{(-\infty,m]}}) ;\,  S^x_{\sigma_{(-\infty,m]}} \notin A] \quad\mbox{if}\quad   E[S_{\sigma_{[1,\infty)}}]=\infty
 \eeqn 
analogously
to the corresponding relation for $a(x)$ (see Appendix (B)). It follows that  $C_A^+>0$ whenever  the assumption  (\ref{HA}) is satisfied, for $C^+<\infty$ entails  $E[S_{\sigma_{[1,\infty)}}]=\infty$.

For the following two lemmas we do not 
need to use the assumption 2)  (the strong aperiodicity) and to make this clear we restate a result that  has been shown
 up to Lemma \ref{lem4.4}  as follows: {\it Uniformly for   $|x_n|<\!< \La_n(x)$, as $n\to\infty$} 
\beqn\label{000}
P[\sigma_{\{0\}}^x > n] \sim a^\dagger(x) P[\sigma_{\{0\}}^0 > n].
\eeqn
This  has been shown for $|\ga| <2-\a$ without assuming  2). Since we have  $f^0(n) \leq C c_n/n^2$ (see Remark \ref{rem2.2}) the same is true also for $\ga =\pm (2-\a)$  (see Lemma \ref{lem04}). 

According to Theorem 4a of \cite{KS} for each  $x$, as  $n\to\infty$
$$P[\sigma_A^x > n]/P[\sigma_{\{0\}}^0 > n] \, \longrightarrow\, u_A(x) $$ 
(valid for all irreducible and recurrent walks on $\Z$ with infinite variance so as to be applicable  in  the present setting),
of which we need  the  following uniform version. 
%Lem7.1
\begin{lem}\label{lem5.1} \, Uniformly for   $|x_n|<\!< \La_n(x)$, as $n\to\infty$ 
\beqn\label{eqL7.1}
P[\sigma_A^x > n] \sim u_A(x)P[\sigma_{\{0\}}^0 > n].
\eeqn 
\end{lem}
\n
\pf\, 
We adapt the proof in \cite{KS} of its Theorem 4a, which is somewhat simplified due to  the explicit asymptotic form  of $P[\sigma_{\{0\}}^0 > n]$  available for us. The proof is made by induction on the number
of points in $A$. Suppose  (\ref{eqL7.1}) holds for the sets of $N$ points and  let  $A$ consist of $N+1$ points.  
Put $A'= A\setminus \{w\}$ with a point $w\in  A$.  Then
\beqn\label{pL9.1} 
P[\sigma^x_A>n] =P[\sigma_{A'}^x>n] -\sum_{k=1}^n Q^k_{A}(x,w) P[\sigma^w_{A'}>n-k].
\eeqn
Let $x_n <\!< \La_n(x)$,
take  a constant $\e \in (0,\frac12)$ and  put $m= \lfloor \e n\rfloor$.  Then, observing 
$$Q^k_A(x,z) \leq f^x_{\{z\}}(k)\leq C a^\dagger(x)r^0_m/m\quad  \mbox{for}\quad  k\geq m, z\in A$$
(where $r^0_m = P[\sigma^0_{\{0\}}>m]$  as in the proof of Lemma \ref{lem05}), 
we have
$$\sum_{k=m}^n Q^k_{A}(x,w) P[\sigma^w_{A'}>n-k] \leq C' a^\dagger(x) \frac{r^0_m}{m}\sum_{k=m}^n r^0_{n-k}= O\big(a^\dagger(x)[r_{m}^0]^2\big).$$
As for the other sum, noting  that   
$P[\sigma^w_{A'}>n-k] \sim u_{A'}(w)r^0_{n-k}$ uniformly for $k<m$  owing to the induction hypothesis and that $\sum_{k=m}^\infty Q^k_{A}(x,w) \leq C a^\dagger(x) r^0_m$, we deduce
$$\sum_{k=1}^{m-1} Q^k_{A}(x,w) P[\sigma^w_{A'}>n-k] = [G_A(x,w)-\1(w=x)]u_{A'}(w)r^0_n\{1+o(1)\} + O\big(a^\dagger(x) [r^0_n]^2\big)$$
as $n\to \infty$ and $\e\downarrow 0$ in this order.
%, where $\de_{w,x}$ is the Kronecker's delta kernel. 
 Observe   $G_A(x,w)- \1(w=x)=P[S^x_{\sigma_A}=w]$ on one hand and 
 $$G_{A'}(x,y)- G_A(x,y) = P[S^x_{\sigma_A}=w]G_{A'}(w,y) \quad \mbox{for} \;\; y\neq w$$
on the other hand, and   then  letting     $y\to\infty$ you obtain  
$$u_{A}(x) = u_{A'}(x) - [G_A(x,w)- \1(w=x)]u_{A'}(w).$$
Now on returning to (\ref{pL9.1}) substitution of the estimates obtained above leads to 
$$P[\sigma_A^x>n] = u_A(x) r^0_n + a^\dagger(x)r^0_n\times o(1),$$ 
which shows (\ref{eqL7.1}), for $a^\dagger (x)/u_A(x)$ is bounded.
 \qed

\v2
The next lemma,  valid for all $\ga$, extends Lemma \ref{lem05} to the general case of $A$. 
%Lem7.2
\begin{lem}\label{lem5.2} \, For any interval  $I\subset \R$,  uniformly for   $|x_n|<\!< \La_n(x)$, as $n\to\infty$ $$P[S^x_n/c_n \in I\,|\, \sigma_A^x >n] \,\longrightarrow\, \int_Ih(\xi)d\xi.$$
\end{lem}
\n
\pf\, As in \cite{B1} we make the decomposition
\beq
Q^n_A(x,y) &=& p^n(y-x) -\sum_{k=1}^n \sum_{z\in A} Q_A^k(x,z) p^{n-k}(y-z)\\
&=&p^n(y-x) - \sum_{k=1}^n  f_A^x(k) p^{n-k}(y) \\
&& +\sum_{k=1}^n \sum_{z\in A} Q_A^k(x,z) \{p^{n-k}(y) - p^{n-k}(y-z )\}\\
&=& J_1(y) + J_{2}(y),  
\eeq
where $J_1(y)$ and $J_2(y)$ designates, respectively,  the difference of the first two terms  and the  double sum of the third member above.  Let  $r_n^x = P[\sigma^x_{\{0\}} >n] $ as before. We claim  that for any finite interval $I$,
$$\sum_{y: y_n\in I}  J_2(y) = o(r^x_n), $$
so that $J_2(y)$ is negligible. To this end,  
splitting the sum that defines $J_2(y)$ at $k=n/2$ and employing the bound $|p^k(y)-p^k(y-z)| \leq C_1 |z|/c_k^2$ (Lemma \ref{lem01}(i)), we see that 
$$\sum_{k=n/2}^n \sum_{z\in A}  \leq C \frac{a^\dagger(x) c_n}{n^2} \sum_{k=1}^{n/2}\sum_{z\in A}|p^k(y)-p^k(y-z)| \leq C_A \frac{a^\dagger(x)c_n}{n^2}
$$  
and
$$\sum_{k=1}^{n/2}\sum_{z\in A}  \leq  C'_A P[\sigma^x_A\leq n/2] /c^2_n\leq C_A'/c_n^2,$$
so that 
$$\sum_{y_n\in I}  J_2(y) \leq C_{ A}\bigg(\frac{ a^\dagger(x)c_n}{n^2} + \frac1{c_n^{2}}\bigg) \sum_{y_n\in I}  1 
=  C_{I,A}\bigg(\frac{c_n}{n}r^x_n + \frac1{c_n}\bigg) =  o(r^x_n),$$
showing the claim. 
Using   Lemmas \ref{lem5.1}  in place of Lemma \ref{lem02}(i)  we can follow  the proof of Lemma \ref{lem05} words for words to see 
$$\frac1{P[\sigma^x_A>n]}\sum_{y=-\infty}^\infty  J_1(y) e^{i\th y/c_n}  \to \phi_h(\th),$$
which finishes  the proof.
\qed
\v2

%Proof of Theorem {thm5}.
{\it  Proof of the extensions of Theorems \ref{thm1} to \ref{thm4}.}
With Lemmas \ref{lem5.1} and \ref{lem5.2} at hand  the same arguments that prove Theorem \ref{thm2} and Proposition \ref{prop1}
apply to the estimation of $Q_A^n(x,y)$ so as to conclude the extension of  them to a general finite set $A$, which in turn allows us  to follow the proof of 
  Theorems \ref{thm3} and \ref{thm4} to obtain the extensions of them.  \qed

%\v2
%For the extension of Theorem \ref{thm4} in which $C^+$ is replaced by $C_A^+$, we need %some additional information concerning $u_A$ and  $C_A^+$.   Suppose $C_+<\infty$ and, 
% for simplicity,  $A\subset (-\infty,0]$.
%Then for $x\geq 1$ and $y < \min A$, we have 
%$G_A(x,y) = \sum_{z\notin A, z\leq 0}H_{(-\infty,0]}^x(z)G_A(z,y)$ of which 
%by $G_A(z,y) \leq G_{\{0\}}(z,z)$ the summability is  uniform in $y$ since $C^+<\infty$. Hence   letting first $y\to-\infty$ and then  $x\to +\infty$  we see that  
%\beqn\label{uA/H}
%u_A(x)=  \sum_{z\notin A, z\leq 0}H_{(-\infty,0]}^x(z)u_A(z)\quad{and} \quad C_A^+=  \sum_{z\notin A, z\leq 0}H_{(-\infty,0]}^{+\infty}(z)u_A(z).
%\eeqn
%With these identities we can follow the proof of Theorem \ref{thm4} word for word except for trivial modifications to obtain the corresponding formula for $Q^n_A(x,y)$.

\v2 
 We conclude this section with a comment  concerning  the deduction of the results for 
 the periodic walks from those for aperiodic ones.
  Suppose the walk is not strongly aperiodic with period $\nu\geq 2$.  When $A=\{0\}$ the  problem is addressed   in Remark \ref{rem2.2}. As   for  the  general case,  recall that our proof of Theorem \ref{thm2} and Proposition \ref{prop1}  is based on the results for $A={\{0\}}$ corresponding to  Lemmas \ref{lem5.1} and   \ref{lem5.2}  that are also valid for periodic walks.  Combined with the asymptotic relation  
   $$\frac{p^{n\nu+k}(x)}{\nu} = \frac{\mathfrak{p}_1(x_{n\nu})}{c_{n\nu}} \1(x\in D_k)\{1+o(1)\} \quad(n\to\infty) \quad  \mbox{for}\quad k=0,\ldots, \nu-1$$
 where $D_k = \{ x: \exists n\geq 1, p^{n\nu+k}(x)>0\}$,  this fact allows us to deduce
   the results for  periodic walks from those which are obtained in Theorem \ref{thm2} and Proposition \ref{prop1}. For $\ga =2-\a$ we can dispose of the regime  $x_n \downarrow 0$ and  $y_n \asymp 1$ by using Lemma \ref{L7.3} below  (with $B=\{0\}$)  which applies  to periodic walks so that we also obtain the extension of Theorem \ref{thm3} to the periodic walks.
Thus,  e.g.,  if $\ga =2 - \a$,  uniformly for  $ |x|\vee |y| <Mc_n$ satisfying $p^n(y-x)>0$,
$$\frac{Q_A^n(x,y)}{\nu} = \bigg\{u_A(x)f_*(n) + \frac{U_{{\rm d}}(x_+)\mathfrak{p}_1(x_n)}{U_{{\rm d}}(c_n)n}\bigg\}u_{-A}(-y) \quad (-\La_n(-y)<\!< y_n <\!<1).
$$  
Here $f_*(n) =  \k_{\a,\ga}c_n /n^{2}$ as in Remark \ref{rem2.2}.  
Since $f_A^x(n)= \sum_{y\in A} Q_A^n(x,y) $, it follows immediately that uniformly for $|x|<Mc_n$ and  for $k=0,\ldots, \nu-1$, as $n\to\infty$ under $n \in \nu\Z$
  $$\frac{f^{x}_A(n +k)}{\nu} = 
\k^{A}_{k, x} \bigg\{u_A(x)f_*(n) +  \frac{U_{{\rm d}}(x_+)\mathfrak{p}_1(x_n)}{U_{{\rm d}}(c_n)n}\bigg\}(1+o(1)),
$$
  where $\kappa^{A}_{k,x} = \sum_{y\in A:  y-x\in D_k}  u_{-A}(-y)$.

The next lemma, virtually the same as Corollary \ref{cor5} with $A$ and $B$ interchanged,  holds  without assuming 2).  We suppose $y>0$ in it for simplicity; a similar result also  holds  for $y<0$. 
%Lem7.3
\begin{lem}\label{L7.3}\,Let $\ga> -2+\alpha$. For any non-empty subset  $B\subset A$,  as $n\wedge x\wedge y\to\infty$ under  $y <Mc_n$, $x<\!< c_n$ and $x, y, n \in \nu\Z$
\beqn\label{Q/A/u}
Q^n_{B}(x,y)- Q^n_A(x,y) = \left\{
\begin{array}{lr}
o(Q^n_{B}(x,y)) \quad  &\mbox{if}\quad C^+=\infty,\\[1mm]
H_A^\infty\{u^0_B\}f^{-y}(n)\{1+o(1)\} +O(x_n/n)\quad &\mbox{if}    \quad C^+<\infty. 
\end{array} \right.
\eeqn
where $u^0_B(x) :=u_B(x)\1(x \notin B)$.
 
 [\! In the special case $B=\{0\}$,    $u^0_B=a$ and  if $C^+<\infty$,  $H_A^\infty\{a\}= C^+ - C_A^+$.]
\end{lem}
\n
\pf\,  The LHS of (\ref{Q/A/u})   equals   
$P[\sigma_A^x \leq n <\sigma_{B}^x, S^x_n =y]$ and  is written as
$$\sum_{k=0}^n\sum_{w\in A\setminus B} Q^k_A(x,w) Q_B^{n-k}(w,y).$$
Split the outer sum at $k=\e n$. Then for each $\e>0$, the sum over $k\geq \e n$ 
  is at most
 $$\sum_{w\in A\setminus B} \sup_{k\geq \e n} Q^k_{B}(x,w)  \leq C f^x_B(n)  = o(Q^n_{B}(x,y))
 $$ 
 as $n\wedge  y\to\infty$. 
 For the other sum we substitute  from the relation $Q_{B}^{n-k}(w,y) \sim u_B(w) f_{-B}^{-y}(n-k)$ to  see that as $n\to\infty$ and $\e\downarrow 0$ in this order, 
$$\sum_{k=0}^{\e n}\sum_{w\in A\setminus B} \sim  H^x_A\{u^0_B\}f_{-B}^{-y}(n)$$
(see (\ref{H_B}) for $H_A^x$). Since  
$H^x_A\{u_B^0\} \to H_A^\infty\{u^0_B\}$  ($|x|\to\infty$) and, owing to (\ref{U/x}),  $f^x_B(n) \leq C[u_B(x)f^0(n)+ x_n/n]$, we conclude the result. \qed
\v2

%Section8
\section{Some properties of $\mathfrak{f}^{\,\xi}$ and $\mathfrak{p}^{\{0\}}_t$ }
In this section we present  properties of $\mathfrak{f}^{\,\xi}$ and $\mathfrak{p}^{\{0\}}_t$
 that are relevant to our estimate of  $Q^n_A(x,y)$. We give  proofs of them, which are devised easily from  the  known  facts as given in \cite{Bt} or \cite{Sk} except, perhaps,  for Lemmas \ref{L8.6} and
  \ref{L8.10}.    

 By specializing the series expansion  of  $\mathfrak{p}_1(x)$ as is  found in, 
 e.g.,  \cite[Lemma 17.6.1]{F} one deduces 
 \beqn\label{p0}
\mathfrak{p}_1(0) = \frac{\Ga(1/\a)}{\pi\a} \sin \frac{\pi(\a-\ga)}{2\a};
\eeqn
 if $|\ga|=2-\a$, in particular,  $\mathfrak{p}_1(0)= -1/ \Ga(-1/\a)$.

%Lem8.1
\begin{lem}\label{lem7.0} \,
 \beqn\label{eq3.16}
\mathfrak{f}^{\, \pm1}(t)= \frac{\sin (\pi/\a)}{\pi \mathfrak{p}_1(0)}\cdot\frac{\pm 1}{\a t^{1+1/\a}}\int_0^1(1-u)^{-1+1/\a}u^{-2/\a}\mathfrak{p}_1'(\mp (tu)^{-1/\a})du.
\eeqn
\end{lem}
\n
\pf\, According to   \cite[Lemma 8.13]{Bt}  
$$\int_0^t \mathfrak{f}^{\, \pm 1}(s) ds=\frac{\sin (\pi/\a)}{\pi \mathfrak{p}_1(0)}\int_0^t(t-s)^{-1+1/\a}\mathfrak{p}_s(\mp1)ds.$$
Substitution from $\mathfrak{p}_s(\mp1)= s^{- 1/\a}\mathfrak{p}_1(\mp s^{-1/\a})$ and the change  of variable $u=s/t$ transform  the integral on the RHS into
$$\int_0^1(1-u)^{-1+1/\a} u^{-1/\a} \mathfrak{p}_1(\mp (tu)^{-1/\a})du.$$
On  noting that $\int_0^1u^{-1/\a-1}|\mathfrak{p}_1'(\mp u^{-1/\a})|du = \a\int_1^\infty |\mathfrak{p}_1'(\mp x)|dx<\infty$ differentiation leads to the formula of the lemma. \qed

%Lem8.2
\begin{lem}\label{lem7.01} \, Let $\ga\neq 2-\a$. Then  $\mathfrak{f}^1(t) $ admits the following asymptotic expansion as $t\downarrow 0$\,{\rm :}
$$\mathfrak{f}^1(t) \sim \frac{\Ga(1/\a) \sin(\pi/\a)}{\pi^2\mathfrak{p}_1(0)}\sum_{k=1}^\infty (-1)^{k-1} \frac{\Ga(\a k+1)}{\Ga(k+\frac1{\a})} \sin[{\textstyle \frac12} k\pi (\a+\ga)] t^{k+1/\a}.$$
\end{lem}
\v2\n
\pf\, The  result is obtained in \cite{KKPW} (the normalization therein differs from ours only by the factor $1/\a$ to the L\'evy measure).  Here we present a different  proof. By the same manner as the asymptotic expansion of $\mathfrak{p}_1(x)$ as $|x|\to\infty$ is derived (cf. \cite{Zol})
one can show that $\mathfrak{p}_1'(x)x^2$ admits the asymptotic expansion in powers of $|x|^{-\a}$, so that the asymptotic expansion of $\mathfrak{p}_1'(x)$ is obtained by formally differentiating  the expansion  of $\mathfrak{p}_1(x)$ (given in \cite[Eq(14,34-35)]{Sk}).  This results in     
$$\mathfrak{p}'_1(-x) =  \frac{\a}{\pi}\sum_{k=1}^{n-1} (-1)^{k-1}\frac{(k\a+1)\Ga(\a k)}{\Ga(k)} \sin[{\textstyle \frac12} k\pi (\a+\ga)] x^{-\a k-2} + O(x^{-\a n-2})$$
as $x\to\infty$ for any $n\geq 2$.
After substitution  into (\ref{eq3.16}) an easy computation leads to the asymptotic expansion of the lemma.  \qed
  
\v2
We have stated the  asymptotic form of $\mathfrak{f}^1(t)$ as $t\to\infty$  in (\ref{cor0}) with two expressions of the constant factor $\kappa^{\mathfrak{f}}_{\a,\ga}$. 
By another application of (\ref{eq3.16})  we obtain it with the second expression of $\kappa^{\mathfrak{f}}_{\a,\ga}$  for $\gamma \neq 2-\alpha$ in the following lemma.

%Lem8.3 
\begin{lem}\label{lem7.1} \, 
$$
\lim_{t\to\infty} t^{2- 1/\a}\, \mathfrak{f}^1(t) =  \bigg[\frac{\sin (\pi/\a)}{\pi \mathfrak{p}_1(0)}\int_0^\infty u^{1-\a} \mathfrak{p}_1'(-u)du\bigg].
$$
and
$\int_0^\infty u^{1-\a}\mathfrak{p}_1'(-u)du>0 \; \mbox{or}  \;  =0 \;\; \mbox{according as} \;\;
\ga < 2-\a\;  \mbox{ or} \; \;  \ga =2-\a.
$
\end{lem}
\n
\pf\,  On performing  the change of variable $u=1/tx^\alpha$ (\ref{eq3.16}) becomes
$$\mathfrak{f}^1(t)= \frac{\sin (\pi/\a)}{\pi \mathfrak{p}_1(0)}\cdot \frac{1}{ t^{2-1/\a}}\int_{t^{-1/\a}}^\infty \Big(1-\frac{1}{tx^\alpha}\Big)^{-1+1/\a} x^{1-\a}\mathfrak{p}_1'(-x)dx,$$
which shows  the relation asserted by the lemma, for  the above integral  restricted to $x\in [1/t^{1/\a}, 1/(t\e)^{1/\a}]$ tends to zero for any $\e>0$, hence the integral itself
  is  asymptotically equivalent to  
  $\int_{0}^\infty  x^{1- \a}\mathfrak{p}_1'(-x)dx$ as $t\to\infty$.  The second half  follows from (\ref{cor0}).
     \qed

%Lem8.4
\begin{lem}\label{lem7.2} If  $\fa(t)$ is a continuous function on $t\geq 0$, then for  $T>0$
$$  \kappa^a_{\a,\ga,\mp} =\a\int_0^\infty \frac{\mathfrak{p}_1 (\pm u) -\mathfrak{p}_1 (0)}{u^{\a}}du \quad  \mbox{and}\quad    \lim_{y\to \pm 0}\int_0^T \frac{\mathfrak{p}_t (y) -\mathfrak{p}_t(0)}{|y|^{\a-1}}\fa(t)dt  =  \kappa^a_{\a,\ga,\mp}  \fa(0). $$
(See Lemma \ref{lem3.1}, (\ref{kapp}) for $\kappa^a_{\a,\ga,\mp}$).
\end{lem}
\n
\pf\,  Let $w_y(t) = [\mathfrak{p}_t (y) -\mathfrak{p}_t(0)]/|y|^{\a-1}$. For any $\e>0$,
$$\int_0^\e w_y(t)dt =\int_0^\e [\mathfrak{p}_1(y/t^{1/\a}) -\mathfrak{p}_1(0)] \frac{|y|dt}{|y|^{\a} t^{1/\a}} = \a\int_{|y| /\e^{1/\a}}^\infty \frac{\mathfrak{p}_1(\pm u)-\mathfrak{p}_1(0)}{|u|^\a}du,$$
where $\pm$ accords to the sign of $y$. The last member converges to a constant, say $b^\pm_{\a,\ga}$ 
and   $w_y(t) = O(|y|^{2-\a}) \to 0$  ($y\to0$) uniformly for $t>\e$ (since  $\mathfrak{p}_t'$ is  bounded), and hence the result follows.  It remains to show  $b_{\a,\ga}^{\pm}= \kappa^a_{\a,\ga,\mp}$, which is done in the following lemma. \qed

\v2
%lem8.5
\begin{lem}\label{lem7.3} %Let  $b_{\alpha,\gamma}^\pm$ be given as in Lemma \ref{lem7.2}.  
Uniformly for $x>0$, as $y\to \pm 0$
$$\mathfrak{p}_t^{\{0\}}(x,y)/|y|^{\a-1} \;\longrightarrow\,  \kappa^a_{\a,\ga,\mp} \mathfrak{f}^{\, x}(t).
$$
%in particular if  $\ga =2-\a$,  $b^-_{\a,\ga} =0$ (the trivial case) and  $b^+_{\a,\ga} = 1/\Ga(\a)$.
\end{lem}
\n
\pf\, Although the result follows from Theorems \ref{thm2} and \ref{thm3}, we use them only for the identification of the constant $b^\pm_{\a,\ga} = \a\int_0^\infty  \frac{\mathfrak{p}_1(\pm u)-\mathfrak{p}_1(0)}{|u|^\a}du$  in this proof that  is based  on the identity 
$$\mathfrak{p}_t^{\{0\}}(x,y) =\mathfrak{p}_t(y-x) -\int_0^t \mathfrak{f}^{\,x}(t-s)\mathfrak{p}_s(y)ds.$$
 On subtracting from this equality that for $y=0$ when  the LHS  vanishes, and then dividing by $|y|^{\a-1}$
$$\frac{\mathfrak{p}_t^{\{0\}}(x,y)}{|y|^{\a-1}} =\frac{\mathfrak{p}_t(y-x) - \mathfrak{p}_t(-x)}{|y|^{\a-1}} - \int_0^t \frac{\mathfrak{p}_s(y) -\mathfrak{p}_s(0)}{|y|^{\a-1}} \mathfrak{f}^{\,x}(t-s)ds.$$
 As $y\to 0$, the first term on the RHS tends to zero and Lemma \ref{lem7.2} applied to the the second term yields  the equality of the lemma.  The uniformity of the convergence is checked by
 noting that the above  integral restricted to $s>t/2$ is negligible.
  %In the case $\ga=2-\a$ (when there is no downward jump) we have $q_t^{\{0\}}(x,y) = q_t^{(-\infty,0]}(x,y) $, so  $C^-_{\a,\ga}$ must  vanish. 
By applying Theorems \ref{thm2} and \ref{thm3} with $L\equiv 1$ it follows that
$$b^\pm_{\a,\ga}  =\frac1{ \mathfrak{f}^x(1)}\lim_{y_n \to \pm 0} \frac{\mathfrak{p}_1^{\{0\}}(x,y_n)}{|y|^{\a-1}/n^{1-1/\alpha}} = \lim_{y\to\pm \infty}\frac{ a(-y)}{|y|^{\a-1}}$$
 (compare with  (\ref{eqR(d)1})), of which the last limit  is evaluated in Lemma \ref{lem3.1}(i) as asserted.
 \qed
 
\v2
Let $Q_t(y)$ denote the distribution function of a stable meander, which may be expressed as
\beqn\label{Bt}
Q_t(y) = \lim_{\e\downarrow 0} P[ Y_t \leq y\,|\, \sigma_{(-\infty,-\e]} > t ]
\eeqn
(cf. \cite[Theorem 18]{Bt}) (so that $q=Q_1'$) and  satisfies the scaling relation $Q_t(y) =Q_1(y/t^{1/\a})$. The distribution function of the meander of $-Y_t$ is denoted by $\hat Q_t(y)$.
%Lem8.6
\begin{lem}\label{L8.6}  Let  $\ga=2-\a$ and $t>0$.  Then  for $y>0$ 
\beqn\label{eq_lem6.4}
K_t(y) :=\lim_{x\downarrow 0}\mathfrak{p}_t^{\{0\}}(x,y)/x = \a \mathfrak{p}_t(0)Q'_t(y)
\eeqn
and for $x>0$
\beqn\label{eq_lem6.41}
\lim_{y\downarrow 0} \frac{\mathfrak{p}_t^{\{0\}}(x,y)}{y^{\a-1}} =\frac{\mathfrak{f}^x(t)}{\Ga(\a)} =  \frac{\hat Q'_t(x)}{\Ga(\a)\Ga(1/\a) t^{1-1/\a}}.
\eeqn
The convergences in (\ref{eq_lem6.4}) and (\ref{eq_lem6.41}) are uniform in $y>0$ and  $x>0$, respectively.
 \end{lem}
 
 Combined with  (\ref{R3f1}) and (\ref{cor0})  the equalities  above entail that if $\ga=2-\a$, 
 $$Q_t'(x) \sim [\Ga(\a+1)]^{-1} x^{\a-1}/t^{1+1/\a}   \quad\mbox{and}\quad \hat Q_t'(x) \sim [-\Ga(1/\a)/\Ga(-1/\a)] x/t^{2/\a}.$$
\n
\pf\, For the proof of  (\ref{eq_lem6.4}) first we show that for any $0<\delta<y$,
\beqn\label{6.3}
\lim_{x\downarrow 0}\frac1{x} \int_\delta^y \mathfrak{p}_t^{\{0\}}(x,z)dz =  \a \mathfrak{p}_t(0)[Q_t(y)-Q_t(\delta)].
\eeqn
For  $\ga=2-\a$,   $\sigma_{(-\infty,-\e]}^Y$ agrees with $\sigma_{\{-\e\}}^Y$ a.s. Hence for $x>0$,  the integral in (\ref{6.3}) which equals 
$  P[\delta-x <Y_t\leq y-x, \sigma_{\{-x\}}^Y>t ]$ (since $\sigma^{x+Y}_{\{0\}} = \sigma^Y_{\{-x\}}$) is expressed as
$$P[\delta-x< Y_t\leq y-x\,|\, \sigma_{(-\infty,-x]}^Y > t] P[ \sigma_{(-\infty,-x]}^Y > t]. $$
The first factor  converges to $Q_t(y)-Q_t(\delta)$ as $x\downarrow 0$. For the second one, recalling  $\mathfrak{f}^x(s)= xs^{-1}\mathfrak{p}_s(x)= xs^{-1-1/\a}\mathfrak{p}_1(xs^{-1/\a})$ and making a change of variable we have
$$P[ \sigma^Y_{(-\infty,-x]} > t] =  \int_t^\infty \mathfrak{f}^x(s) ds = \a\int_0^{x/t^{1/\a}} \mathfrak{p}_1(u)du.$$
Thus dividing by $x$ and passing to the limit conclude the required formula (\ref{6.3}) since $\mathfrak{p}_1(0)t^{-1/\a}=\mathfrak{p}_t(0)$. 
In order to conclude  (\ref{eq_lem6.4}) it suffices to show that $\lim_{x\downarrow 0}\mathfrak{p}_1^{\{0\}}(x,y)/x $ exists and the  convergence is uniform in $y$ on any compact set of $(0, \infty)$. To this end we use (\ref{Doney})  and postpone the  proof to that of Lemma \ref{L8.8}, although the proof can be done directly from the Fourier representation of $\mathfrak{p}_1^{\{0\}}(x,y)$.  

 As for (\ref{eq_lem6.41}) we make use of the duality relation and write (\ref{Bt})
as
$$\hat Q_t(x) = \lim_{\e\downarrow 0}
\frac{\int_0^{x+\e}\mathfrak{p}_t^{[0, \infty)}(-\e,-\xi)d\xi}{P[\sigma_{[\e,\infty)}>t]}
=\lim_{\e\downarrow 0}
\frac{\int_0^{x+\e}\mathfrak{p}_t^{(-\infty,0]}(\xi,\e)d\xi}{P[\sigma_{[\e,\infty)}>t]}.
$$
The first equality of (\ref{eq_lem6.41}) follows from the preceding lemma  and is written as  
$\mathfrak{p}_t^{(-\infty,0]}(\xi,\e) = \mathfrak{p}_t^{\{0\}}(\xi,\e)\sim \mathfrak{f}^{\,\xi}(t)\e^{\a-1}/\Gamma(\alpha)$ ($\xi>0$).  By $\gamma= 2-\alpha$ we have $P[Y_t>0] =1-1/\alpha$   (cf. (\ref{rho})) which entails 
 $P[\sigma^Y_{[\e,\infty)}>t]= P[\sigma^Y_{[1,\infty)}>t/\e^\a] \sim C_* (t/\e^\a)^{-1+1/\a}$ \cite[Proposition VIII.2]{Bt} and  accordingly  we obtain 
$$\hat Q_t(x) = \frac{t^{1-1/\a}}{C_*\Ga(\a)} \int_0^x \mathfrak{f}^{\,\xi}(t)d\xi.$$
We derive $C_* = 1/\Ga(\a)\Ga(1/\a)$   from  $\hat Q_t(+\infty)=1$ with the help of the next lemma (cf. Remark \ref{R9}).  
Finally  differentiation concludes  the second equality of  (\ref{eq_lem6.41}). \

The  uniformity of convergence in (\ref{eq_lem6.4}) is shown by using  (\ref{R3eq})  and the fact that 
$\sup_{0<x<1} \mathfrak{p}_t^{\{0\}}(x,y)/x \to 0$ as $y\to\infty$ (the latter  can be  shown in the same way as  Lemma \ref{L6.1})),  and similarly for the convergence in  (\ref{eq_lem6.41}).
\qed

\v2
%Lem8.7
\begin{lem}\label{lem7.5}
$$\int_{-\infty}^\infty \mathfrak{p}_t(x)|x|dx= \frac{2 t^{1/\a} }{\pi }\Ga(1-1/\a)\sin[\tst12 \pi (\a-\ga)/\a],$$
in particular  if $\ga=2-\a$, $\int_0^\infty \mathfrak{f}^x(t)dx = t^{-1}\int_0^\infty \mathfrak{p}_1(-x)xdx = t^{-1+1/\a}/\Ga(1/\a) $.
\end{lem}
\n
\pf\, Put $\chi_\la(x) =|x|e^{-\la |x|}$ ($\la>0, -\infty <x<\infty$). 
 By Parseval equality
$$\int_{-\infty}^\infty \mathfrak{p}_t(x)|x|dx = \lim_{\la\downarrow 0}\int_{-\infty}^\infty \mathfrak{p}_t(x)\chi_\la(x)dx
= \frac1{\pi}\lim_{\la\downarrow 0} \int_{-\infty}^\infty e^{- t \psi(\th)} C_\la(\th) d\th, $$
where $C_\la(\th) =\int_0^\infty \chi_\la(x) \cos \th x \, dx,$ or explicitly 
$$C_\la(\th) =  \frac{\la^2-\th^2}{(\la^2+\th^2)^2}.$$
 Observing $\int_0^\infty C_\la(\th)d\th =0$, we infer that as $\la \downarrow 0$
$$\int_{0}^\infty e^{-t \psi( \th)} C_\la(\th) d\th= \int_{0}^\infty [e^{- t \psi(\th)}-1] C_\la(\th) d\th \; \longrightarrow  \int_{0}^\infty \frac{1- \exp\{-t e^{i\ga\pi/2}\th^\a\}}{ \th^2} d\th.$$
The last integral equals $(t e^{i\ga\pi/2})^{1/\a}\Ga(1-1/\a)$ \cite[p.313 (18)]{E}, and we find the first formula of the lemma obtained. If $\ga=2-\a$,  then $\Ga(1-1/\a)\sin[\tst12 \pi (\a-\ga)/\a]= \pi/\Ga(1/\a)$, which together with $\mathfrak{f}^x(t)=xt^{-1}\mathfrak{p}_t(-x)$ and $\int_{-\infty}^\infty \mathfrak{p}_t(x)xdx=0$ shows the second 
formula.
\qed

%REM8.1
\begin{rem} \label{R9}\, %The equality in Lemma \ref{lem7.5} is a special case of the results of \cite{MP}, although the absolute moment is expressed quite differently. 
 We have used Lemma \ref{lem7.5} for identification of the constant factor in (\ref{eq_lem6.41}). Alternatively we could have applied either the exact formula for $P[\sup_{s\leq t} Y_s\in d\xi]/d\xi$ obtained in  \cite{BDP} (cf. also  \cite{DS} for related results) or a known expression of  the absolute  first moment of the  stable law (cf. \cite{P}, \cite{MP}) from which the formula of the lemma is derived by elementary algebraic manipulations.
%the formula $P[\sigma_{[0,+\infty)}>n]\sim c_1^{1/\a}n^{1-1/\a}/E|\hat Z| \Ga(1/\a)$, $c_1 = -B\Ga(1-\a)[-c_\circ (\cos \ga\pi/2)/(\cos \a\pi/2)]$ (cf. \cite{Unote}) together with Theorem \ref{thm5} and Doney's formula (\ref{Doney}). 
       %known result on the exact formula
%for $\mathfrak{h}(s,x)$ %$\int_0^t \mathfrak{h}(s,x)ds$
\end{rem}

The  results given in the next lemma are well known except for the expressions of some of constants involved. Put  $\ell^*(x) =\int_0^x P[-\hat Z>t]dt$ as in Remark \ref{R_D}.
%Lem8.8\mathfrak{f}^{\,\sgn \, x}
\begin{lem}\label{L8.8}\, Let $\ga=2-\a$.  Then
\v2
{\rm (i)} \quad $  U_{{\rm d}}(c_n)V_{{\rm a}}(c_n) \sim n/\Ga(\a), 
\;\; P[\sigma^0_{(-\infty,\,-1]}>n] U_{{\rm d}}(c_n) \to  1/\Ga(1-1/\a)  $ \;\; and
\v2
\qquad\,   $P[\sigma^0_{[0,\infty)}>n]  V_{{\rm a}}(c_n)  \to 1/\Ga(\a)\Ga(1/\a)$ \quad as  \quad $n\to\infty$; 
\v2
{\rm (ii)} \quad $ U_{{\rm d}}(x)\sim x/\ell^*(x) \quad\mbox{and}\quad   V_{{\rm a}}(x)\sim [\Ga(\a)]^{-1}x^{\a-1}\ell^*(x)/L(x)$\quad as\quad $x\to\infty$.
\end{lem}
  \n
  \pf\, It is known \cite[Eq(15) and Eq(31)]{VW} that for some positive constant $b$,
\beqn\label{b/U}
P[\sigma^0_{(-\infty,-1]}>n] \sim b/U_{{\rm d}}(c_n); \mbox{and}
\eeqn
\beqn\label{+/-}
 P[\sigma^0_{(-\infty,\,  -1]}>n]  P[\sigma^0_{[0, +\infty)}>n]  \sim \b/n \quad\mbox{with}\;\;  \b := 1/\Ga(\rho^+)\Ga(1-\rho^+),
 \eeqn
where  $\rho^+ =\frac12(1-\ga/\a)$ (cf. (\ref{rho})).
Let $\ga=2-\a$. 
Combined with (\ref{b/U})  the third and fourth cases of (\ref{Doney}) 
show that locally uniformly for  $\eta>0$,  $\mathfrak{p}_1^{\{0\}}(x_n, \eta)/x_n \sim U_{{\rm d}}(x)Q_1'(\eta)/U_{{\rm d}}(c_n)$ 
as $n\to \infty$ along with  $x_n\to \xi>0$ and then $\xi \downarrow 0$.  Thus,  locally uniformly for $y>0$,
\beqn\label{U/b1}
\mathfrak{p}_1^{\{0\}}(x,y)/x \to  bQ_1'(y) \quad  \mbox{as} \quad  x\downarrow 0.
\eeqn 
  This completes the proof of (\ref{6.3}) as is notified in the proof of Lemma \ref{L8.6} and accordingly shows 
 $b= \a\mathfrak{p}_{1}(0)$,   hence the second relation of (i) since $\a\mathfrak{p}_{1}(0) = 1/\Ga(1-1/\a)$.

 %and we verify the second relation of (i). To this end,  employing  (\ref{eq_lem6.4}) that says $\mathfrak{p}^{\{0\}}_{1}(\xi,\eta)\sim \a \mathfrak{p}_{1}(0)Q_{1}'(\eta)\xi$ ($\xi \downarrow 0, \eta>0$)  we deduce from  the third and fourth cases of (\ref{Doney}) (see also  (\ref{D1})) that as $n\to\infty$ and $\xi \downarrow 0$ in turn
% $$bU_{{\rm d}}(x)/U_{{\rm d}}(c_n) \sim  \a \mathfrak{p}_{1}(0)x_n \qquad (\xi =x_n \downarrow 0),$$ which immediately leads to  $b= \a\mathfrak{p}_{1}(0) = 1/\Ga(1-1/\a)$,  

 In a similar  way,  employing (\ref{eq_lem6.41}) 
and the second formula of (\ref{Doney}),  we obtain  the third relation of (i). 
 The first one of (i) then follows from (\ref{+/-}).

The first relation of (ii) is shown by  Rogozin\cite[Theorems 2 and 9]{R} and  then the second  follows from the first of (i). 
     \qed

In the following lemma the condition 2) is not assumed. 
%Lem8.9
\begin{lem} \label{L8.9} Let $\nu$ be the period of $S$ and $\tilde U_{{\rm d}}$  the renewal function of the descending ladder process for $\tilde S_n:=\nu^{-1}S_{\nu n}$.  Then  
$ \tilde U_{{\rm d}}(x)/  U_{{\rm d}}(x)$ tends to a positive constant as $x\to  \infty$. 
\end{lem}
\n
\pf\, By Gnedenko's theorem one deduces that $P[S_{n+1}>0] -P[S_n>0] =O(1/c_n)$, which shows that  there exists
 $\k^* := \sum_{n=1}^\infty \sum_{k=0}^{\nu-1}  \Big(\frac{P[S_{\nu n +k}>0]}{\nu n+k}-\frac{P[S_{\nu n}>0]}{\nu n}\Big)$. Observe that
 $$\sum_{n=1}^\infty \frac1{n}P[S_{\nu n}>0]t^n = -\k^* +o(1) + \sum_{n=1}^\infty \frac1{n}P[S_{n}>0]t^n  \quad (t\uparrow 1),$$
 so that $1- E[t^{\tilde \sigma_{[1,\infty)}}] = e^{-\sum n^{-1} P[S_{\nu n}>0]t^n} \sim e^{\kappa^*}(1- E[t^{\sigma_{[1,\infty)}}])$ (see \cite[Theorem XII.7.1]{F} for the  equality), where  $\tilde \sigma$ denotes the first hitting time for $\tilde S$. Hence by Tauberian theorem  $P[\tilde \sigma_{[1,\infty)}>n] \sim \kappa^* P[\sigma_{[1,\infty)}>n]$. On the other hand   
     $ nP[\sigma_{[1, \infty)}>n]/U_{{\rm d}}(c_n)$ approaches a positive constant  (\cite[Eq(31)]{VW}).  Now  we can  conclude    the assertion of the lemma.
 \qed 

%By the well-know relation between the distribution of $\sigma_{(-\infty,-1]}$ and its generating function  the assertion of the lemma follows. \qed
%for  $\sigma_{(-\infty,-1]}$ and $U_{{\rm d}}$   (cf. \cite{R}, \cite{VW},  \cite[XIV(3.4)]{F}) it then follows that $ \tilde U_{{\rm d}}(x) \sim e^{-\kappa^*}  U_{{\rm d}}(x)$. \qed

% in view of Lemma \ref{L8.8}  if $\ga=2-\a$  and of the equality  $\lim U_{{\rm d}}(x)P[\hat Z< -x] =(\sin \frac{\a+\ga}{2}\pi)/[\frac{\a+\ga}{2}\pi]$
% otherwise.  \qed

\v2
The next lemma  verifies the identity $h(\xi) = \mathfrak{f}^{-\xi}(1)/\kappa$ asserted in  Lemma \ref{lem06}.
%we consider  the probability density $h$  appearing  in (\ref{Belkin}) and (\ref{eL04}).%\beqn\label{eqL03}\hat h_0(\th) = 1- \psi(\th)\int_0^1 (1- t)^{\frac1{\a}-1}e^{ -  \psi(\th)t}dt.\eeqn

%Lem8.10 
\begin{lem}\label{L8.10} Let  $\kappa = \kappa_{\a,\ga} /(1-\frac1{\alpha})= (\sin \frac\pi{\alpha})/ 
\mathfrak{p}_1(0)\pi$. Then 
\beqn\label{eL8.7}
\frac1{\kappa}\int_{-\infty}^\infty  \mathfrak{f}^{\,-x} (1) e^{i\th x}dx = 1- \psi(\th)\int_0^1 (1- t)^{\frac1{\a}-1}e^{ -  \psi(\th)t}dt \qquad (\th\in \R).
\eeqn
\end{lem}

 Recall  that the RHS is  the characteristic function of $h$ which is denoted by $\phi_h(\th)$.

\v2\n
\pf\, From the identity  $\mathfrak{f}^{ \, x} (1)= \mathfrak{f}^{\,\sgn \, x} (|x|^{-\a})|x|^{-\a} $ and the integral representation  of $\mathfrak{f}^{\, \pm 1} (t)$ given in  (\ref{eq3.16})   we deduce  
\beqn\label{pL8.6} 
 \frac{\mathfrak{f}^{\, -x}(1)}{\kappa}=\frac{-x}{\a}\int_0^1(1-u)^{-1+1/\a}u^{-2/\a}\mathfrak{p}_1'(u^{-1/\a}x)du.
\eeqn
After a change of variable the Fourier transform of this identity  is written as
\beqn\label{pL8.61}
\int_{-\infty}^\infty \frac{ \mathfrak{f}^{\, -x}(1)}{\kappa}e^{i\th x}dx 
= \frac{1}{\a}\int_0^1(1-u)^{-1+1/\a}\om(u^{1/\a}\th) du,
\eeqn
where
$$\om(\zeta) =- \int_{-\infty}^\infty x \mathfrak{p}_1'(x)e^{i\zeta x}dx \quad (\zeta\in \R). $$
On integrating by parts 
\beq
\om(\zeta) =  \int_{-\infty}^\infty (1 + i\zeta x) \mathfrak{p}_1(x)e^{i\zeta x}dx &=& e^{-\psi(\zeta)} 
+ \zeta \frac{d}{d\zeta}e^{-\psi(\zeta)}\\
&=& e^{-\psi(\zeta)} - \alpha\psi(\zeta)e^{-\psi(\zeta)},
\eeq
and by substitution the RHS of (\ref{pL8.61}) becomes
 $$  \frac{1}{\a}\int_0^1(1-u)^{-1+1/\a}e^{-\psi(\th)u}  du + [- \psi(\th)]\int_0^1u (1-u)^{-1+1/\a}e^{-\psi(\th)u}  du.$$
Decomposing  $u=1 -(1-u)$, we can  write the second term  as 
$$\phi_h(\th)-  \bigg[1 - \psi(\th)\int_0^1 (1-u)^{1/\a}e^{-\psi(\th)u}  du\bigg],$$
of which the quantity in the  square brackets  equals the first term as is   inferred  by  integration by parts again. This results in  the identity of the lemma. \qed

\section{Appendix}

 (A) Here we state
 condition   (\ref{f_hyp}) in terms of the tails of the distribution function
$F(t) := P[X\leq t]$, and  provide the explicit expressions for the constants  relevant to the present paper and an estimate of the derivative $\phi'(\th)$.  Let $E X=0$.

The assumption (\ref{f_hyp}) on the characteristic function  $\phi(\th)$ is equivalent to the condition  
\beqn\label{A7.1}
P[X>x] \sim q^+ Bx^{-\a}L(x) \quad\mbox{and}\quad  P[X< -x] \sim  q^-Bx^{-\a}L(x)
\eeqn
as $x\to\infty$ with  a positive constant $B$ and  two non-negative constants $q^+$ and $ q^-$ such that
 $q^++q^-=1$ ($L$ is the same slowly varying function as in (\ref{f_hyp})).   The characteristic exponent and L\'evy measure  $M\{dx\}$ of the limiting stable variable $Y_1$ is given by
\[
 \psi(\th) =  |\th|^{\a}B\Ga(1-\a)\{\cos \tst12 \a\pi - i(\sgn\, \th) (q^+ - q^-)\sin \tst12 \a\pi\}
 \]
and 
 $$M\{(-x, 0]\} = Bq^- x^{2-\a}, \quad M\{(0,x]\}=  Bq^+ x^{2-\a}\quad (x>0)$$
 respectively (cf. \cite[Section XVII.3]{F}; the first equality is immediate from  the asymptotic form of $E e^{i\th X}$  ($\th \to 0$) as given by e.g.  \cite[Theorem 1.3]{B1}, \cite[Eq(6.5)]{Upot}). From the above expression of $\psi(\th)$ we read off
  \beqn\label{eq8.1}
   B\Ga(1-\a)[(\cos \tst12\a \pi)/(\cos \tst12 \ga\pi)] =1 \quad \mbox{and}\quad    \tan\tst12 \ga\pi = (q^+ - q^-)(- \tan \tst12 \a\pi )
  \eeqn
  (which reduce to $B=-/\Ga(1-\a)$ and $q^+ =1$, respectively,  if $\ga =2-\a$) and hence
 \beqn\label{eq8.2}
 \psi(\th) = (\cos \tst12 \ga\pi)  |\th|^\a\{1- i (\sgn\, \th)(q^+- q^-) \tan \tst12 \a\pi \}.
 \eeqn
According to Zolotarev \cite{Zol}  (cf. \cite[Section 8.9.2]{BGT},  \cite[Section VIII.1]{Bt}) Spitzer's  constant  $\rho^+:= \lim_{n\to \infty} n^{-1}\sum_{k=1}^n P[S_k\geq 0]$ is given by
 \beqn\label{rho}
   \rho^+= \frac12(1-\ga/\a).
   \eeqn
  On putting $\rho^-=1-\rho^+$ we deduce from (\ref{eq8.1})  that
  $$Bq^\pm = \pi^{-1} \Ga(\a)\sin (\a \rho^\pm \pi).$$
% By Theorems 1.2 and 1.3 of \cite{Unote}
%$$\lim_{n\to\infty} n^{1-1/\a}P[\sigma^0_{[0, +\infty)}>n ]   =   c_\circ^{1/\a}/\Ga(1/\a)E|\hat Z|.$$

 By the second equality of (\ref{eq8.1}) $\sin[ \frac12 \pi(\a\pm\ga)] = (\sin\frac12 \a\pi)
(\cos \frac12\ga\pi)[1\pm (q^--q^+)]$ so that the constant  $\kappa^a_{\a,\ga,\pm} $ defined in Lemma \ref{lem3.1} is expressed as follows:
\begin{eqnarray}\label{kapp}
\kappa^a_{\a,\ga,\pm} &=& -\pi^{-1}\Ga(1-\a)(\sin {\textstyle \frac12} \a\pi) 
(\cos {\textstyle \frac12}  \ga\pi)[1\pm (q^--q^+) ] \nonumber\\
&=& \frac{- \cos\,  { \textstyle \frac12} \ga\pi}{\Ga(\a) \cos {\textstyle \frac12}  \a\pi }\, q^{\mp}.
\end{eqnarray}

 From (\ref{A7.1}) it follows  that
 \beqn\label{F/psi}
 \phi'(\th)/L(1/|\th|)  \sim -\psi'(\th) =\mp \a  e^{\pm i\pi\ga/2}|\th|^{\a-1} \quad \mbox{as}\;\; \th \to  \pm0
 \eeqn
(which  is used  in the proof of Lemma \ref{lem3.1}(ii)).
Indeed,  on writing $ \phi'(\th) = i\int_{-\infty}^\infty (e^{i\th t}-1)tdF(t)$ the integration by parts  yields
 \beqn\label{integral}
 \phi'(\th) %= i\int_{-\infty}^\infty (e^{i\th t}-1)tdF(t)
 =  i\int_{-\infty}^\infty \{ e^{i\th t}-1 + i \th t e^{i\th t}\}[-F(t)\1(t<0) + (1-F(t))\1(t>0)]dt;
 \eeqn
the integral on the RHS restricted to $\{|t|<\e/|\th|\}\cup \{|t|>1/\e |\th|\}$ tends to zero as  $\th\to 0$ and $\e \downarrow 0$ in turn and, scaling by the factor $1/|\th|$,  we find that $ \phi'(\th)  \sim \pm \zeta |\th|^{\alpha-1}L(1/|\th|)$, where
\[
 \zeta=  iB \int_{-\infty}^\infty \{1- e^{\pm iu} \mp i u e^{\pm iu} \}\frac{q^-\1(u<0) - q^+\1(u>0)}{|u|^\a}du.
 \]
 Since $\zeta$ depends on the  tails of $P[\pm X>x]/L(|x|)$ only and $-\psi'(\th)$ is given by  the integral in (\ref{integral}) with $d F$ replaced by the Levy measure $M\{dx\}$,  $\pm \zeta|\th|^{\a-1}$ must be equal to $-\psi'(\th)$.

\v2\v2
(B) \,  Let $Z$ be the ascending ladder height: 
$Z=S_{\sigma_{{[1,\infty)}}}$. Recall  $H^x_B\{\fa\} =\sum_{y\in B}H_B^x(y)  \fa(y)$ for 
a function $\fa\geq 0$ on $B$ (see (\ref{H_B}) for $H_B^x$). Then    
 \beqn\label{H/a}
 H^x_{(-\infty,0]}\{a\} =a^\dagger(x) \quad (x\geq 0)\quad \mbox{if}\quad E Z =\infty. 
 \eeqn
This is shown in \cite[Corollary 1]{Uladd} for every recurrent walk irreducible on $\Z$. Under the present setting the proof is much simplified as given below.  By  standard arguments, one can see that for $x>0$,  the process $M^x(n) := a(S^x_{n\wedge \sigma_{(-\infty,0]}})$ is a martingale and $h(x):=a(x) -  H^x_{(-\infty,0]}\{a\}$ is harmonic on 
$[1,\infty)$, i.e.,  $E[h(S_1^x); S^x_1\geq 1] =h(x)$ (cf. \cite[Lemma 2.1]{U1dm}).  It follows that $H^x_{(-\infty,0]}\{a\} \leq \lim E[M^x(n)] =a(x)$ so that $h\geq 0$, entailing that  $h$ is a constant multiple of $U_{{\rm d}}(x)$ by uniqueness.   This concludes (\ref{H/a}) if $x>0$, for   we know that $a(x)/U_{{\rm d}}(x) \to 0$  because of the growth property of  $U_{{\rm d}}(x)$ as $x\to\infty$ that varies regularly with index $\a \rho^-$ (see Lemma \ref{L8.8} in case $|\ga|=2-\a$; by (\ref{rho})   $\a\rho^- = \frac12(\a +\ga)> \a-1$ if $\ga > -2+\a$). If  $x=0$,  use $\sum p(z)a(z)=1$ as well as what has just been verified.

The same argument as above applies to $u_A(x)$ for verification of (\ref{u/A}).
 
\v2\v2
%\subsection{8.3}
%\subsection{Escape  probabilities from the origin} 
(C) \; Here we give an estimate of $P[ \sigma^x_{[R,\infty)} <\sigma^x_{\{0\}}]$ as $R\to\infty$ valid uniformly  for $x<R$. The problem is much easier than the classical exit problem for 
 intervals. 
By \cite[Proposition 29.4]{S} 
$$G_{\{0\}}(x,y) = a^\dagger(x)+a^\dagger(-y)-a(x-y),$$
which entails the subadditivity $a(x+y) \leq a(x)+ a(y)$ and 
 \beqn  \label{id_h_p}   
P[\sigma^x_{\{y\}}<\sigma^x_{\{0\}}]= \frac{G_{\{0\}}(x,y)}{G_{\{0\}}(y,y)} = \frac{a^\dagger(x)+a(-y)-a(x-y)}{a(y)+a(-y)} \qquad (y\neq0, x).
 \eeqn

\v2
%Lem9.1
\begin{lem}\label{lem9.3} \, If $\ga>-2 +\a$, 
\beqn \label{R/q}
\liminf_{R\to\infty} \inf_{x\in \Z} P[\sigma^x_{\{R\}}<\sigma^x_{\{0\}}\,|\, \sigma^x_{[R,\infty)}<\sigma^x_{\{0\}}] =: q>0
\eeqn
with $q=1$ for $\ga= 2- \alpha$.  %For $\ga =\a-2$, (\ref{R/q})   remains true with $q\geq 1/2$ if $x$ is restricted to the interval $|x|<R$.
 \end{lem}
 \v2\n
 \pf\,   In view of     (\ref{id_h_p})  and the decomposition
 $$ P[\sigma^x_{\{R\}}<\sigma^x_{\{0\}}\,|\, \sigma^x_{[R,\infty)}<\sigma^x_{\{0\}}] =\sum_{z\geq R} P[S^x_{\sigma [R,\infty)} =z\,|\, \sigma^x_{[R,\infty)}<\sigma^x_{\{0\}}] P[\bar\sigma^z_{\{R\}}<\sigma^z_{\{0\}}],$$
where $ \bar\sigma^z_{\{R\}}$ is defined to be  zero if $z=R$ and agree with  $\sigma^z_{\{R\}}$ otherwise,  for the first half of the lemma    it suffices  to show that
$$\lim_{R\to\infty} \inf_{z\geq R} \frac{a(z)+a(-R)-a(z-R)}{a(R)+a(-R)} = \frac{\k^a_{\a,\ga,-}}
{\k^a_{\a,\ga,-}+ \k^a_{\a,\ga,+}};$$
the last ratio is positive if $\ga>-2+\a$ and equals unity if $\ga=2-\a$. If $\ga<2-\a$, by   Lemma \ref{lem3.1}(ii)  $a(z) -a(z-R) >0$ for $R$ large enough and the equality above follows immediately  from Lemma \ref{lem3.1}(i). The case $\ga =2-\a$   also follows from  Lemma \ref{lem3.1}(i) and (ii),  the latter    showing  $\sup_{z\geq R}|a(z)- a(z-R) |= o(R^{\a-1}/L(R)).$ \qed
\v2

%Lem9.2
\begin{lem}\label{lem9.4} \,For any $\ga$ there exists  a constant $C$ such that for $R>1$,
\[
P[ \sigma^x_{[R,\infty)} <\sigma^x_{\{0\}}] \leq  C\Big[\frac{a^\dagger(x)L(R)}{R^{\a-1}} + \frac{x_+}{R}\Big] \quad (x\leq R).
\]
\end{lem}
\v2\n
\pf\,
For $\ga>-2 +\a$, on using   Lemma \ref{lem3.1}(ii)  the result is derived   from the preceding lemma: indeed for $x<R$,
\begin{eqnarray}\label{(L9.2)} 
P[ \sigma^x_{[R,\infty)} <\sigma^x_{\{0\}}] \leq C P[\sigma^x_{\{R\}} <\sigma^x_{\{0\}}] & = &C \frac{a^\dagger(x)+a(-R)-a(x-R)}{a(R)+a(-R)} \\
&\leq& C'[a^\dagger(x)/a(-R)  + x_+ R^{-1}], \nonumber
\end{eqnarray}
where    Lemma \ref{lem9.3} is used for the first inequality  and Lemma \ref{lem3.1}(ii) is applied to  estimate   the increment of $a$  for the last inequality (as for the equality see (\ref{id_h_p})).
In case  $\ga= -2+\alpha $ see \cite[Lemma 5.5]{Upot} (use the fact that $M^x_n = a(S^x_{\sigma_{\{0\}\cup[R,\infty)}})$ is a martingale).  \qed

\v2\v2\v2

{\bf Acknowledgments.}\, 
 I wish to thank the anonymous  referees for valuable comments and suggestions. One of them  especially suggests the way to improve the original manuscript where the random walk is supposed to belong to the normal attraction of a stable law which is relaxed to the general attraction in revision. % although some results in the original one .  

\end{document}